\documentclass[11pt,a4paper]{article}
\usepackage[T1]{fontenc}
\usepackage[utf8]{inputenc}
\usepackage{authblk}
\usepackage{amssymb,amsmath}
\usepackage{amsthm,amsfonts,mathrsfs}
\usepackage{times}
\usepackage{graphicx}
\usepackage{cite}
\usepackage{xcolor}
\usepackage{hyperref}

\hypersetup{ colorlinks=true, citecolor=., linkcolor=., urlcolor=., pdfauthor={Xueru Liu, Zhenxin Liu, Wei Wang}, pdftitle={Gaussian Fluctuations of Lyapunov--Perron Reductions for Mean Field Interacting SPDEs}
}
\newtheorem{theorem}{Theorem}[section]
\newtheorem{corollary}[theorem]{Corollary}

\newtheorem{lemma}[theorem]{Lemma}
\newtheorem{proposition}[theorem]{Proposition}

\newtheorem{remark}[theorem]{Remark}
\newtheorem*{remark*}{Remark}
\newtheorem{example}[theorem]{Example}
\numberwithin{equation}{section}

\title{\textbf{Gaussian Fluctuations of Lyapunov--Perron Reductions for Mean Field Interacting SPDEs}}

\author[1]{Xueru Liu\thanks{dg21210009@smail.nju.edu.cn}}
\author[1]{Zhenxin Liu\thanks{zxliu@dlut.edu.cn}}
\author[2]{~Wei Wang\thanks{Corresponding author: wangweinju@nju.edu.cn}}

\affil[1]{School of Mathematical Sciences, Dalian University of Technology, P. R. China}
\affil[2]{School of Mathematics, Nanjing University, Nanjing 210093, P. R. China}

\date{}

\ifdefined\pdfmapline
  \ifnum\pdfoutput>0
    \pdfmapline{=eufm10 EUFM10 <eufm10.pfb}
  \fi
\fi

\begin{document}
\maketitle

\begingroup
\renewcommand{\thefootnote}{}
\footnotetext[0]{This work is supported by National Key R\&D Program of China (No. 2023YFA1009200) and NSFC (Grant 12531009).}
\endgroup

\noindent\textbf{Abstract~~~} 
We prove a functional central limit theorem for finite particle reductions of mean field interacting stochastic partial differential equations on random invariant 
manifolds. Under global smoothness and a uniform Lyapunov--Perron gap, the backward orbit and invariant graph fluctuate at order $N^{-1/2}$. Their limit 
solves a particle dependent linear fixed point equation driven by the joint Gaussian limit of the empirical fields and averaged noise. Passing from the transformed
 average and centered coordinates back to the original SPDE variables produces an additional correction at the same scale as the direct empirical fluctuation of the
  reduced drift. Although the fluctuation associated with a fixed particle is only conditionally Gaussian, averaging over an independent limiting particle yields a 
  correction jointly Gaussian with the limiting drift fluctuation and averaged noise. We obtain 
a mean-square error of order $N^{-1}$
 and a Gaussian limit in path space for reduced trajectories with a common deterministic initial value. An explicit nonlinear example computes the corrected 
 covariance and shows that omitting this correction changes the limiting trajectory variance.
\par\medskip
\noindent\textbf{Keywords~~~} interacting SPDEs, random invariant manifolds, center unstable reduction, mean field interaction, Gaussian approximation, functional CLT
\par\medskip

\section{Introduction}

We are concerned with the following system consisting of $N$-interacting stochastic partial differential equations~(SPDEs) defined on interval $(0, \pi)$
\begin{equation}\label{eq:intro-finite-particle-model}
du_i^N=\left[Au_i^N+\nabla U(u_i^N)+\frac\gamma N\sum_{j=1}^N\nabla\widetilde K(u_j^N-u_i^N)\right]dt +\sigma\,dW_i, \qquad i=1,\ldots,N,
\end{equation}
where $A=\Delta+\gamma I$  with zero Dirichlet boundary conditions on $(0, \pi)$\,, $\nabla\widetilde K(z)=z+\nabla K(z)$, $\nabla K$ is odd, and 
$\{(W_i(t))_{t\in\mathbb R}\}_{i=1}^{N}$ are independent $H=L^{2}(0, \pi)$ valued two side $Q$-Wiener processes with trace class covariance operator $Q$.  
The parameters $\gamma$ and $\sigma$  are positive constants.

System (\ref{eq:intro-finite-particle-model}) describes $N$-interacting particles which depends on spatial parameter~\cite[e.g.]{Morien1996,Criens2023}. 
 Usually $N$ is very large that one cannot treat each particle.    A natural idea is to   derive a limit system as $N\rightarrow\infty$\,, which  is  called mean field limit.  
 To derive the limit system, usually we are concerned with the empirical measure
\begin{equation*}
\mu_{N}(t, \cdot)=\frac{1}{N}\sum_{j}\delta_{u_{j}^{N}(t)}(\cdot)\,,
\end{equation*}
and  have the following  McKean--Vlasov differential equation 
\begin{equation*}
du=\left [Au+\nabla U(u)+\gamma \nabla \tilde{K}*\mu(u)\right]dt+\sigma dW(t)
\end{equation*}
with 
\begin{equation*}
\mu_N\rightarrow \mu\,,\quad N\rightarrow\infty
\end{equation*}
in some appropriate sense.
The mean-field limit has been extensively studied for both finite dimensional $N$-interacting particle system~\cite[e.g.]{McKean1966,Gartner1988,Sznitman1984,Sznitman1991} and infinite dimensional  $N$-interacting particle system~\cite[e.g.]{Morien1996,Criens2023, E, BailleulMoench2026}

However, all the results we have mentioned are on finite time interval, and few result  describes the longtime path wise behavior of the system for both  large $N$ and $N=\infty$\,.   
   An early work~\cite{Dawson} studied the number of the stationary measure for an one dimensional McKean--Vlasov equation~(that is $N=\infty$) with cubic 
   nonlinearity and linear interaction. Recent work~\cite{Bashiri1} obtained same result by  a gradient flow method and describe the  stability of the stationary 
   measures and Gess  et al.~\cite{GGH25} constructed random dynamical system from  McKean--Vlasov SDEs via rough path theory.  Very few result shows the 
   effects of number of the particles to the dynamics of the $N$-interacting particle system as $N\rightarrow\infty$\,, especially in pathwise.
   In fact the longtime  behavior of $N$-interacting particle system is different from that of the corresponding McKean--Vlasov equation, especially in the case of 
   strong interaction~\cite[e.g.]{Dawson,GomesPavliotis2018}.  
   
 For  infinite dimensional $N$-interacting particle system~(\ref{eq:intro-finite-particle-model}), recent work~\cite{Cheng} constructed  random attractor for
  McKean--Vlasov SPDEs. For finite $N\gg 1$,  Zhao and Wang~constructed an amplitude equation~\cite{ZhaoWangNonlinearity2026} and  a dimensional  $N$ 
  independent random invariant manifold~\cite{ZhaoWangSD2026} which showed that the dynamics of interacting system is determined  by a deterministic 
  equation~(see (\ref{eq:physical-limit-reduced-equation}) for system reduced on the random invariant manifold)  with interactions playing a significant role for large
   $N$.  However, the original  system is stochastic, the more detail dynamics can not be described  by  a deterministic equation. 
 
 In this paper we intend to give a more detail analysis of the random invariant manifold~(RIM) and the reduced system on RIM for large $N$.  In fact we obtain a 
 central limit theorem that describe  Gaussian fluctuation for both the  RIM  and the reduced model on the RIM for system (\ref{eq:intro-finite-particle-model}).  
 Define
\[\bar u^N=\frac1N\sum_i u_i^N,\qquad \hat u_i^N=u_i^N-\bar u^N,\qquad \bar W^N=\frac1{\sqrt N}\sum_iW_i,\qquad \hat W_i^N=W_i-\frac1{\sqrt N}\bar W^N.\]
For $y=(\bar u,h_1,\ldots,h_N)$, define
\[
\mathcal F_i^N(y)=\nabla U(\bar u+h_i)-\frac1N\sum_j\nabla U(\bar u+h_j)+\frac\gamma N\sum_j\nabla K(h_j-h_i),\qquad 1\le i\le N.
\]
The oddness of $\nabla K$ implies that the interaction term cancels in the equation for $\bar u^N$. Hence
\begin{align}
d\bar u^N&=\left[A\bar u^N+\frac1N\sum_{j=1}^N\nabla U(\bar u^N+\hat u_j^N)\right]dt+\sigma N^{-1/2}d\bar W^N, \label{eq:bar-u-mean-fluctuation-system}\\ 
d\hat u_i^N&=[(A-\gamma I)\hat u_i^N+\mathcal F_i^N(Y^N)]dt+\sigma d\hat W_i^N. \label{eq:hat-u-equation}
\end{align}
 Let $Y^N=(\bar u^N,\hat u_1^N,\ldots,\hat u_N^N)$ and
$$\mathcal F^N(y)=\left(\frac1N\sum_{j=1}^N\nabla U(\bar u+h_j),\mathcal F_1^N(y),\ldots,\mathcal F_N^N(y)\right),$$
then 
\begin{equation}\label{eq:vector-particle-system}
dY^N=[L^NY^N+\mathcal F^N(Y^N)]dt+d\mathbf M^N,  
\end{equation}
with $L^N=\operatorname{diag}(A,A-\gamma I,\ldots,A-\gamma I)$ and
\begin{equation*}
d\mathbf M^N=(\sigma N^{-1/2}d\bar W^N,\sigma d\hat W_1^N,\ldots,\sigma d\hat W_N^N).
\end{equation*}
   
Then standard Lyapunov--Perron approach construct a RIM, which is a graph of  a random function,  for the abstract system (\ref{eq:vector-particle-system}) under some Lipschitz assumptions.  Moreover one can further reduce the system (\ref{eq:vector-particle-system}) onto the RIM~\cite[e.g.]{ZhaoWangSD2026}. 
  For   $N$-interacting particle  system, however, the 
reduced vector field depends on the particle system in two different ways. The projected nonlinearity contains empirical averages, and the invariant graph on which the dynamics is reduced is 
itself constructed from the same particles. In the limiting system~($N\rightarrow\infty$) these empirical quantities are replaced by their expectations. At the central limit scale, the fluctuation of the invariant 
graph contributes at the same order as the direct empirical fluctuation of the reduced drift. Therefore, a Gaussian approximation of the reduced dynamics requires the joint fluctuation of these 
two quantities.


A central limit theorem for the reduced dynamics cannot be obtained from a central limit theorem for the empirical drift alone. The invariant graph is defined by a random Lyapunov--Perron 
equation on $(-\infty,0]$, and the graph and the empirical drift depend on the same particles. Therefore, their first order fluctuations are correlated. Jourdain and Tse \cite{JourdainTse2021} 
establish central limit theorems for nonlinear functionals of empirical measures under differentiability assumptions with respect to the measure variable. In the present setting, the relevant 
functional is given implicitly by a random fixed point equation on an infinite backward interval. We obtain its first order fluctuation directly from the Lyapunov--Perron equation, without assuming 
differentiability of the map from the empirical measure to the invariant graph. {
The main additional step is to control this implicit expansion on the infinite backward interval, jointly with the empirical drift and the averaged noise.
}
Our first result concerns the backward particle states that determine the invariant graph. Let $\mathbf U_i^{N,*}$ and $\mathbf U_i^*$ denote the finite particle and limiting backward states. 
Coupling the two systems by the same particle noises, we prove the functional central limit theorem
\[\sqrt N\bigl(\mathbf U_i^{N,*}-\mathbf U_i^*\bigr) \Longrightarrow \mathbf Y_i,\]
locally uniformly in backward time and in the center parameter, jointly with the empirical fluctuation fields and for every fixed finite collection of particle indices. The limit $\mathbf Y_i$ is the 
unique solution of a linear Lyapunov--Perron equation. 
The second result concerns the reduced dynamics in the original variables. Let $X_t^N=P_{cu}\bar u_t^N$ be the center component of a trajectory on the finite particle random graph, and let 
$X$ be its limiting reduced path. For a common deterministic initial value $X_0^N=X_0$, the limit $X$ is deterministic and, for every $T>0$,
\[\mathbb E\sup_{0\le t\le T}|X_t^N-X_t|^2\le \frac{C_T}{N}.\]
Moreover,
\[\sqrt N\,(X^N-X)\Longrightarrow Z \qquad\text{in }C([0,T],H_{cu}),\]
where $Z$ solves
\[
\begin{aligned}
 dZ_t={}&\left\{L_{cu}
 +\mathbb E\left[
 P_{cu}D\nabla U\bigl(\widehat{\mathcal U}_1(X_t,\omega)\bigr)
 D_x\widehat{\mathcal U}_1(X_t,\omega)
 \right]\right\}Z_t\,dt\\
 &+\bigl[\Gamma_t(X_t)+\mathcal R_t(X_t)\bigr]dt
 +\sigma P_{cu}\,d\bar W_t,
 \qquad Z_0=0.
\end{aligned}
\]  
The proof of the backward fluctuation theorem has two main parts. We first establish functional central limit theorems and
 uniform moment estimates for the empirical fields indexed by 
backward time, center parameters, and particle states. We then linearize the Lyapunov--Perron equation and use the uniform 
contraction gap to control the infinite backward interval. For the 
reduced dynamics, the particle dependent fluctuation is averaged with respect to an independent limiting particle.
 This identifies $\mathcal R$ and the joint law of $(\Gamma,\mathcal R,\bar W)$ appearing in \eqref{eq:physical-cu-Gaussian-limit-equation}.

Our results shows the longtime random dynamics of infinite dimensional $N$-interacting particles  is described  by a deterministic  system with a $\mathcal{O}(\frac{1}{\sqrt{N}})$ fluctuation. This description is pathwise on microscopic scale.
One related work  is the Gaussian fluctuation  of the empirical measure~\cite{WangZhaoZhu2023} as $N\rightarrow\infty$ which describes the  behavior of the system on  mesocoscopic scale. However this loses the information of   dynamics at the microscopic scale.

The rest of the paper is organized as follows. Section~\ref{sec:preliminaries} introduces the decomposition into the particle average 
and centered components, the Lyapunov--Perron equations, and the main results. 
Section~\ref{sec:gaussian-fields-linearization} proves the empirical field central limit theorem and the first order fluctuation of the backward orbit 
and invariant graph. Section~\ref{sec:gaussian-reduced-dynamics} studies the reduced dynamics in the original variables and gives the linear and nonlinear 
examples.  The appendices contain the regularity estimates for the limiting 
backward system and the moment estimates for the empirical fluctuation fields.

\section{Preliminaries and main results}
\label{sec:preliminaries}

Let $H=L^2(0,\pi)$, with inner product and norm
\[
\langle u,v\rangle_H=\int_0^\pi u(x)v(x)\,dx, \qquad |u|=\langle u,u\rangle_H^{1/2}.
\]
Let $A=\Delta+\gamma I$, with domain $H^2(0,\pi)\cap H_0^1(0,\pi)$, choose $\gamma>0$ such that  every eigenvalue of $A-\gamma I$ 
is negative. Let $\sigma\in\mathbb R$ and $\nabla\widetilde K(z)=z+\nabla K(z)$, where $\nabla K$ is odd.  The spectral decomposition,
 in the form used for invariant manifolds in \cite[Section~2]{DuanLuSchmalfuss2003}, is
\begin{equation}\label{eq:single-particle-cu-s-splitting}
H=H_{cu}\oplus H_s,\qquad A=L_{cu}\oplus A_s,\quad H_s=H_{cu}^{\perp}.
\end{equation}
where $A_{s}=A|_{H_{s}}$, $L_{cu}=A|_{H_{cu}}$ and $L_s^N=A_s\oplus(A-\gamma I)^{\oplus N}$.
Let $L^N=\operatorname{diag}(A,A-\gamma I,\ldots,A-\gamma I)$, then
\begin{equation}\label{eq:fixed-N-splitting}
H^{\oplus(N+1)}=H_{cu}\oplus\mathbb H_s^N,
\qquad \mathbb H_s^N=H_s\oplus H^{\oplus N},\qquad L^N=L_{cu}\oplus L_s^N,
\end{equation}
For $y=(\bar u,h_1,\ldots,h_N)\in H^{\oplus(N+1)}$, set
\begin{equation}\label{eq:empirical-norm}
\|y\|_N^2=|\bar u|^2+\frac1N\sum_i|h_i|^2, \qquad \mathbb E_N=\left\{y:\sum_i h_i=0\right\}.
\end{equation}
Here $\|\cdot\|_N$ is a weighted norm on $H^{\oplus(N+1)}$. Then $N^{-1}\sum_i|\bar u+h_i|^2=\|y\|_N^2$ on $\mathbb E_N$.
The spectral representation gives, with the norm $\|\cdot\|_N$ defined in \eqref{eq:empirical-norm}, for constants $\beta<0\le\alpha$,
\begin{equation}\label{eq:fixed-N-uniform-dichotomy}
\|e^{L_{cu}t}\|\le e^{\alpha t}\quad(t\le0),\qquad
\|e^{L_s^Nt}\|_{\mathcal L(\mathbb H_s^N,\|\cdot\|_N)}\le e^{\beta t}\quad(t\ge0).
\end{equation}
These bounds are independent of $N$.
%
The oddness of $\nabla K$ gives $\sum_i\mathcal F_i^N(y)=0$, and $\sum_i\hat W_i^N=0$.
Since $e^{tL^N}\mathbb E_N\subset\mathbb E_N$ for $t\ge0$, the system preserves the constraint $\sum_i\hat u_i^N=0$.

For our purpose, we assume the following conditions.

\noindent{\rm$(\mathbf{H}1)$}
The processes $\{W_i(t)\}_{i\ge1}$ are i.i.d.  Wiener processes in $H$ on a metric dynamical system 
$(\Omega,\mathcal F,\mathbb P,(\theta_t)_{t\in\mathbb R})$, with $W_i(r,\theta_t\omega)=W_i(r+t,\omega)-W_i(t,\omega)$ for $r,t\in\mathbb R$
 on a common invariant set of full probability.
Their covariance operator $Q$ is nonnegative and self adjoint and satisfies $Q^{1/2}\in\mathcal L_2(H)$.

\noindent{\rm$(\mathbf{H}2)$} The maps $\nabla U,\nabla K:H\to H$ are twice continuously Fr\'echet differentiable with bounded first and second derivatives, and $\nabla K$ is odd. Set
\[L_U=\sup_{u\in H}\|D\nabla U(u)\|_{\mathcal L(H)},\qquad L_K=\sup_{u\in H}\|D\nabla K(u)\|_{\mathcal L(H)}.\]

\noindent{\rm$(\mathbf{H}3)$} Let $\ell_{\mathrm{LP}}=8(L_U+\gamma L_K)$ and choose $\eta$ with $\beta<\eta<\min\{\alpha,0\}$ such that  
\begin{equation}\label{eq:standing-LP-gap}
q:=\ell_{\mathrm{LP}}\left(\frac1{\alpha-\eta}+\frac1{\eta-\beta}\right)<1.
\end{equation}

For $ y= (\bar u,h_1,\ldots,h_N)\in H^{\oplus(N+1)}$, $y'=(\bar u',h_1',\ldots,h_N')\in H^{\oplus(N+1)}$, by $(\mathbf{H}2)$ and Jensen inequality,
\begin{align}
\|\mathcal F^N(y)-\mathcal F^N(y')\|_N&\le\ell_{\mathrm{LP}}\|y-y'\|_N, \label{eq:FN-uniform-Lipschitz}\\ |\mathcal F_i^N(y)-\mathcal F_i^N(y')|&\le\ell_{\mathrm{LP}}\bigl(|h_i-h_i'|+\|y-y'\|_N\bigr),\qquad i\ge1. \label{eq:FNi-particle-Lipschitz}
\end{align}
Let $z^N$ be the stationary Ornstein--Uhlenbeck process
\begin{equation}\label{eq:joint-stationary-OU}
dz^N=\operatorname{diag}(\Delta,\ldots,\Delta)z^Ndt+d\mathbf M^N, \qquad z^N=(\bar z^N,\hat z_1^N,\ldots,\hat z_N^N),
\end{equation}
where
\begin{equation}\label{eq:OU-mean-fluctuation-components}
\bar z^N(t)=\frac\sigma{\sqrt N}\int_{-\infty}^t e^{\Delta(t-r)}d\bar W^N(r),\qquad \hat z_i^N(t)=\sigma\int_{-\infty}^t e^{\Delta(t-r)}dW_i(r)-\bar z^N(t).
\end{equation}
Since $\bar W^N$ is a $Q$ Wiener process with the same law for every $N$,
\[\mathbb E\sup_{|t|\le T}|\bar z^N(t)|^p\le C_{T,p}N^{-p/2},\qquad p\ge2.\]
Define
\begin{equation}\label{eq:VN-transform}
V^N(t)=Y^N(t)-z^N(\theta_t\omega).
\end{equation}
Then $ \dot V^N=L^NV^N+F^N(V^N,\theta_t\omega), $ where
\begin{equation}\label{eq:transformed-random-nonlinearity}
F^N(v,\omega)=\mathcal F^N(v+z^N(\omega))+ \bigl(L^N-\operatorname{diag}(\Delta,\ldots,\Delta)\bigr)z^N(\omega).
\end{equation}
The second term is independent of $v$. Therefore, the Lipschitz estimates in \eqref{eq:FN-uniform-Lipschitz}--\eqref{eq:FNi-particle-Lipschitz} remain valid.

Let $P_{cu}$ and $P_s^N$ be the projections associated with \eqref{eq:fixed-N-splitting}, and let $\mathbb E_{N,s}=\mathbb E_N\cap\mathbb H_s^N$.
 Then $\mathbb E_N=H_{cu}\oplus\mathbb E_{N,s}$, and
\begin{equation}\label{eq:cu-s-random-system-fixed-N}
\dot V_{cu}^N=L_{cu}V_{cu}^N+P_{cu}F^N(V^N,\theta_t\omega),\qquad \dot V_s^N=L_s^NV_s^N+P_s^NF^N(V^N,\theta_t\omega).
\end{equation}
For $\xi\in H_{cu}$, the Lyapunov--Perron method gives the integral equations with $V_{cu}^N(0)=\xi$:
\begin{equation}\label{eq:LP-orbit-fixed-N}
\begin{aligned} V_{cu}^N(t)&=e^{L_{cu}t}\xi+\int_0^t e^{L_{cu}(t-r)}P_{cu}F^N(V^N(r),\theta_r\omega)dr,\\
 V_s^N(t)&=\int_{-\infty}^t e^{L_s^N(t-r)}P_s^NF^N(V^N(r),\theta_r\omega)dr,\qquad t\le0. \end{aligned}
\end{equation}
For a Hilbert space $E$ and $\eta<0$, let
\[C_\eta^-(E)=\left\{f\in C((-\infty,0],E): \|f\|_{C_\eta^-(E)}:=\sup_{r\le0}e^{-\eta r}\|f(r)\|_E<\infty\right\}.\]
For $E=\mathbb E_N$, use the norm $\|\cdot\|_N$ and set
\[
\|V\|_{\eta,N}:=\|V\|_{C_\eta^-(\mathbb E_N)}=\sup_{r\le0}e^{-\eta r}\|V(r)\|_N.
\]
By \eqref{eq:FN-uniform-Lipschitz}, \eqref{eq:fixed-N-uniform-dichotomy}, and $(\mathbf{H}3)$,
the Lyapunov--Perron argument in \cite{DuanLuSchmalfuss2003,DuanLuSchmalfuss2004}
gives a unique solution
$V^{N,*}(\cdot,\xi,\omega)\in C_\eta^-(\mathbb E_N)$
of \eqref{eq:LP-orbit-fixed-N}.
The corresponding invariant manifold  to \eqref{eq:cu-s-random-system-fixed-N} is
\[
\mathcal M^N(\omega)=\{V^{N,*}(0,\xi,\omega):\xi\in H_{cu}\}.
\]
Moreover,
\begin{equation}\label{eq:fixed-N-graph-empirical-Lipschitz}
\|V^{N,*}(\cdot,\xi)-V^{N,*}(\cdot,\zeta)\|_{\eta,N}
\le\frac{C}{1-q}|\xi-\zeta|.
\end{equation}
Define
\begin{equation}\label{eq:original-backward-components}
\begin{aligned}
 \bigl(\bar u^{N,*},\hat u_1^{N,*},\ldots,\hat u_N^{N,*}\bigr)(r,\xi,\omega):=V^{N,*}(r,\xi,\omega)+z^N(\theta_r\omega),\qquad r\le0,
\end{aligned}
\end{equation}
and set
\[
u_i^{N,*}=\bar u^{N,*}+\hat u_i^{N,*},\qquad\mathbf U_i^{N,*}=(\bar u^{N,*},\hat u_i^{N,*}),\qquad z_{cu}^N(\omega)=P_{cu}\bar z^N(0,\omega).
\]
Let $\widetilde{\mathbb E}$ denote expectation with respect to an independent copy. With $u_i=\bar u+\hat u_i$, the limiting particle system
 associated with \eqref{eq:bar-u-mean-fluctuation-system}--\eqref{eq:hat-u-equation} is
\begin{equation}\label{eq:limiting-particle-mean-fluctuation-system}
\begin{aligned}
 d\bar u&=[A\bar u+\mathbb E\nabla U(u_i)]dt,\\ d\hat u_i&=\bigl[(A-\gamma I)\hat u_i+\nabla U(u_i) -\mathbb E\nabla U(u_i) 
+\gamma\widetilde{\mathbb E}\nabla K(\widetilde u_i-u_i)\bigr]dt+\sigma dW_i.
 \end{aligned}
\end{equation}
For fixed $\xi\in H_{cu}$, consider the backward equations associated with \eqref{eq:limiting-particle-mean-fluctuation-system}
\begin{equation}\label{eq:limiting-backward-LP-system}
\begin{aligned}
P_{cu}\bar u(t,\xi)&=e^{L_{cu}t}\xi+\int_0^t e^{L_{cu}(t-r)}P_{cu}\mathbb E\nabla U(u_1(r,\xi))dr,\\
P_s\bar u(t,\xi)&=\int_{-\infty}^t e^{A_s(t-r)}P_s\mathbb E\nabla U(u_1(r,\xi))dr,\\
\hat u_i(t,\xi)&=\sigma\int_{-\infty}^t e^{(A-\gamma I)(t-r)}dW_i(r)+\int_{-\infty}^t e^{(A-\gamma I)(t-r)}\Bigl[\nabla U(u_i)-\mathbb E\nabla U(u_1)\\
&\hspace{40mm}+\gamma\widetilde{\mathbb E}\nabla K(\widetilde u_1-u_i)\Bigr](r,\xi)dr,\qquad t\le0.
\end{aligned}
\end{equation}
Set $\mathbf H=H\times H$, with $|(v,w)|_{\mathbf H}=(|v|^2+|w|^2)^{1/2}$ and
\begin{equation*}
\mathbf H=\mathbf H_{cu}\oplus\mathbf H_s, \qquad
\mathbf H_{cu}=H_{cu}\times\{0\}, \qquad
\mathbf H_s=H_s\times H.
\end{equation*}
where $0$ denotes the zero element of $H$.
Define $\mathbf P_{cu}(\bar y,\hat y)=(P_{cu}\bar y,0), ~ \mathbf P_s(\bar y,\hat y)=(P_s\bar y,\hat y)$. 
Let $\mathbf L_{cu}=L_{cu}$ on $\mathbf H_{cu}\simeq H_{cu}$ and $\mathbf L_s=A_s\oplus(A-\gamma I)$ on $\mathbf H_s$.

Fix $p>4$. We work in $L^p(\Omega;C_\eta^-(\mathbf H))$. By \eqref{eq:fixed-N-uniform-dichotomy} and $(\mathbf{H}3)$, 
the integral operator associated with \eqref{eq:limiting-backward-LP-system} is a contraction in this weighted $L^p$ space. 
Hence \eqref{eq:limiting-backward-LP-system}, driven by $W_i$,
has a unique solution in $L^p(\Omega;C_\eta^-(\mathbf H))$, denoted by $\mathbf U_i^*(\cdot,\xi)=\bigl(\bar u^*(\cdot,\xi),\hat u_i^*(\cdot,\xi)\bigr).$
Define $u_i^*(r,\xi,\omega):=\bar u^*(r,\xi)+\hat u_i^*(r,\xi,\omega),~ r\le0.$
Taking expectations in \eqref{eq:limiting-backward-LP-system} and using the oddness of $\nabla K$, we obtain
\begin{equation}\label{eq:limiting-centered-zero-mean}
\mathbb E\hat u_i^{*}(t,\xi)=0,\qquad \bar u^{*}(t,\xi)=\mathbb E u_i^{*}(t,\xi).
\end{equation}
For deterministic $\xi$, the path $t\mapsto\bar u^{*}(t,\xi)$ is deterministic, and the limiting particles are i.i.d. For each $i$, the family indexed by $\xi$ is driven
 by the same Wiener process $W_i$. Define the fluctuation processes by
\begin{equation}\label{eq:particle-state-fluctuation}
\mathbf Y_i^N :=\sqrt N\,\bigl(\mathbf U_i^{N,*}-\mathbf U_i^{*}\bigr), \qquad y_i^N :=\sqrt N\,\bigl(u_i^{N,*}-u_i^{*}\bigr).
\end{equation}
Let $\mathbf Y_i$ denote the unique solution of \eqref{eq:gaussian-linear-fluctuation-equation} with $\mathbf U=\mathbf U_i^{*}$.

Zhao and Wang \cite{ZhaoWangSD2026} construct a finite dimensional random invariant manifold for the interacting SPDE system by decomposing it into its 
empirical average and deviations from this average. They also obtain a mean field approximation of the reduced dynamics for sufficiently large interaction strength.
The first theorem establishes the $N^{-1/2}$ fluctuation limit for the backward trajectories $\mathbf U_i^{N,*}$  which describes the fluctuation of the  
 random invariant manifold $\mathcal M^N(\omega)$.

Throughout, $\Longrightarrow$ denotes convergence in distribution in the stated space with its  topology. For probability laws, it denotes weak convergence.
For $R>0$, define $B_R^{cu}:=\{\xi\in H_{cu}:|\xi|\le R\}.$

\begin{theorem}\label{thm:intro-LP-fluctuation}
Assume $(\mathbf{H}1)$--$(\mathbf{H}3)$. For every fixed $m\ge1$ and $T,R>0$,
\begin{equation}\label{eq:intro-joint-particle-LP-fluctuation}
(\mathbf Y_i^N)_{1\le i\le m}\Longrightarrow(\mathbf Y_i)_{1\le i\le m}\quad\text{in }C\bigl(B_R^{cu},C([-T,0],\mathbf H)\bigr)^m.
\end{equation}
For each fixed $\xi\in H_{cu}$ and $1\le i\le m$, $\mathbf Y_i(\cdot,\xi)$ is Gaussian conditionally on the backward trajectory $\mathbf U_i^{*}(\cdot,\xi)$.
\end{theorem}

Next main result shows the  fluctuation of the solution on $\mathcal M^N(\omega)$.
\begin{theorem}\label{thm:physical-cu-gaussian}
Assume $(\mathbf{H}1)$--$(\mathbf{H}3)$. Let $X^N$ and $X$ be the finite particle and limiting reduced processes with values in $H_{cu}$ 
satisfying \eqref{eq:finite-N-physical-reduced-equation} and \eqref{eq:physical-limit-reduced-equation}, respectively. Suppose that $X_0\in H_{cu}$ is 
deterministic and $X_0^N=X_0$ for $N\ge2$. Then $X$ is deterministic. Let $\Gamma^N$ and $\mathcal R^N$ be defined by \eqref{eq:Gamma-field-definition} 
and \eqref{eq:averaged-particle-state-correction-N}, respectively. Then, for every $T,R>0$,
\begin{equation}\label{eq:intro-joint-direct-particle-state-fluctuation-limit}
(\Gamma^N,\mathcal R^N,\bar W^N) \Longrightarrow (\Gamma,\mathcal R,\bar W)
\end{equation}
in $ C([0,T]\times B_R^{cu},H_{cu})^2\times C([0,T],H), $ where $\Gamma$ and $\mathcal R$ are given by \eqref{eq:Gamma-limit-definition} 
and \eqref{eq:graph-average-fluctuation-formula}, respectively. The limiting fields $(\Gamma,\mathcal R,\bar W)$ are jointly Gaussian. 
For every $T>0$, there exists $C_T>0$ such that
\begin{equation}\label{eq:intro-reduced-mean-square-rate}
\mathbb E\sup_{0\le t\le T}|X_t^N-X_t|^2 \le \frac{C_T}{N}.
\end{equation}
Moreover,
\begin{equation}\label{eq:physical-cu-Gaussian-convergence}
Z^N:=\sqrt N\,(X^N-X) \Longrightarrow Z \qquad\text{in }C([0,T],H_{cu}),
\end{equation}
where $Z$ is the unique solution of
\begin{equation}\label{eq:physical-cu-Gaussian-limit-equation}
\begin{aligned} dZ_t={}&\left\{L_{cu} +\mathbb E\left[ P_{cu}D\nabla U\bigl(u_1^{*}(0,X_t,\omega)\bigr) D_xu_1^{*}(0,X_t,\omega) \right]\right\}Z_t\,dt\\ 
&+\bigl[\Gamma_t(X_t)+\mathcal R_t(X_t)\bigr]dt +\sigma P_{cu}\,d\bar W_t, \qquad Z(0)=0. \end{aligned}
\end{equation}
In particular, $Z$ is a Gaussian process.
\end{theorem}

\begin{remark}
For fixed $(t,x)$, define the empirical measures on $H$ by
\[
\begin{aligned}
\mu_{t,x}^N =\frac1N\sum_{i=1}^N\delta_{u_i^{N,*}(0,x-z_{cu}^N(\theta_t\omega),\theta_t\omega)},\qquad
\widetilde\mu_{t,x}^N=\frac1N\sum_{i=1}^N\delta_{u_i^*(0,x,\theta_t\omega)},\qquad
\mu_x =\operatorname{Law}(u_1^*(0,x,\omega)).
\end{aligned}
\]
Since $\theta_t$ preserves $\mathbb P$, $\operatorname{Law}(u_1^*(0,x,\theta_t\omega))=\mu_x.$
By \eqref{eq:Gamma-field-definition} and
\eqref{eq:averaged-particle-state-correction-N},
we have the  decomposition
\begin{equation}
\label{eq:empirical-drift-fluctuation-decomposition}
\begin{aligned}
\sqrt N\langle\mu_{t,x}^N-\mu_x,P_{cu}\nabla U\rangle
&=\sqrt N\langle\widetilde\mu_{t,x}^N-\mu_x,P_{cu}\nabla U\rangle+\sqrt N\langle\mu_{t,x}^N-\widetilde\mu_{t,x}^N, P_{cu}\nabla U\rangle\\
&=\Gamma_t^N(x)+\mathcal R_t^N(x).
\end{aligned}
\end{equation}
The first term is the empirical fluctuation of the independent limiting particles. The second term is the difference between
the finite particle and limiting reduced nonlinear drifts at the same parameter $x$, multiplied by $\sqrt N$.
Its first order expansion and limit are given by \eqref{eq:averaged-particle-state-linearization} and \eqref{eq:graph-average-fluctuation-formula} respectively.
Therefore, Theorem~\ref{thm:physical-cu-gaussian} gives
\[
\bigl((t,x)\mapsto\sqrt N\langle\mu_{t,x}^N-\mu_x,P_{cu}\nabla U\rangle\bigr) \Longrightarrow \Gamma+\mathcal R
\]
in $C([0,T]\times B_R^{cu},H_{cu})$ for every $T,R>0$.

We next compare this decomposition with the linearized McKean--Vlasov equation for empirical measure fluctuations at the formal level.
For related distribution valued limits in interacting diffusion systems, see Wang, Zhao and Zhu \cite[Theorem~1.4 and Proposition~1.5]{WangZhaoZhu2023}.
Along the reduced trajectories, set
\[
\mathcal E_t^N=\sqrt N\bigl(\mu_{t,X_t^N}^N-\mu_{X_t}\bigr),\qquad \widetilde{\mathcal E}_t^N=\sqrt N\bigl(\widetilde\mu_{t,X_t}^N-\mu_{X_t}\bigr).
\]
These are the empirical fluctuations of the interacting particles and the independent limiting particles, respectively.
For a cylindrical function $\psi(v)=f(\langle v,h_1\rangle_H,\ldots,\langle v,h_m\rangle_H),~f\in C_b^2(\mathbb R^m),~ h_1,\ldots,h_m\in D(A^*),$
define
\[
\begin{aligned}
\mathcal L_\mu\psi(v)=\langle v,A^*D\psi(v)\rangle_H +\left\langle\nabla U(v)+\gamma\int_H\nabla\widetilde K(w-v)\,\mu(dw),D\psi(v)\right\rangle_H\frac{\sigma^2}{2}\operatorname{Tr}\bigl(QD^2\psi(v)\bigr).
\end{aligned}
\]
By uniqueness in \eqref{eq:limiting-backward-LP-system}, the limiting particle trajectories on the invariant manifold are $u_i^*(0,X_t,\theta_t\omega)$.
Therefore their common law $\mu_{X_t}$  satisfies
\[
\frac{d}{dt}\langle\mu_{X_t},\psi\rangle=\langle\mu_{X_t},\mathcal L_{\mu_{X_t}}\psi\rangle.
\]
For the following comparison, assume  that a joint empirical measure central limit theorem holds
\[
(\mathcal E^N,\widetilde{\mathcal E}^N,\bar W^N) \Longrightarrow (\mathcal E,\widetilde{\mathcal E},\bar W),
\]
jointly with the limits in Theorem~\ref{thm:physical-cu-gaussian}.  For every cylindrical test function $\psi$ specified above, and also
for $\psi(v)=\langle v,h\rangle_H$ with $h\in H_{cu}$, assume that
the processes
$$
\langle\mathcal E_\cdot^N,\psi\rangle, \quad\langle\widetilde{\mathcal E}_\cdot^N,\psi\rangle,\quad \int_0^\cdot\langle\mathcal E_s^N,\mathcal L_{\mu_{X_s}}\psi\rangle\,ds, \quad
\int_0^\cdot\langle\widetilde{\mathcal E}_s^N,\mathcal L_{\mu_{X_s}}\psi\rangle\,ds,$$
 and $\int_0^\cdot
\left\langle\mathcal E_s^N,
\int_H
\left\langle\nabla\widetilde K(\,\cdot-v),D\psi(v)\right\rangle_H
\,\mu_{X_s}(dv)\right\rangle\,ds
$
converge jointly in distribution, each in $C([0,T])$, to the corresponding expressions with $\mathcal E^N$ and $\widetilde{\mathcal E}^N$ replaced by $\mathcal E$ and
$\widetilde{\mathcal E}$. Assume the same convergence of the pairing processes for $\psi(v)=\langle\nabla U(v),h\rangle_H$, $h\in H_{cu}$.
For every finite collection of these test functions, all these convergences are joint with the preceding convergence.
 The initial joint limit is assumed to be Gaussian, with the joint law determined by the finite particle and limiting invariant manifold initial data.
For every cylindrical test function $\psi$ given above, assume also that
\[
\frac{1}{\sqrt N}\sup_{0\le t\le T}\left|\int_0^t\int_{H\times H}\left\langle\nabla\widetilde K(w-v),D\psi(v)\right\rangle_H\,\mathcal E_s^N(dv)\mathcal E_s^N(dw)\,ds \right|
\xrightarrow{\mathbb P}0.
\]
Then the limiting fluctuations satisfy
\[
\begin{aligned}
d\langle\mathcal E_t,\psi\rangle
={}&\langle\mathcal E_t,\mathcal L_{\mu_{X_t}}\psi\rangle\,dt+\gamma\left\langle\mathcal E_t,\,\int_H \left\langle\nabla\widetilde K(\,\cdot-v),D\psi(v)
\right\rangle_H\,\mu_{X_t}(dv)\right\rangle dt+dM_t(\psi),\\
d\langle\widetilde{\mathcal E}_t,\psi\rangle={}&\langle\widetilde{\mathcal E}_t,\mathcal L_{\mu_{X_t}}\psi\rangle\,dt+dM_t(\psi).
\end{aligned}
\]
By \eqref{eq:physical-particle-state-fixed-center-N-rate}, \eqref{eq:finite-physical-particle-state-empirical-Lipschitz}
and \eqref{eq:reduced-mean-square-rate}, the Burkholder--Davis--Gundy inequality implies that,
for each test function $\psi$ above, the difference between the two empirical martingales
converges to zero in probability in $C([0,T])$.
Thus both equations have the same Gaussian martingale $M$, with $M_0=0$. Its covariance and its joint covariance with $\bar W$ are
\[
\begin{aligned}
\mathbb E[M_t(\psi)M_s(\varphi)]&=\sigma^2\int_0^{t\wedge s}\int_H\left\langle Q^{1/2}D\psi(v),Q^{1/2}D\varphi(v)\right\rangle_H\,\mu_{X_r}(dv)\,dr,\\\mathbb E[M_t(\psi)\langle\bar W_s,h\rangle_H]
&=\sigma\int_0^{t\wedge s}\int_H\langle D\psi(v),Qh\rangle_H\,\mu_{X_r}(dv)\,dr, \qquad h\in H.
\end{aligned}
\]
In particular, for every $h\in H_{cu}$, $M_t\bigl(v\mapsto\langle v,h\rangle_H\bigr)=\sigma\langle\bar W_t,h\rangle_H.$
The decomposition \eqref{eq:empirical-drift-fluctuation-decomposition} and the first order expansion in $x$ give
\[
\begin{aligned}
Z_t&=\langle\mathcal E_t,P_{cu}\mathrm{id}_H\rangle,\\
\Gamma_t(X_t)&=\langle\widetilde{\mathcal E}_t,P_{cu}\nabla U\rangle,\\
\mathcal R_t(X_t)&=\langle\mathcal E_t-\widetilde{\mathcal E}_t, P_{cu}\nabla U\rangle-\mathbb E\left[ P_{cu}D\nabla U\bigl(u_1^*(0,X_t,\omega)\bigr) D_xu_1^*(0,X_t,\omega) \right]Z_t.
\end{aligned}
\]
The last term accounts for the difference between $X_t^N$ and $X_t$, since $\mathcal R_t^N(x)$ compares the two states
at the same parameter $x$. Taking $\psi(v)=\langle v,h\rangle_H$ with $h\in H_{cu}$,
the two interaction contributions in the equation for $\mathcal E$ cancel by the oddness of $\nabla\widetilde K$.
Therefore,
\[
dZ_t=L_{cu}Z_t\,dt+\langle\mathcal E_t,P_{cu}\nabla U\rangle\,dt+\sigma P_{cu}\,d\bar W_t.
\]
Substituting the preceding identities yields \eqref{eq:physical-cu-Gaussian-limit-equation}.
The initial joint fluctuation must be the one determined by the finite and limiting invariant manifolds.
In particular, $X_0^N=X_0$ implies $\langle\mathcal E_0,P_{cu}\mathrm{id}_H\rangle=0$. With this initial joint law and the noise coupling above,
the empirical measure formulation and Theorem~\ref{thm:physical-cu-gaussian}
yield the same equation for the projected fluctuation $Z$.
\end{remark}

\begin{remark}
Since Theorem~\ref{thm:intro-LP-fluctuation} gives the convergence
in \eqref{eq:intro-joint-particle-LP-fluctuation} for every $R>0$,
and $H_{cu}$ is finite dimensional, we also have, for every $T>0$,
\[
(\mathbf Y_i^N)_{1\le i\le m}
\Longrightarrow
(\mathbf Y_i)_{1\le i\le m}
\quad\text{in }
C_{\mathrm{loc}}(H_{cu},C([-T,0],\mathbf H))^m.
\]
Here $C_{\mathrm{loc}}(H_{cu},C([-T,0],\mathbf H))$
is equipped with the topology of uniform convergence
on compact subsets of $H_{cu}$. By \cite[Proposition~16.6]{Kallenberg2002}. the convergence of $(\Gamma^N,\mathcal R^N,\bar W^N)$ in 
Theorem~\ref{thm:physical-cu-gaussian} also holds in $C_{\mathrm{loc}}(H_{cu},C([0,T],H_{cu}))^2\times C([0,T],H).$

\end{remark}

\section{Gaussian fields and Lyapunov--Perron linearization}\label{sec:gaussian-fields-linearization}
\subsection{Regularity of the backward solutions}

In this subsection, fix $ p>4$, choose $\frac14<\vartheta<\frac12$ such that $p\vartheta>2$ and choose
\begin{equation}\label{eq:rho-gaussian}
\frac1p<\rho<\vartheta-\frac1p.
\end{equation}
Let $V_\vartheta=D((I-\Delta)^\vartheta),~ \|u\|_{V_\vartheta}=\|(I-\Delta)^\vartheta u\|_H.$
 Since $H_{cu}$ is finite dimensional and spanned by eigenfunctions of $A$, $H_{cu}\subset V_\vartheta$. 
 If $S_s(t)$ denotes either $e^{tA_s}$ or $e^{t(A-\gamma I)}$, then
\begin{align}
\|(I-\Delta)^\vartheta S_s(t)u\|_H
&\le C_\vartheta t^{-\vartheta}e^{-ct}|u|,
\qquad 0<t\le1,
\label{eq:positive-fractional-semigroup-gaussian}\\
\|(S_s(h)-I)u\|_H &\le C_\vartheta h^\vartheta\|u\|_{V_\vartheta},
\qquad 0<h\le1.
\label{eq:semigroup-increment-gaussian}
\end{align}
For $S_s(t)=e^{tA_s}$, the first estimate holds for $u\in H_s$ and the second for $u\in V_\vartheta\cap H_s$.
For $S_s(t)=e^{t(A-\gamma I)}$, they hold for $u\in H$ and $u\in V_\vartheta$, respectively.

For $M>0$, set
\[K_M:=\{z\in V_\vartheta:\|z\|_{V_\vartheta}\le M\}.\]
For a Banach space $E$, a compact interval $I\subset\mathbb R$, and $a\in(0,1)$, define
$$
 [f]_{C^a(I,E)}=\sup_{\substack{s,t\in I\\s\ne t}} \frac{\|f(t)-f(s)\|_E}{|t-s|^a},\qquad \|f\|_{C^a(I,E)}= \sup_{t\in I}\|f(t)\|_E+[f]_{C^a(I,E)}.
$$

The following Lemmas give some regularity estimates of the solutions. Their proofs are given in Appendix~\ref{app:backward-regularity}.

\begin{lemma}\label{lem:backward-regularity-gaussian}
Assume $(\mathbf{H}1)$--$(\mathbf{H}3)$,  fix $p>4$ and  $1/p<\rho<\vartheta-1/p$.  For every $T,R>0$, there exists $C_{T,R,p}>0$ such that
\begin{equation}\label{eq:backward-Vtheta-moment}
\mathbb E \left[\sup_{|\xi|\le R}\sup_{-T\le r\le0}\|u_1^{*}(r,\xi)\|_{V_\vartheta}^p\right] \le C_{T,R,p},
\end{equation}
\begin{equation}\label{eq:backward-time-holder}
\mathbb E\left[\sup_{|\xi|\le R}[u_1^{*}(\cdot,\xi)]_{C^\rho([-T,0],H)}^p\right]\le C_{T,R,p},
\end{equation}
and
\begin{equation}\label{eq:backward-xi-Lipschitz}
\sup_{-T\le r\le0}\left|u_1^{*}(r,\xi)-u_1^{*}(r,\zeta)\right|\le C_{T,R}|\xi-\zeta|
\end{equation}
almost surely for $|\xi|,|\zeta|\le R$. In addition,
\begin{equation}\label{eq:backward-global-weighted-moment}
\mathbb E\left[\sup_{|\xi|\le R}\|u_1^{*}(\cdot,\xi)\|_{C_\eta^-}^{p}\right]<\infty,
\end{equation}
and
\begin{equation}\label{eq:backward-xi-weighted-Lipschitz}
\sup_{r\le0}e^{-\eta r}\left|u_1^{*}(r,\xi)-u_1^{*}(r,\zeta)\right| \le C_R|\xi-\zeta|.
\end{equation}
Estimates \eqref{eq:backward-global-weighted-moment} and \eqref{eq:backward-xi-weighted-Lipschitz} also hold with $u_1^{*}$ replaced by 
$\mathbf U_1^{*}$.

The map from $H_{cu}$ to $C([-T,0],H)$ defined by $\xi\mapsto u_1^{*}(\cdot,\xi)|_{[-T,0]}$ is continuously Fr\'echet differentiable. Moreover,
\begin{equation}\label{eq:backward-xi-derivative-Vtheta-moment}
\mathbb E \left[\sup_{|\xi|\le R}\sup_{-T\le r\le0}\|D_\xi u_1^{*}(r,\xi)\|_{\mathcal L(H_{cu},V_\vartheta)}^p\right]\le C_{T,R,p},
\end{equation}
\begin{equation}\label{eq:backward-xi-derivative-time-holder}
\mathbb E \left[\sup_{|\xi|\le R} [D_\xi u_1^{*}(\cdot,\xi)]_{C^\rho([-T,0],\mathcal L(H_{cu},H))}^p \right] \le C_{T,R,p}.
\end{equation}
\end{lemma}

\begin{lemma}\label{lem:backward-theta-regularity-gaussian}
Assume the hypotheses of Lemma~\ref{lem:backward-regularity-gaussian}. Let $T_0,T,R>0$, then
\begin{equation}\label{eq:backward-theta-holder}
\mathbb E\left[\sup_{|\xi|\le R}\sup_{-T\le r\le0}\bigl[t\mapsto u_1^{*}(r,\xi,\theta_t\omega)\bigr]_{C^{\rho_0}([0,T_0],H)}^p\right] <\infty,
\end{equation}
and
\begin{equation}\label{eq:backward-theta-weighted-holder}
\mathbb E\left[\sup_{|\xi|\le R}\bigl[t\mapsto u_1^{*}(\cdot,\xi,\theta_t\omega)\bigr]_{C^{\rho_0}([0,T_0],C_\eta^-)}^p\right]<\infty.
\end{equation}
Moreover,
\begin{equation}\label{eq:backward-joint-theta-regularity}
\mathbb E\sup_{\substack{0\le t\le T_0\\|\xi|\le R}} \left\{ [r\mapsto u_1^{*}(r,\xi,\theta_t\omega)]_{C^\rho([-T,0],H)}^p 
+\sup_{-T\le r\le0}\|u_1^{*}(r,\xi,\theta_t\omega)\|_{V_\vartheta}^p \right\}<\infty.
\end{equation}
There exist i.i.d. nonnegative random variables $\mathcal L_i=\mathcal L_i(T_0,T,R)$ with $\mathbb E\mathcal L_i^p<\infty$ such that,
 for $0\le t,v\le T_0$, $-T\le r,s\le0$, and $|\xi|,|\zeta|\le R$,
\begin{equation}\label{eq:backward-joint-increment-bound}
|u_i^{*}(r,\xi,\theta_t\omega)-u_i^{*}(s,\zeta,\theta_v\omega)|
\le\mathcal L_i\bigl(|t-v|^{\rho_0}+|r-s|^\rho+|\xi-\zeta|\bigr).
\end{equation}
Estimates \eqref{eq:backward-theta-holder}--\eqref{eq:backward-joint-theta-regularity} also hold with $u_1^{*}$ replaced 
by $\mathbf U_1^{*}$, and \eqref{eq:backward-joint-increment-bound} holds with $u_i^{*}$ replaced by $\mathbf U_i^{*}$.
 If $\nu,a$ are chosen as in \eqref{eq:neighboring-weight-choice}, then
\begin{equation}\label{eq:derivative-theta-reduced-holder}
\mathbb E\left[\sup_{|\xi|\le R} [t\mapsto D_\xi\mathbf U_1^{*}(\cdot,\xi,\theta_t\omega)]_{C^{a\rho_0}([0,T_0],\mathcal L(H_{cu},C_\eta^-))}^{p}\right]<\infty.
\end{equation}
Hence $(t,\xi)\mapsto D_\xi\mathbf U_1^{*}(\cdot,\xi,\theta_t\omega)$ is almost surely jointly continuous with values in $\mathcal L(H_{cu},C_\eta^-)$.
\end{lemma}


\subsection{A functional central limit theorem}

We introduce the empirical fluctuation fields associated with the
nonlinear drift terms and prove their joint convergence with the
averaged Wiener process. The limiting fields will be used in the
linear fluctuation equation.
 
Define
\begin{equation}\label{eq:empirical-GU-definition}
\mathbb G_U^N(r,\xi):=\frac1{\sqrt N}\sum_{j=1}^N\left[\nabla U(u_j^{*}(r,\xi))-\mathbb E\nabla U(u_1^{*}(r,\xi))\right].
\end{equation}
For $z\in H$, let $\widetilde u_1^{*}$ be an independent copy of $u_1^{*}$ and define
\begin{equation}\label{eq:empirical-GK-definition}
\mathbb G_K^N(r,\xi,z):=\frac1{\sqrt N}\sum_{j=1}^N\left[\nabla K(u_j^{*}(r,\xi)-z)
-\widetilde{\mathbb E}\nabla K(\widetilde u_1^{*}(r,\xi)-z)\right].
\end{equation}
Comparing the drift terms in \eqref{eq:bar-u-mean-fluctuation-system}--\eqref{eq:hat-u-equation}
with those in \eqref{eq:limiting-particle-mean-fluctuation-system},
we obtain the following decompositions:
\begin{align}
&\frac1N\sum_{j=1}^N\nabla U(u_j^{N,*}(r,\xi))-\mathbb E\nabla U(u_1^{*}(r,\xi)) \notag\\
&\quad=\frac1N\sum_{j=1}^N\left[\nabla U(u_j^{N,*}(r,\xi))-\nabla U(u_j^{*}(r,\xi))\right]+\frac1{\sqrt N}\mathbb G_U^N(r,\xi),
\label{eq:single-drift-fluctuation-decomposition}
\end{align}
and
\begin{align}
&\quad\frac1N\sum_{j=1}^N\nabla K(u_j^{N,*}(r,\xi)-u_i^{N,*}(r,\xi))-\widetilde{\mathbb E}\nabla K(\widetilde u_1^{*}(r,\xi)-u_i^{*}(r,\xi)) \notag\\
&=\frac1N\sum_{j=1}^N\left[\nabla K(u_j^{N,*}(r,\xi)-u_i^{N,*}(r,\xi))-\nabla K(u_j^{*}(r,\xi)-u_i^{*}(r,\xi))\right]
+\frac1{\sqrt N}\mathbb G_K^N(r,\xi,u_i^{*}(r,\xi)).
\label{eq:interaction-fluctuation-decomposition}
\end{align}

\begin{proposition}\label{prop:driving-field-CLT}
Assume \textup{$(\mathbf{H}1)$--$(\mathbf{H}3)$}. For every $T,R,M>0$,
\[
\left(\mathbb G_U^N,\mathbb G_K^N,\bar W^N\right)\Longrightarrow\left(\mathbb G_U,\mathbb G_K,\bar W\right)
\]
in $C([-T,0]\times B_R^{cu},H)\times C([-T,0]\times B_R^{cu}\times K_M,H)\times C([-T,0],H).$
The limit is jointly Gaussian. For any finite families of parameter points
$\{(r_a,\xi_a)\}_{a=1}^{m_1}$, $\{(s_b,\zeta_b,z_b)\}_{b=1}^{m_2}$, and $\{t_c\}_{c=1}^{m_3}$, the Gaussian vector
$\left(\bigl(\mathbb G_U(r_a,\xi_a)\bigr)_{a=1}^{m_1},\bigl(\mathbb G_K(s_b,\zeta_b,z_b)\bigr)_{b=1}^{m_2},\bigl(\bar W(t_c)\bigr)_{c=1}^{m_3}\right)$
has mean zero and the same covariance operator as the random vector
\begin{equation}\label{eq:driving-field-covariance}
\left(
\begin{aligned}
&\bigl(\nabla U(u_1^{*}(r_a,\xi_a))-\mathbb E\nabla U(u_1^{*}(r_a,\xi_a))\bigr)_{1\le a\le m_1},\\
&\bigl(\nabla K(u_1^{*}(s_b,\zeta_b)-z_b)-\widetilde{\mathbb E}\nabla K(\widetilde u_1^{*}(s_b,\zeta_b)-z_b)\bigr)_{1\le b\le m_2},\\
&\bigl(W_1(t_c)\bigr)_{1\le c\le m_3}
\end{aligned}
\right).
\end{equation}
\end{proposition}

\begin{proof}
For any finite families of parameter points as in the statement, the Hilbert space central limit theorem \cite{AraujoGine1980}, applied to the i.i.d. vectors formed by the three 
groups in the proposition, gives the finite dimensional convergence and the stated covariance.

For tightness of $\mathbb G_K^N$, equip $\Theta_{T,R,M}:=[-T,0]\times B_R^{cu}\times K_M$
with the metric $d((r,\xi,z),(s,\zeta,y)):=|r-s|^\rho+|\xi-\zeta|+|z-y|.$
Set
\[
F_j^K(r,\xi,z):=\nabla K(u_j^{*}(r,\xi)-z),\qquad H_j(T,R):=\sup_{|\xi|\le R}[u_j^{*}(\cdot,\xi)]_{C^\rho([-T,0],H)}.
\]
Lemma~\ref{lem:backward-regularity-gaussian} gives
\begin{equation}\label{eq:gaussian-Hj-moment}
\mathbb EH_1(T,R)^p<\infty.
\end{equation}
By \textup{$(\mathbf{H}2)$} and \eqref{eq:backward-xi-Lipschitz},
\[
|F_j^K(\theta)-F_j^K(\theta')|\le C\bigl(1+H_j(T,R)\bigr)d(\theta,\theta'), \qquad \theta,\theta'\in\Theta_{T,R,M}.
\]
Define $d_{N,K}(\theta,\theta'):=\left[\frac1N\sum_{j=1}^N |F_j^K(\theta)-F_j^K(\theta')|^2\right]^{1/2}.$
Then
\begin{equation}\label{eq:gaussian-LN-bound}
d_{N,K}(\theta,\theta')\le L_N^K d(\theta,\theta'), \qquad L_N^K:=C\left[N^{-1}\sum_{j=1}^N(1+H_j(T,R))^2\right]^{1/2},
\end{equation}
and \eqref{eq:gaussian-Hj-moment} implies $\sup_{N\ge1}\mathbb E(L_N^K)^p<\infty.$

For a set $S$ equipped with a pseudometric $d_*$, let $N(S,d_*,u)$ denote the smallest number of balls of radius $u$,
with centers in $S$, needed to cover $S$. By \eqref{eq:gaussian-LN-bound}, every $d$-ball of radius
$u/L_N^K$ is contained in the $d_{N,K}$-ball of radius $u$ with the same center. Hence,
\begin{equation}\label{eq:empirical-covering-comparison}
N(\Theta_{T,R,M},d_{N,K},u) \le N\!\left(\Theta_{T,R,M},d,\frac{u}{L_N^K}\right),
\end{equation}
Let $(\epsilon_j)_{j=1}^N$ be independent Rademacher random
variables, independent of $(F_j^K)_{j=1}^N$, and set
\[
R_N^K(\theta):=\frac1{\sqrt N}\sum_{j=1}^N\epsilon_jF_j^K(\theta),\qquad \theta\in\Theta_{T,R,M}.
\]
Let $\mathbb E_\epsilon$ denote expectation with respect to
$(\epsilon_j)_{j=1}^N$, with $(F_j^K)_{j=1}^N$ fixed. The symmetrization inequality
\cite[Lemma~6.3 and the subsequent remark]{LedouxTalagrand1991} gives
\[
\begin{aligned}
 \mathbb E\sup_{d(\theta,\theta')\le\delta}\left|\mathbb G_K^N(\theta)-\mathbb G_K^N(\theta')\right|\le2\mathbb E\mathbb E_\epsilon
\sup_{d(\theta,\theta')\le\delta}\left|R_N^K(\theta)-R_N^K(\theta')\right|.
\end{aligned}
\]
Conditionally on $(F_j^K)_{j=1}^N$, $R_N^K$ is a Hilbert space valued Rademacher process, and
\[
\left(\mathbb E_\epsilon\left|R_N^K(\theta)-R_N^K(\theta')\right|^2\right)^{1/2}=d_{N,K}(\theta,\theta').
\]
By \eqref{eq:gaussian-LN-bound}, $d(\theta,\theta')\le\delta$ implies $d_{N,K}(\theta,\theta')\le L_N^K\delta$.
Applying the local entropy estimate obtained in the proof of Theorem~11.6, together with the extension to Banach space valued
processes in Remark~11.5 of \cite{LedouxTalagrand1991}, we obtain
\[
\begin{aligned}
\mathbb E_\epsilon\sup_{d(\theta,\theta')\le\delta}\left|R_N^K(\theta)-R_N^K(\theta')\right| &
\le\mathbb E_\epsilon\sup_{d_{N,K}(\theta,\theta')\le L_N^K\delta}\left|R_N^K(\theta)-R_N^K(\theta')\right|\\
&\qquad\le C\int_0^{2L_N^K\delta}\sqrt{\log N(\Theta_{T,R,M},d_{N,K},u)}\,du.
\end{aligned}
\]
By \eqref{eq:empirical-covering-comparison} and the substitution $u=L_N^K v$,
\[
\begin{aligned}\int_0^{2L_N^K\delta}\sqrt{\log N(\Theta_{T,R,M},d_{N,K},u)}\,du& \le\int_0^{2L_N^K\delta}\sqrt{\log N\!\left(\Theta_{T,R,M},d,\frac{u}{L_N^K}\right)}\,du\\
&=L_N^K\int_0^{2\delta}\sqrt{\log N(\Theta_{T,R,M},d,v)}\,dv.
\end{aligned}
\]
Taking expectations and combining these inequalities yield
\[
\begin{aligned}
 \mathbb E\sup_{d(\theta,\theta')\le\delta}
\left|\mathbb G_K^N(\theta)-\mathbb G_K^N(\theta')\right| \le C\,\mathbb E L_N^K\int_0^{2\delta}\sqrt{\log N(\Theta_{T,R,M},d,v)}\,dv.
\end{aligned}
\]
Since $\sup_{N\ge1}\mathbb E L_N^K<\infty$, Markov's inequality gives
\begin{equation}\label{eq:gaussian-GK-increment-bound}
\sup_{N\ge1}\mathbb P\left(\sup_{d(\theta,\theta')\le\delta}\left|\mathbb G_K^N(\theta)-\mathbb G_K^N(\theta')\right|>\varepsilon\right)
\le \frac C\varepsilon\int_0^{2\delta}\sqrt{\log N(\Theta_{T,R,M},d,u)}\,du,
\end{equation}
where $C$ is independent of $N$, $\delta$, and $\varepsilon$. To estimate the entropy integral, write $z=\sum_{k\ge1}z_ke_k\in K_M$. The definition of $V_\vartheta$ gives
\[
\sum_{k>m}|z_k|^2\le M^2\bigl(1+(m+1)^2\bigr)^{-2\vartheta}.
\]
For $0<\varepsilon\le M$, choose $m=\lceil(2M/\varepsilon)^{1/(2\vartheta)}\rceil$, so that the norm of this tail is at most $\varepsilon/2$.
The projection of $K_M$ onto $\operatorname{span}\{e_1,\ldots,e_m\}$ is contained in the Euclidean ball of radius $M$.
The volume bound in \cite[Proposition~4.2.12]{Vershynin2018} gives an $\varepsilon/2$ net in this projection with at most $(1+4M/\varepsilon)^m$ points.
These points belong to $K_M$ and form an $\varepsilon$ net for $K_M$. Hence
\begin{equation}\label{eq:KM-covering-bound}
\log N(K_M,|\cdot|,\varepsilon)
\le C_\vartheta\left(\frac M\varepsilon\right)^{1/(2\vartheta)}
\log\left(2+\frac M\varepsilon\right),\qquad 0<\varepsilon\le M.
\end{equation}
For $\varepsilon>M$, one ball centered at zero suffices. Since $\vartheta>1/4$,
\begin{equation}\label{eq:KM-entropy-integral}
\int_0^1\sqrt{\log N(K_M,|\cdot|,\varepsilon)}\,d\varepsilon<\infty.
\end{equation}
Since $H_{cu}$ is finite dimensional, the time and center variables satisfy a logarithmic entropy bound. Together with \eqref{eq:KM-covering-bound},
\[
\log N(\Theta_{T,R,M},d,u)\le C_{T,R}\log(2+u^{-1})+C_\vartheta\left(\frac Mu\right)^{1/(2\vartheta)}
\log\left(2+\frac Mu\right), \qquad 0<u\le1.
\]
Hence the right hand side of \eqref{eq:gaussian-GK-increment-bound} converges to zero as $\delta\downarrow0$. For $\nabla U(u_j^{*}(r,\xi))$,
$(r,\xi)\in[-T,0]\times B_R^{cu}$, the bounds \eqref{eq:gaussian-Hj-moment}--\eqref{eq:gaussian-GK-increment-bound} are obtained without the $z$ variable. 
They give the corresponding increment estimate for $\mathbb G_U^N$.

For each fixed $(r,\xi)\in[-T,0]\times B_R^{cu}$ and $z\in K_M$, the finite dimensional convergence implies that
$(\mathbb G_U^N(r,\xi))_N$ and $(\mathbb G_K^N(r,\xi,z))_N$ are tight in $H$. Therefore, the tightness criterion for continuous processes on
compact parameter sets \cite[Theorem~16.5]{Kallenberg2002} gives tightness of $(\mathbb G_U^N)_N$ and $(\mathbb G_K^N)_N$ in
$C([-T,0]\times B_R^{cu},H)~\text{and}~C([-T,0]\times B_R^{cu}\times K_M,H),$respectively.
Since $\bar W^N$ has the same law as $W_1$ in $C([-T,0],H)$ for every $N$, the sequence
$(\mathbb G_U^N,\mathbb G_K^N,\bar W^N)_N$ is tight in $C([-T,0]\times B_R^{cu},H)\times C([-T,0]\times B_R^{cu}\times K_M,H)\times C([-T,0],H).$
The joint finite dimensional convergence identifies every
subsequential limit and proves the proposition.

\end{proof}

For $t\in\mathbb R$, define
\[
\mathbb G_{U,\theta}^N(t,r,\xi):=\mathbb G_U^N(r,\xi,\theta_t\omega),\qquad
\mathbb G_{K,\theta}^N(t,r,\xi,z) :=\mathbb G_K^N(r,\xi,z,\theta_t\omega).
\]
The following corollary extends the joint convergence in Proposition~\ref{prop:driving-field-CLT} to
$\mathbb G_{U,\theta}^N$ and $\mathbb G_{K,\theta}^N$, which also depend on the time parameter $t$.

\begin{corollary}\label{cor:theta-driving-field-CLT}
Assume $(\mathbf{H}1)$--$(\mathbf{H}3)$. For every $T_0,T,R,M>0$ and every compact interval $I\subset\mathbb R$,
\begin{equation}
\left(\mathbb G_{U,\theta}^N,\mathbb G_{K,\theta}^N,\bar W^N|_I\right)
\Longrightarrow \left(\mathbb G_{U,\theta},\mathbb G_{K,\theta},\bar W|_I\right) \label{eq:theta-driving-field-CLT}
\end{equation}
in $C([0,T_0]\times[-T,0]\times B_R^{cu},H)\times C([0,T_0]\times[-T,0]\times B_R^{cu}\times K_M,H)\times C(I,H).$
The limit is jointly Gaussian.
\end{corollary}

\begin{proof}
Finite dimensional convergence follows from the Hilbert space central limit theorem applied to the i.i.d. summands defining
the three processes. Using \eqref{eq:backward-joint-increment-bound} with $1/p<\rho_0<\rho$, the tightness argument in the proof of
Proposition~\ref{prop:driving-field-CLT} applies, since the additional time parameter, with distance $|t-v|^{\rho_0}$,
contributes only a logarithmic term to the entropy bound.
\end{proof}

The next two lemmas give uniform $p$th moment bounds for the empirical fields and for $(r,\xi)\mapsto\mathbb G_K^N(r,\xi,u_i^{*}(r,\xi))$. 
The second lemma also gives $L^p$ equicontinuity. Their proofs are given in Appendix~\ref{app:empirical-moments}.

\begin{lemma}\label{lem:driving-field-local-pmoment}
For the fixed $p>4$ and every $T,R,M>0$,
\begin{align}
\sup_{N\ge1}\mathbb E \|\mathbb G_U^N\|_{C([-T,0]\times B_R^{cu},H)}^p&<\infty,\label{eq:GU-local-pmoment}\\
\sup_{N\ge1}\mathbb E\|\mathbb G_K^N\|_{C([-T,0]\times B_R^{cu}\times K_M,H)}^p&<\infty.\label{eq:GK-local-pmoment}
\end{align}
For every $T_0>0$,
\begin{align*}
\sup_{N\ge1}\mathbb E
\|\mathbb G_{U,\theta}^N\|_{C([0,T_0]\times[-T,0]\times B_R^{cu},H)}^p&<\infty,\\
\sup_{N\ge1}\mathbb E\|\mathbb G_{K,\theta}^N\|_{C([0,T_0]\times[-T,0]\times B_R^{cu}\times K_M,H)}^p&<\infty.
\end{align*}
\end{lemma}

\begin{lemma}\label{lem:particle-driving-forcing-pmoment}
Let $\rho_0$ satisfy $1/p<\rho_0<\rho<\vartheta-1/p$. For every $T,R>0$,
\begin{equation}\label{eq:particle-GK-random-state-pmoment}
\sup_{N\ge2}\max_{1\le i\le N} \mathbb E\sup_{\substack{-T\le r\le0\\|\xi|\le R}}\left|\mathbb G_K^N\bigl(r,\xi,u_i^{*}(r,\xi)\bigr)\right|^p<\infty.
\end{equation}
For every $T_0>0$,
\begin{equation}\label{eq:theta-particle-GK-random-state-pmoment}
\sup_{N\ge2}\max_{1\le i\le N} \mathbb E\sup_{\substack{0\le t\le T_0,\,-T\le r\le0\\|\xi|\le R}} \left|\mathbb G_{K,\theta}^N
\bigl(t,r,\xi,u_i^{*}(r,\xi,\theta_t\omega)\bigr)\right|^p<\infty.
\end{equation}
Moreover,
\begin{equation}\label{eq:particle-GK-Lp-equicontinuity}
\lim_{\delta\downarrow0}\sup_{N\ge2}\max_{1\le i\le N}\mathbb E\!\left[\sup_{\substack{-T\le r,s\le0,\ |\xi|,|\zeta|\le R\\
|r-s|^\rho+|\xi-\zeta|\le\delta}} \left| \mathbb G_K^N(r,\xi,u_i^{*}(r,\xi))-\mathbb G_K^N(s,\zeta,u_i^{*}(s,\zeta))\right|^p\right]=0,
\end{equation}
and
\begin{equation}\label{eq:theta-particle-GK-Lp-equicontinuity}
\lim_{\delta\downarrow0}
\sup_{N\ge2}\max_{1\le i\le N}\mathbb E\!\left[\sup_{\substack{0\le t,v\le T_0,\ -T\le r,s\le0\\|\xi|,|\zeta|\le R\\
|t-v|^{\rho_0}+|r-s|^\rho+|\xi-\zeta|\le\delta}}\left|
\begin{aligned}
&\mathbb G_{K,\theta}^N(t,r,\xi,u_i^{*}(r,\xi,\theta_t\omega))\\
&\quad-\mathbb G_{K,\theta}^N(v,s,\zeta,u_i^{*}(s,\zeta,\theta_v\omega))
\end{aligned}
\right|^p
\right]=0.
\end{equation}
\end{lemma}

\subsection{Linearization of the backward equations}

 We now prove Theorem~\ref{thm:intro-LP-fluctuation} by comparing
the finite particle backward solutions defined by \eqref{eq:LP-orbit-fixed-N} and \eqref{eq:original-backward-components}
with the limiting backward solutions of \eqref{eq:limiting-backward-LP-system}.
 We first establish a mean-square error bound of order $N^{-1}$. We then linearize the backward equations and use the joint convergence
of the empirical fields and the averaged noise to identify the limit of the fluctuations defined in \eqref{eq:particle-state-fluctuation}.

For $E\in\{H,\mathbf H\}$, a path $h:(-\infty,0]\to E$, and $T>0$, set
\[
\|h\|_\eta=\sup_{r\le0}e^{-\eta r}\|h(r)\|_E,
\qquad
\|h\|_{\eta,T}=\sup_{-T\le r\le0}e^{-\eta r}\|h(r)\|_E.
\]
Let $\hat V_i^{N,*}$ be defined by 
$\hat V_i^{N,*}(r,\xi,\omega)=\hat u_i^{N,*}(r,\xi,\omega)-\hat z_i^N(\theta_r\omega)$. Set $q_s=\frac{\ell_{\mathrm{LP}}}{\eta-\beta}\le q<1.$
By \eqref{eq:LP-orbit-fixed-N}, \eqref{eq:fixed-N-uniform-dichotomy}, and \eqref{eq:FNi-particle-Lipschitz},
\begin{align*}
\|\hat V_i^{N,*}(\cdot,\xi)-\hat V_i^{N,*}(\cdot,\zeta)\|_{\eta} \le q_s\|\hat V_i^{N,*}(\cdot,\xi)-\hat V_i^{N,*}(\cdot,\zeta)\|_{\eta} 
+ q_s\|V^{N,*}(\cdot,\xi)-V^{N,*}(\cdot,\zeta)\|_{\eta,N}.
\end{align*}
By \eqref{eq:fixed-N-graph-empirical-Lipschitz}, the preceding inequality, and the cancellation of the Ornstein--Uhlenbeck terms,
there exists a deterministic constant $C>0$, independent of $N,\xi$, and $\zeta$, such that, almost surely,
\begin{equation}\label{eq:fixed-N-particle-center-Lipschitz}
\sup_{1\le i\le N} \|\mathbf U_i^{N,*}(\cdot,\xi) -\mathbf U_i^{N,*}(\cdot,\zeta)\|_{\eta} \le C|\xi-\zeta|.
\end{equation}

\begin{proposition}
\label{prop:backward-graph-N-rate}
Assume $(\mathbf{H}1)$--$(\mathbf{H}3)$. Let $V^{N,*}$ solve \eqref{eq:LP-orbit-fixed-N}, and let
$\mathbf U_i^{*}=(\bar u^{*},\hat u_i^{*})$ solve \eqref{eq:limiting-backward-LP-system}.
Define $\mathbf U_i^{N,*}=(\bar u^{N,*},\hat u_i^{N,*})$
from $V^{N,*}$ through \eqref{eq:original-backward-components}. For every $T,R>0$, there exists $C_{T,R}<\infty$ such that
\begin{equation}\label{eq:backward-state-N-rate}
\max_{1\le i\le N} \sup_{|\xi|\le R} \mathbb E \left\| \mathbf U_i^{N,*}(\cdot,\xi)-\mathbf U_i^{*}(\cdot,\xi) \right\|_{\eta,T}^{2} \le \frac{C_{T,R}}{N}, \qquad N\ge2.
\end{equation}
\end{proposition}

\begin{proof}
Fix $T,R>0$ and $|\xi|\le R$, and set $\Delta_i(r,\xi):=\mathbf U_i^{N,*}(r,\xi)-\mathbf U_i^{*}(r,\xi),~ r\le0.$
Let $z^N$ be the stationary solution of \eqref{eq:joint-stationary-OU}. The representation \eqref{eq:OU-mean-fluctuation-components}
and the maximal inequality for stochastic convolutions \cite{DaPratoZabczyk2014} give $\mathbb E\sup_{-1\le r\le0} \|z^N(\theta_r\omega)\|_N^2<\infty.$
By stationarity and $\eta<0$,
\[
\begin{aligned}
\mathbb E\|z^N(\theta_\cdot\omega)\|_{\eta,N}^2&\le\sum_{k=0}^{\infty}e^{2\eta k}\mathbb E\sup_{-1\le r\le0}\|z^N(\theta_r\omega)\|_N^2<\infty.
\end{aligned}
\]
Equations \eqref{eq:LP-orbit-fixed-N} and \eqref{eq:transformed-random-nonlinearity}, together with
\eqref{eq:fixed-N-uniform-dichotomy} and $(\mathbf{H}2)$--$(\mathbf{H}3)$, give
\[
(1-q)\|V^{N,*}(\cdot,\xi)\|_{\eta,N}\le C_R\left(1+\|z^N(\theta_\cdot\omega)\|_{\eta,N}\right).
\]
Hence $\mathbb E\|V^{N,*}(\cdot,\xi)\|_{\eta,N}^2<\infty$. Using \eqref{eq:original-backward-components} and
\eqref{eq:backward-global-weighted-moment}, we conclude that
\[
\mathbb E\|\Delta_i(\cdot,\xi)\|_\eta^2<\infty
\]
for every fixed $N$ and $1\le i\le N$.
For deterministic $\xi$, the solutions $\mathbf U_j^{*}(\cdot,\xi)$ of \eqref{eq:limiting-backward-LP-system}, driven respectively
by $W_j$, are independent and identically distributed. For fixed $r\le0$, write $u_j^{*}=u_j^{*}(r,\xi)$
 and $\widetilde u_1^{*}=\widetilde u_1^{*}(r,\xi)$.
Set $\mathcal G_i=\sigma(W_i(s):s\in\mathbb R)$ and define
\[
F_{j|i}:=\nabla K(u_j^{*}-u_i^{*})-\widetilde{\mathbb E}\nabla K(\widetilde u_1^{*}-u_i^{*}),\qquad 1\le j\le N.
\]
Conditionally on $\mathcal G_i$, the variables $F_{j|i}$, $j\ne i$, are independent and satisfy $\mathbb E[F_{j|i}\mid\mathcal G_i]=0.$
The term $F_{i|i}$ is $\mathcal G_i$-measurable.
By $(\mathbf{H}2)$,
\[
\mathbb E|F_{j|i}|^2 \le L_K^2\mathbb E|u_1^{*}|^2\quad(j\ne i),\qquad\mathbb E|F_{i|i}|^2\le4L_K^2\mathbb E|u_1^{*}|^2.
\]
Thus \eqref{eq:empirical-GU-definition} and
\eqref{eq:empirical-GK-definition} give
\[
\begin{aligned}
\mathbb E|\mathbb G_U^N(r,\xi)|^2&=\mathbb E\left|\nabla U(u_1^{*})-\mathbb E\nabla U(u_1^{*})\right|^2\le L_U^2\mathbb E|u_1^{*}|^2,\\
\mathbb E\left|\mathbb G_K^N(r,\xi,u_i^{*})\right|^2&=\frac1N\sum_{j=1}^N\mathbb E|F_{j|i}|^2\le\frac{N+3}{N}L_K^2\mathbb E|u_1^{*}|^2.
\end{aligned}
\]
By \eqref{eq:backward-global-weighted-moment}, it follows that
\begin{equation}\label{eq:backward-comparison-field-bound}
\mathbb E|\mathbb G_U^N(r,\xi)|^2+\max_{1\le i\le N}\mathbb E\left|
\mathbb G_K^N(r,\xi,u_i^{*}(r,\xi))
\right|^2
\le C_Re^{2\eta r},
\qquad r\le0.
\end{equation}
Define $D_N(\xi):=\max_{1\le i\le N}\sup_{r\le0}e^{-\eta r}\|\Delta_i(r,\xi)\|_{L^2(\Omega;\mathbf H)}.$
The weighted moment estimates above imply $D_N(\xi)<\infty$. For the nonlinear drift differences, set
\[
\begin{aligned}
a_N(r,\xi)&:=\frac1N\sum_{j=1}^N\left[\nabla U(u_j^{N,*}(r,\xi))-\nabla U(u_j^{*}(r,\xi))\right],\\
b_{N,i}(r,\xi)&:=\nabla U(u_i^{N,*}(r,\xi))-\nabla U(u_i^{*}(r,\xi))-a_N(r,\xi)\\
&\quad+\frac{\gamma}{N}\sum_{j=1}^N\left[\nabla K(u_j^{N,*}(r,\xi)-u_i^{N,*}(r,\xi))-\nabla K(u_j^{*}(r,\xi)-u_i^{*}(r,\xi))\right].
\end{aligned}
\]
Since $|u_j^{N,*}(r,\xi)-u_j^{*}(r,\xi)|\le\sqrt2\,|\Delta_j(r,\xi)|_{\mathbf H},$  $(\mathbf{H}2)$ and Minkowski's inequality yield
\begin{equation}\label{eq:backward-comparison-drift-bound}
\begin{aligned}
\|a_N(r,\xi)\|_{L^2(\Omega;H)}
&\le \sqrt2 L_U e^{\eta r}D_N(\xi),\\\|b_{N,i}(r,\xi)\|_{L^2(\Omega;H)}
&\le 2\sqrt2(L_U+\gamma L_K)e^{\eta r}D_N(\xi).
\end{aligned}
\end{equation}
In particular, using  $3\sqrt2(L_U+\gamma L_K)\le\ell_{\mathrm{LP}}$
\[
\left\|
\begin{pmatrix}
a_N(r,\xi)\\ b_{N,i}(r,\xi)
\end{pmatrix}
\right\|_{L^2(\Omega;\mathbf H)}\le\ell_{\mathrm{LP}}e^{\eta r}D_N(\xi).
\]
By \eqref{eq:original-backward-components}, \eqref{eq:LP-orbit-fixed-N}, and \eqref{eq:limiting-backward-LP-system}, the contribution of
$\bar W^N$ to the difference of the backward equations is
\[
\begin{aligned}
h_N(r):={}&
\begin{pmatrix}
e^{L_{cu}r}z_{cu}^N(\omega)\\0
\end{pmatrix}+\frac{\sigma}{\sqrt N}
\int_0^r e^{\mathbf L_{cu}(r-s)}\mathbf P_{cu}
\begin{pmatrix}
d\bar W^N(s)\\-d\bar W^N(s)
\end{pmatrix}\\
&+\frac{\sigma}{\sqrt N}\int_{-\infty}^r e^{\mathbf L_s(r-s)}\mathbf P_s
\begin{pmatrix}
d\bar W^N(s)\\-d\bar W^N(s)
\end{pmatrix},
\qquad r\le0.
\end{aligned}
\]
Equation \eqref{eq:OU-mean-fluctuation-components}, the It\^o
isometry, and \eqref{eq:fixed-N-uniform-dichotomy} give the first estimate below, the maximal inequality for stochastic
convolutions gives the second:
\begin{equation}\label{eq:backward-comparison-noise-bound}
\sup_{r\le0}e^{-2\eta r}
\mathbb E|h_N(r)|_{\mathbf H}^2
\le\frac C N,
\qquad
\mathbb E\|h_N\|_{\eta,T}^2
\le\frac{C_T}{N}.
\end{equation}
Subtract \eqref{eq:limiting-backward-LP-system} from the finite particle backward equations obtained from
\eqref{eq:LP-orbit-fixed-N} and \eqref{eq:original-backward-components}. By \eqref{eq:single-drift-fluctuation-decomposition} and
\eqref{eq:interaction-fluctuation-decomposition}, the nonlinear drift difference is
\[
\begin{pmatrix}
a_N(r,\xi)\\ b_{N,i}(r,\xi)\end{pmatrix} +\frac1{\sqrt N}
\begin{pmatrix}\mathbb G_U^N(r,\xi)\\-\mathbb G_U^N(r,\xi)+\gamma\mathbb G_K^N(r,\xi,u_i^{*}(r,\xi))
\end{pmatrix}.
\]
Applying Minkowski's inequality to the center and stable integrals, and using \eqref{eq:fixed-N-uniform-dichotomy} and
\eqref{eq:backward-comparison-field-bound}--\eqref{eq:backward-comparison-noise-bound}, we obtain
\begin{equation}\label{eq:backward-comparison-contraction}
\begin{aligned}
D_N(\xi)\le\ell_{\mathrm{LP}}\left(\frac1{\alpha-\eta}+\frac1{\eta-\beta}\right)D_N(\xi)+C_RN^{-1/2} =qD_N(\xi)+C_RN^{-1/2}.
\end{aligned}
\end{equation}
Since $D_N(\xi)<\infty$ and $q<1$, \eqref{eq:backward-comparison-contraction} implies
\begin{equation}\label{eq:backward-state-weighted-L2-rate}
\max_{1\le i\le N}\sup_{|\xi|\le R}\sup_{r\le0}
e^{-2\eta r}
\mathbb E|\Delta_i(r,\xi)|_{\mathbf H}^2
\le\frac{C_R}{N},
\qquad N\ge2.
\end{equation}
To prove \eqref{eq:backward-state-N-rate}, set $D_{N,T}(\xi):=\max_{1\le i\le N}\left(\mathbb E\|\Delta_i(\cdot,\xi)\|_{\eta,T}^2\right)^{1/2}.$
For $-T\le r\le0$, rewrite the stable equation with initial time $-T$. Its initial term is $e^{\mathbf L_s(r+T)}\mathbf P_s\Delta_i(-T,\xi)$.
By \eqref{eq:backward-state-weighted-L2-rate} and \eqref{eq:fixed-N-uniform-dichotomy},
\[
\left(
\mathbb E\sup_{-T\le r\le0}e^{-2\eta r}\left|e^{\mathbf L_s(r+T)}\mathbf P_s\Delta_i(-T,\xi)\right|_{\mathbf H}^2\right)^{1/2} \le C_RN^{-1/2}.
\]
Moreover,  the stochastic contribution is $h_N(r)-e^{\mathbf L_s(r+T)}\mathbf P_sh_N(-T).$ Its norm $\bigl(\mathbb E\|\cdot\|_{\eta,T}^2\bigr)^{1/2}$ is bounded by $C_TN^{-1/2}$,
 by \eqref{eq:backward-comparison-noise-bound} and \eqref{eq:fixed-N-uniform-dichotomy}. Since $p>2$, 
the estimates \eqref{eq:GU-local-pmoment} and \eqref{eq:particle-GK-random-state-pmoment} give
\[
\begin{aligned}
&\left(\mathbb E\|\mathbb G_U^N(\cdot,\xi)\|_{\eta,T}^2\right)^{1/2} +\max_{1\le i\le N}\left(\mathbb E\left\| r\mapsto\mathbb G_K^N(r,\xi,u_i^{*}(r,\xi))\right\|_{\eta,T}^2
\right)^{1/2} \le C_{T,R}.
\end{aligned}
\]
Applying $(\mathbf{H}2)$ and Minkowski's inequality to $a_N$ and $b_{N,i}$ in the norm $\bigl(\mathbb E\|\cdot\|_{\eta,T}^2\bigr)^{1/2}$
gives the constants in \eqref{eq:backward-comparison-drift-bound}, with $e^{\eta r}D_N(\xi)$ replaced by $D_{N,T}(\xi)$.
The center and stable integrals over $[-T,0]$ are bounded by the same exponential kernel integrals as in
\eqref{eq:backward-comparison-contraction}. Therefore,
\[
D_{N,T}(\xi) \le qD_{N,T}(\xi)+C_{T,R}N^{-1/2}.
\]
Since $q<1$ and the constants are uniform for $|\xi|\le R$, this proves \eqref{eq:backward-state-N-rate}.
\end{proof}

Recall the fluctuations in \eqref{eq:particle-state-fluctuation}. Write $\mathbf Y_i^N=(\bar y^N,\hat y_i^N)$, so that
$y_i^N=\bar y^N+\hat y_i^N$, and set $\mathbf Y^N=(\bar y^N,\hat y_1^N,\ldots,\hat y_N^N)$.
For $F\in C_\eta^-(\mathbf H)$, define
\begin{equation}\label{eq:LP-convolution-operator}
(\mathscr K F)(t)=\int_0^t e^{\mathbf L_{cu}(t-r)}\mathbf P_{cu}F(r)\,dr +\int_{-\infty}^t e^{\mathbf L_s(t-r)}\mathbf P_sF(r)\,dr, \qquad t\le0.
\end{equation}
Using \eqref{eq:original-backward-components}, \eqref{eq:LP-orbit-fixed-N}, and \eqref{eq:limiting-backward-LP-system}, 
subtract the limiting backward equation from the finite particle equation and multiply by $\sqrt N$. Adding and subtracting 
the empirical drift at the limiting backward particles gives
\begin{equation}\label{eq:exact-mean-value-LP-fluctuation}
\mathbf Y_i^N =\mathscr K\left[ \mathfrak D_{N,i}^{\xi}(\mathbf Y^N)+\mathcal B_i^N \right]+\mathbb Z^N.
\end{equation}
The terms $\mathfrak D_{N,i}^{\xi}$, $\mathcal B_i^N$, and $\mathbb Z^N$ are defined by  \eqref{eq:exact-mean-value-drift}, 
\eqref{eq:driving-particle-forcing-N}, and \eqref{eq:wiener-backward-contribution} in the following.

For $1\le j\le N$, define
\[A_{U,j}^N(r,\xi) :=\int_0^1 D\nabla U\!\left( u_j^{*}(r,\xi) +\theta\bigl(u_j^{N,*}(r,\xi)-u_j^{*}(r,\xi)\bigr) \right)\,d\theta,\]
and, for $1\le i,j\le N$,
\[
\begin{aligned}
A_{K,ij}^N(r,\xi):=\int_0^1 D\nabla K\Bigl(&u_j^{*}(r,\xi)-u_i^{*}(r,\xi)+\theta\bigl[(u_j^{N,*}-u_j^{*})-(u_i^{N,*}
-u_i^{*})\bigr](r,\xi)\Bigr)\,d\theta.
\end{aligned}
\]
For $\mathbf y^N=(\bar y,\hat y_1,\ldots,\hat y_N)
\in C_\eta^-(H^{\oplus(N+1)})$, write
$\mathbf y_j=(\bar y,\hat y_j)$ and $y_j=\bar y+\hat y_j$, and set
\begin{equation}\label{eq:exact-mean-value-drift}
\mathfrak D_{N,i}^{\xi}(r,\mathbf y^N) := \begin{pmatrix} \displaystyle \frac1N\sum_{j=1}^N A_{U,j}^N(r,\xi)y_j(r) \\[3mm] \begin{aligned} 
&A_{U,i}^N(r,\xi)y_i(r) -\frac1N\sum_{j=1}^N A_{U,j}^N(r,\xi)y_j(r)\\ &\qquad +\frac{\gamma}{N}\sum_{j=1}^N A_{K,ij}^N(r,\xi)\bigl(y_j(r)-y_i(r)\bigr). 
\end{aligned} \end{pmatrix}.
\end{equation}
The empirical fluctuation term is
\begin{equation}\label{eq:driving-particle-forcing-N}
\mathcal B_i^N(r,\xi) := \begin{pmatrix} \mathbb G_U^N(r,\xi) \\[1mm] -\mathbb G_U^N(r,\xi) 
+\gamma\mathbb G_K^N\bigl(r,\xi,u_i^{*}(r,\xi)\bigr) \end{pmatrix}.
\end{equation}
By \eqref{eq:bar-u-mean-fluctuation-system}, the contribution of
$\bar W^N$ to the average equation is $\sigma N^{-1/2}d\bar W^N$.
Since $\hat W_i^N=W_i-N^{-1/2}\bar W^N$, its contribution to
\eqref{eq:hat-u-equation} is $-\sigma N^{-1/2}d\bar W^N$.
After subtraction of the limiting equations and multiplication by
$\sqrt N$, then
\begin{equation}\label{eq:wiener-backward-contribution}
\begin{aligned} \mathbb Z^N(t,\omega):={}& \begin{pmatrix} e^{L_{cu}t}\sqrt N\,z_{cu}^N(\omega)\\ 0 \end{pmatrix} +\sigma\int_0^t
 e^{\mathbf L_{cu}(t-r)}\mathbf P_{cu} \begin{pmatrix} d\bar W^N(r,\omega)\\ 
-d\bar W^N(r,\omega) \end{pmatrix}\\ &+\sigma\int_{-\infty}^t e^{\mathbf L_s(t-r)}\mathbf P_s \begin{pmatrix} d\bar W^N(r,\omega)\\ 
-d\bar W^N(r,\omega) \end{pmatrix}, \qquad t\le0. \end{aligned}
\end{equation}
The integral form of the mean value theorem gives the drift difference in \eqref{eq:exact-mean-value-drift}.
With the coefficients $A_{U,j}^N$ and $A_{K,ij}^N$ fixed,
$\mathfrak D_{N,i}^{\xi}(r,\cdot)$ is linear.
By $(\mathbf{H}2)$, $\|A_{U,j}^N\|\le L_U$ and $\|A_{K,ij}^N\|\le L_K$.
By \eqref{eq:FN-uniform-Lipschitz}, \eqref{eq:fixed-N-uniform-dichotomy}, and \eqref{eq:standing-LP-gap}, the operator
$\mathbf y^N\mapsto\mathscr K\mathfrak D_{N,\cdot}^{\xi}(\mathbf y^N)$ has contraction constant at most $q$ on
$C_\eta^-(H^{\oplus(N+1)})$, with norm $\sup_{r\le0}e^{-\eta r}\|\mathbf y^N(r)\|_N$,
where $\|\cdot\|_N$ is defined in \eqref{eq:empirical-norm}. By \eqref{eq:wiener-backward-contribution},
\begin{equation}\label{eq:wiener-center-at-zero}
\mathbf P_{cu}\mathbb Z^N(0) = \begin{pmatrix} \sqrt N\,z_{cu}^N(\omega)\\ 0 \end{pmatrix}.
\end{equation}
By \eqref{eq:OU-mean-fluctuation-components}, all terms in
\eqref{eq:wiener-backward-contribution} depend linearly on the same
process $\bar W^N$. Since $\bar W^N$ has the same distribution as $W_1$,
the distribution of $\mathbb Z^N$ does not depend on $N$.

The next lemma bounds the particle fluctuations and the averaged Wiener contribution uniformly in $N$ and in the center and time parameters.
 The weighted bounds control the integrals over $(-\infty,0]$, while the local path bounds will be used for tightness and for the Taylor remainder.

\begin{lemma}
\label{lem:particle-fluctuation-weighted-parameter-bound}
Under the hypotheses of Proposition~\ref{prop:backward-graph-N-rate}, let $T_0,R>0$, $J=[0,T_0]\times B_R^{cu}$, and $2\le s\le p$ for the fixed $p>4$. Then
\begin{equation}\label{eq:particle-fluctuation-weighted-higher-moment}
\sup_{N\ge2}\max_{i\le N}\sup_{r\le0}e^{-s\eta r} \mathbb E\sup_{(t,\xi)\in J} |\mathbf Y_i^N(r,\xi,\theta_t\omega)|^s<\infty,
\end{equation}
and, for every $T>0$,
\begin{equation}\label{eq:particle-fluctuation-local-higher-moment}
\sup_{N\ge2}\max_{i\le N}\mathbb E\sup_{(t,\xi)\in J} \|\mathbf Y_i^N(\cdot,\xi,\theta_t\omega)\|_{\eta,T}^{s}<\infty.
\end{equation}
For every $2\le s\le p$, $\mathbb Z^N$ satisfies
\begin{equation}\label{eq:averaged-noise-weighted-theta-bound}
\sup_{N\ge2}\sup_{r\le0}e^{-s\eta r} \mathbb E\sup_{0\le t\le T_0}|\mathbb Z^N(r,\theta_t\omega)|^s<\infty,
\end{equation}
and its local path bound is
\begin{equation}\label{eq:averaged-noise-local-theta-pmoment}
\sup_{N\ge2}\mathbb E\sup_{0\le t\le T_0} \|\mathbb Z^N(\cdot,\theta_t\omega)\|_{\eta,T}^s<\infty.
\end{equation}
In particular, both estimates hold at $s=2$.
\end{lemma}

\begin{proof}
It is enough to prove the estimates for $s=p$. For $2\le s\le p$, the bounds follow from Jensen inequality applied to the concave function 
$x\mapsto x^{s/p}$.

 Let $z^N$ be the stationary solution of \eqref{eq:joint-stationary-OU}. By \eqref{eq:OU-mean-fluctuation-components}, stationarity,
and the maximal inequality for stochastic convolutions \cite{DaPratoZabczyk2014},
\[
\mathbb E\sup_{0\le t\le T_0}\|z^N(\theta_{\cdot+t}\omega)\|_{\eta,N}^p<\infty.
\]
The estimate on $(-\infty,0]$ follows by summing the estimates on $[-k-1,-k]$ with weights $e^{p\eta k}$. Equations \eqref{eq:LP-orbit-fixed-N} and
\eqref{eq:transformed-random-nonlinearity}, together with \eqref{eq:fixed-N-uniform-dichotomy} and $(\mathbf{H}2)$--$(\mathbf{H}3)$, give
\[
\begin{aligned}
(1-q)\sup_{(t,\xi)\in J}\|V^{N,*}(\cdot,\xi,\theta_t\omega)\|_{\eta,N}\le C_R\left(1+\sup_{0\le t\le T_0}\|z^N(\theta_{\cdot+t}\omega)\|_{\eta,N}\right).
\end{aligned}
\]
Combining this estimate with
\eqref{eq:original-backward-components},
\eqref{eq:particle-state-fluctuation},
\eqref{eq:backward-global-weighted-moment}, and
\eqref{eq:backward-theta-weighted-holder}, we obtain
\[
\mathbb E\sup_{(t,\xi)\in J}
\|\mathbf Y_i^N(\cdot,\xi,\theta_t\omega)\|_\eta^p
<\infty
\]
for every fixed $N$ and $1\le i\le N$. Fix $1/p<\rho_0<\rho$ as in Lemma~\ref{lem:backward-theta-regularity-gaussian}, and set
\[
d_J\bigl((t,\xi),(v,\zeta)\bigr):=|t-v|^{\rho_0}+|\xi-\zeta|.
\]
By \eqref{eq:backward-global-weighted-moment}, \eqref{eq:backward-xi-weighted-Lipschitz}, and \eqref{eq:backward-theta-weighted-holder}, there exist
nonnegative random variables $\Lambda_i$ and satisfy
$\sup_i\mathbb E\Lambda_i^p<\infty$. such that
\begin{equation}\label{eq:weighted-backward-parameter-controls}
\begin{aligned}
\sup_{r\le0}e^{-\eta r}\sup_{(t,\xi)\in J}|u_i^{*}(r,\xi,\theta_t\omega)|&\le\Lambda_i,\\
\sup_{r\le0}e^{-\eta r}|u_i^{*}(r,\xi,\theta_t\omega)-u_i^{*}(r,\zeta,\theta_v\omega)|&\le\Lambda_i d_J\bigl((t,\xi),(v,\zeta)\bigr).
\end{aligned}
\end{equation}
The variables $\Lambda_i$ do not depend on $r$.
Moreover, $\int_0^1\sqrt{\log N(J,d_J,\varepsilon)}\,d\varepsilon<\infty.$
Thus the symmetrization, entropy, and moment estimates in the proof of Lemma~\ref{lem:driving-field-local-pmoment}
apply to the weighted $\nabla U$ terms uniformly in $r$. For the interaction terms, use the conditional estimates
in the proof of Lemma~\ref{lem:particle-driving-forcing-pmoment}, with $\mathcal G_i=\sigma(W_i(s):s\in\mathbb R)$.
By \eqref{eq:weighted-backward-parameter-controls}, both the conditional sum over $j\ne i$ and the term with
$j=i$ satisfy the required moment bounds uniformly in $r$. Therefore,
\begin{align}
&\sup_{N\ge2}\sup_{r\le0}e^{-p\eta r}\mathbb E\sup_{(t,\xi)\in J}|\mathbb G_U^N(r,\xi,\theta_t\omega)|^p<\infty,\notag\\
&\sup_{N\ge2}\max_{i\le N}\sup_{r\le0}e^{-p\eta r}\mathbb E\sup_{(t,\xi)\in J}\left|\mathbb G_{K,\theta}^N\bigl(t,r,\xi,u_i^{*}(r,\xi,\theta_t\omega)\bigr)\right|^p<\infty.
\label{eq:weighted-driving-field-parameter-bound}
\end{align}
We next estimate $\mathbb Z^N$. By $(\mathbf{H}1)$, \eqref{eq:OU-mean-fluctuation-components}, and \eqref{eq:wiener-backward-contribution}, the distribution
of the field $(t,r)\mapsto\mathbb Z^N(r,\theta_t\omega)$ does not depend on $N$. The stable term in \eqref{eq:wiener-backward-contribution} is a stationary
stochastic convolution. Its maximal $p$th moments on intervals of fixed length are bounded uniformly in $N$. For the center component, integration by parts in the
finite dimensional stochastic integral, the Brownian maximal moment inequality, and \eqref{eq:OU-mean-fluctuation-components} give
\[
\sup_{N\ge2} \mathbb E\sup_{\substack{0\le t\le T_0\\-L\le r\le0}}|\mathbf P_{cu}\mathbb Z^N(r,\theta_t\omega)|_{\mathbf H}^p
\le C_{T_0,p}(1+L^{p/2}), \qquad L>0.
\]
Since $\sum_{k=0}^{\infty}e^{p\eta k}\bigl(1+(k+1)^{p/2}\bigr)<\infty,$ summing over $[-k-1,-k]$ proves \eqref{eq:averaged-noise-weighted-theta-bound}.
The estimates on bounded intervals also prove \eqref{eq:averaged-noise-local-theta-pmoment}.

Set
$
M_N:=\max_{i\le N}\sup_{r\le0}e^{-\eta r}\left(\mathbb E\sup_{(t,\xi)\in J}|\mathbf Y_i^N(r,\xi,\theta_t\omega)|_{\mathbf H}^p\right)^{1/p}.
$
Then $M_N<\infty$ for every fixed $N$.
For $\mathbf v^N=(\bar v,\hat v_1,\ldots,\hat v_N)$,
write $\mathbf v_i=(\bar v,\hat v_i)$.
By $(\mathbf{H}2)$ and \eqref{eq:exact-mean-value-drift},
\begin{equation}\label{eq:fluctuation-drift-component-bound}
|\mathfrak D_{N,i}^{\xi}(r,\mathbf v^N)|_{\mathbf H}
\le
\sqrt2(L_U+\gamma L_K)
\left(
|\mathbf v_i(r)|_{\mathbf H}
+\frac2N\sum_{j=1}^N|\mathbf v_j(r)|_{\mathbf H}
\right).
\end{equation}
Minkowski's inequality therefore gives a drift coefficientat most $3\sqrt2(L_U+\gamma L_K)\le\ell_{\mathrm{LP}}$
in the norm defining $M_N$. Apply \eqref{eq:exact-mean-value-LP-fluctuation} at $\theta_t\omega$.
Using \eqref{eq:LP-convolution-operator}, \eqref{eq:fixed-N-uniform-dichotomy}, \eqref{eq:fluctuation-drift-component-bound},
and Minkowski's inequality, we obtain
\begin{equation}\label{eq:fluctuation-weighted-moment-contraction}
\begin{aligned}
M_N&\le\ell_{\mathrm{LP}}\left(\frac1{\alpha-\eta}+\frac1{\eta-\beta}\right)M_N+C_{R,T_0,p}=qM_N+C_{R,T_0,p}.
\end{aligned}
\end{equation}
Here $\mathcal B_i^N$, defined in \eqref{eq:driving-particle-forcing-N}, is controlled by \eqref{eq:weighted-driving-field-parameter-bound},
and $\mathbb Z^N$ is controlled by
\eqref{eq:averaged-noise-weighted-theta-bound}.
Since $M_N<\infty$ and $q<1$, \eqref{eq:fluctuation-weighted-moment-contraction} proves \eqref{eq:particle-fluctuation-weighted-higher-moment}
for $s=p$.
For the local estimate, fix $T>0$ and set
$
M_{N,T}:=\max_{i\le N}\left(\mathbb E\sup_{(t,\xi)\in J}\|\mathbf Y_i^N(\cdot,\xi,\theta_t\omega)\|_{\eta,T}^p\right)^{1/p}.
$
By \eqref{eq:exact-mean-value-LP-fluctuation},
the stable integral over $(-\infty,-T]$ equals
\[
e^{\mathbf L_s(r+T)}\mathbf P_s\left[\mathbf Y_i^N(-T,\xi,\theta_t\omega)-\mathbb Z^N(-T,\theta_t\omega)\right],\qquad -T\le r\le0.
\]
Its norm $\bigl(\mathbb E\sup_{(t,\xi)\in J}\|\cdot\|_{\eta,T}^p\bigr)^{1/p}$ is bounded uniformly in $N$ and $i$, by
\eqref{eq:particle-fluctuation-weighted-higher-moment}, \eqref{eq:averaged-noise-weighted-theta-bound}, and \eqref{eq:fixed-N-uniform-dichotomy}.
Lemma~\ref{lem:driving-field-local-pmoment} and \eqref{eq:theta-particle-GK-random-state-pmoment} control $\mathcal B_i^N$ on $[-T,0]$, uniformly for
$(t,\xi)\in J$, and \eqref{eq:averaged-noise-local-theta-pmoment}
controls $\mathbb Z^N$. Applying Minkowski's inequality to the remaining
integrals over $[-T,0]$, and using
\eqref{eq:fluctuation-drift-component-bound}
and \eqref{eq:fixed-N-uniform-dichotomy}, gives
\[
M_{N,T}\le qM_{N,T}+C_{T,R,T_0,p}.
\]
Since $M_{N,T}<\infty$ and $q<1$, this proves
\eqref{eq:particle-fluctuation-local-higher-moment}
for $s=p$ and completes the proof.
\end{proof}

Set
\begin{equation}\label{eq:driving-field-triples}
\mathbb G^N=(\mathbb G_U^N,\mathbb G_K^N,\bar W^N),
\qquad \mathbb G=(\mathbb G_U,\mathbb G_K,\bar W).
\end{equation}
Let $\bar W$ be the limiting Wiener process in
\eqref{eq:theta-driving-field-CLT}.
The field $\mathbb Z$ is defined by \eqref{eq:wiener-backward-contribution} with $\bar W^N$ replaced by $\bar W$ and $\sqrt N z_{cu}^N(\omega)$
 replaced by $\sigma P_{cu}\int_{-\infty}^0e^{-\Delta r}\,d\bar W(r)$.
The next lemma extends Proposition~\ref{prop:driving-field-CLT} to include $\mathbb Z^N$ and proves that the limiting field
is independent of any fixed finite collection of limiting
particle trajectories.

\begin{lemma}\label{lem:gaussian-augmented-field-convolution}
Under the hypotheses of Proposition~\ref{prop:driving-field-CLT}, the joint convergence stated there extends to $(\mathbb G^N,\mathbb Z^N)$,
 with limit $(\mathbb G,\mathbb Z)$, and holds jointly with $(\mathbf U_i^{*})_{1\le i\le m}$ for every fixed $m\ge1$. The additional factors are $C([-T,0],
 \mathbf H)$ and $C(B_R^{cu},C([-T,0],\mathbf H))^m$, respectively, for every $T,R,M>0$. The fields $\mathbb G_U$, $\mathbb G_K$, $\bar W$, and $\mathbb Z$ are jointly Gaussian, and $(\mathbb G,\mathbb Z)$ is
 independent of $(\mathbf U_i^{*})_{1\le i\le m}$.
\end{lemma}

\begin{proof}
For each $j$, let $\mathcal Z_j$ be the process obtained from \eqref{eq:wiener-backward-contribution} by replacing $\bar W^N$ with 
$W_j$ and $\sqrt N\,z_{cu}^N(\omega)$ with $\sigma P_{cu}\int_{-\infty}^0e^{-\Delta r}dW_j(r)$. Since every term in \eqref{eq:wiener-backward-contribution}
 is linear in the Wiener process, $\mathbb Z^N=N^{-1/2}\sum_j\mathcal Z_j$.  Applying the Hilbert space central limit theorem \cite{AraujoGine1980}
  to the i.i.d. vectors formed by the summands defining $(\mathbb G_U^N,\mathbb G_K^N,\bar W^N,\mathbb Z^N)$ at finitely many parameter points 
  gives the joint finite dimensional convergence.
Approximation of the Wiener integrals in $L^2$, together with the It\^o isometry and \eqref{eq:driving-field-covariance}, identifies the covariance of the limit with that of $(\mathbb G,\mathbb Z)$.
For every $0<\vartheta_1<\vartheta$ and $0<\rho_1<\min\{\rho,\vartheta_1\}$, with $\rho$ as in \eqref{eq:rho-gaussian}, we claim that
\begin{equation}\label{eq:mean-noise-fluctuation-regularity-moments}
\sup_N\mathbb E\left[ \sup_{-T\le r\le0}\|\mathbb Z^N(r)\|_{V_{\vartheta_1}\times V_{\vartheta_1}}^p + [\mathbb Z^N]_{C^{\rho_1}([-T,0],\mathbf H)}^p \right]<\infty.
\end{equation}
Since $\mathbb Z^N$ has the law of $Z_1$, it suffices to estimate one Gaussian convolution. Choose $p'>p$ and $a>1/p'$ with $\vartheta_1+a<1/2$. 
Stochastic factorization and \eqref{eq:positive-fractional-semigroup-gaussian} give a finite $p'$ moment for the time supremum of the stable stochastic 
convolution in $V_{\vartheta_1}\times V_{\vartheta_1}$. For increments, choose $\rho_1<a_0<\vartheta_1$ and enlarge $p'$ so that $\rho_1<a_0-1/p'$. 
Since the stochastic integral is Gaussian, the It\^o isometry and
\eqref{eq:semigroup-increment-gaussian} imply that the $p'$th moments of the stochastic integral over $[s,t]$ and of the semigroup increment of the 
value at time $s$ are bounded by $C|t-s|^{p'/2}$ and $C|t-s|^{p'a_0}$, respectively.
The Garsia--Rodemich--Rumsey inequality \cite{GarsiaRodemichRumsey1970} yields the required $C^{\rho_1}$ moment.
 The center unstable stochastic integral is finite dimensional, and the homogeneous center term in \eqref{eq:wiener-backward-contribution} 
 is smooth on bounded intervals. This proves \eqref{eq:mean-noise-fluctuation-regularity-moments}. The compact embedding 
 $V_{\vartheta_1}\times V_{\vartheta_1}\hookrightarrow\mathbf H$ and the Arzel\`a--Ascoli theorem yield tightness of
  $(\mathbb Z^N)$ in $C([-T,0],\mathbf H)$, which, together with Proposition~\ref{prop:driving-field-CLT}, gives joint tightness.
The finite dimensional limits determine every subsequential limit, proving the joint convergence of $(\mathbb G^N,\mathbb Z^N)$ to $(\mathbb G,\mathbb Z)$.  
  
 Fix $m\ge1$ and $T,R,M>0$. For $N>m$, let $(\mathbb G^{N,m},\mathbb Z^{N,m})$ be the field obtained by
summing over $j=m+1,\ldots,N$, with normalization $(N-m)^{-1/2}$. Let $(\mathbb G^m,\mathbb Z^m)$ denote the
corresponding field formed from the first $m$ summands, with normalization $m^{-1/2}$. Then
\[
(\mathbb G^N,\mathbb Z^N)=\sqrt{\frac mN}\,(\mathbb G^m,\mathbb Z^m)+\sqrt{\frac{N-m}{N}}\,(\mathbb G^{N,m},\mathbb Z^{N,m}).
\]
Since $m$ is fixed, $\sqrt{\frac mN}\,(\mathbb G^m,\mathbb Z^m)\longrightarrow0~\text{almost surely}$
in $C([-T,0]\times B_R^{cu},H)\times C([-T,0]\times B_R^{cu}\times K_M,H) \times C([-T,0],H)\times C([-T,0],\mathbf H).$
The field $(\mathbb G^{N,m},\mathbb Z^{N,m})$ has the same law as
$(\mathbb G^{N-m},\mathbb Z^{N-m})$.  Since $N-m\to\infty$, the joint convergence proved above gives
\[
(\mathbb G^{N,m},\mathbb Z^{N,m})\Longrightarrow(\mathbb G,\mathbb Z).
\]
Moreover, $(\mathbb G^{N,m},\mathbb Z^{N,m})$ depends only on $W_{m+1},\ldots,W_N$, whereas each $\mathbf U_i^{*}$ is determined
by $W_i$. Hence, for every $N>m$, this field is independent of $(\mathbf U_i^{*})_{1\le i\le m}$, and
\[
\begin{aligned}
&\operatorname{Law}\left((\mathbb G^{N,m},\mathbb Z^{N,m}),(\mathbf U_i^{*})_{1\le i\le m}\right) \Longrightarrow\operatorname{Law}(\mathbb G,\mathbb Z)
\otimes\operatorname{Law}\bigl((\mathbf U_i^{*})_{1\le i\le m}\bigr).
\end{aligned}
\]
Finally, $\sqrt{(N-m)/N}\to1$. The decomposition above and Slutsky's theorem give the same joint limit with
$(\mathbb G^{N,m},\mathbb Z^{N,m})$ replaced by $(\mathbb G^N,\mathbb Z^N)$. This proves the joint convergence
and the asserted independence.

\end{proof}

 For fixed $\xi\in H_{cu}$ and $\mathbf y^N=(\bar y,\hat y_1,\ldots,\hat y_N)\in C_\eta^-(H^{\oplus(N+1)})$, write $y_j=\bar y+\hat y_j$ and define the 
 linearization of the drift in \eqref{eq:exact-mean-value-LP-fluctuation} by
\begin{equation}\label{eq:finite-linearized-particle-drift}
\mathcal D_{N,i}(r,\mathbf y^N) := \begin{pmatrix} \displaystyle \frac1N\sum_{j=1}^N D\nabla U(u_j^{*}(r,\xi))\,y_j(r) \\[3mm] \begin{aligned}
 &D\nabla U(u_i^{*}(r,\xi))\,y_i(r) - \frac1N\sum_{j=1}^N D\nabla U(u_j^{*}(r,\xi))\,y_j(r) \\ &\qquad+ \frac{\gamma}{N}\sum_{j=1}^N 
 D\nabla K\bigl( u_j^{*}(r,\xi)-u_i^{*}(r,\xi) \bigr) \bigl(y_j(r)-y_i(r)\bigr). \end{aligned} \end{pmatrix}.
\end{equation}
 By $(\mathbf{H}2)$--$(\mathbf{H}3)$, \eqref{eq:fixed-N-uniform-dichotomy}, and \eqref{eq:standing-LP-gap}, the integral operator generated by $\mathcal D_{N,\cdot}$ 
 has contraction constant at most $q<1$ on $C_\eta^-(H^{\oplus(N+1)})$, with norm $\sup_{r\le0}e^{-\eta r}\|\mathbf y^N(r)\|_N$.

 The next lemma gives the first order expansion of \eqref{eq:exact-mean-value-LP-fluctuation} and proves that the remainder converges to zero.

\begin{lemma}
\label{lem:finite-LP-first-order-expansion}
Assume $(\mathbf{H}1)$--$(\mathbf{H}3)$. Then, for each $N\ge2$ and $1\le i\le N$,
\begin{equation}\label{eq:exact-linearized-LP-N}
\mathbf Y_i^N = \mathscr K\left[\mathcal D_{N,i}(\mathbf Y^N)+\mathcal B_i^N\right] +\mathbb Z^N+\mathbf r_i^N,
\end{equation}
where the remainder $\mathbf r_i^N$ satisfies, for every $T,R>0$,
\begin{equation}\label{eq:gaussian-rem-L1-zero}
\lim_{N\to\infty} \sup_{|\xi|\le R} \frac1N\sum_{i=1}^N \mathbb E\|\mathbf r_i^N(\cdot,\xi)\|_{\eta,T} =0.
\end{equation}
In particular, for every fixed $i$, $\xi\in B_R^{cu}$, and $\varepsilon>0$,
\begin{equation}\label{eq:LP-Taylor-remainder-vanish}
\mathbb P\left( \|\mathbf r_i^N(\cdot,\xi)\|_{\eta,T}>\varepsilon \right) \longrightarrow0.
\end{equation}
\end{lemma}

\begin{proof}
Write $u_j^{N,*}-u_j^{*}=N^{-1/2}y_j^N$. Taylor's formula gives
\begin{equation}\label{eq:gaussian-U-Taylor}
\sqrt N\bigl[\nabla U(u_j^{N,*})-\nabla U(u_j^{*})\bigr] = D\nabla U(u_j^{*})y_j^N+\varepsilon_{U,j}^N,
\end{equation}
and, similarly,
\[\sqrt N\Bigl[ \nabla K(u_j^{N,*}-u_i^{N,*}) -\nabla K(u_j^{*}-u_i^{*}) \Bigr] = D\nabla K(u_j^{*}-u_i^{*})(y_j^N-y_i^N) 
+\varepsilon_{K,ij}^N.\]
The derivative bounds in $(\mathbf{H}2)$ imply
\[|\varepsilon_{U,j}^N| \le C\min\{N^{-1/2}|y_j^N|^2,|y_j^N|\}, \qquad |\varepsilon_{K,ij}^N| \le C\min\{N^{-1/2}|y_j^N-y_i^N|^2,|y_j^N-y_i^N|\}.\]
Set
\[
\varepsilon_i^N=
\begin{pmatrix}
\displaystyle\frac1N\sum_{j=1}^N\varepsilon_{U,j}^N\\[2mm]
\displaystyle\varepsilon_{U,i}^N-\frac1N\sum_{j=1}^N\varepsilon_{U,j}^N+\frac\gamma N\sum_{j=1}^N\varepsilon_{K,ij}^N
\end{pmatrix}.
\]
Substituting these expansions into \eqref{eq:exact-mean-value-LP-fluctuation} gives \eqref{eq:exact-linearized-LP-N} with $\mathbf r_i^N=\mathscr K\varepsilon_i^N$.
Fix $S>T$. The quadratic remainder bounds and \eqref{eq:particle-fluctuation-local-higher-moment} with $s=2$ give
\begin{equation}\label{eq:gaussian-local-Taylor-remainder-L1}
\sup_{|\xi|\le R} \frac1N\sum_i \mathbb E\|\varepsilon_i^N(\xi)\|_{\eta,S} \le C_{S,R}N^{-1/2}.
\end{equation}
 On $[-T,0]$, let $\mathscr K^{(-\infty,-S)}$ denote the stable integral in \eqref{eq:LP-convolution-operator} restricted to $(-\infty,-S]$.
For the contribution from $(-\infty,-S]$, the linear remainder bounds, \eqref{eq:particle-fluctuation-weighted-higher-moment} with $s=2$, 
Minkowski inequality, and \eqref{eq:fixed-N-uniform-dichotomy} yield
\[\sup_N\sup_{|\xi|\le R} \frac1N\sum_i \mathbb E\|\mathscr K^{(-\infty,-S)}\varepsilon_i^N(\xi)\|_{\eta,T} \le C_{T,R}e^{-(\eta-\beta)(S-T)}.\]
Letting first $N\to\infty$ and then $S\to\infty$ proves \eqref{eq:gaussian-rem-L1-zero}. Finally, exchangeability and Markov inequality 
imply \eqref{eq:LP-Taylor-remainder-vanish}.
\end{proof}

For fixed $\xi\in H_{cu}$, let $P_\xi=\operatorname{Law}\bigl(\mathbf U_1^{*}(\cdot,\xi)\bigr) \in\mathcal P\bigl(C_{\mathrm{loc}}((-\infty,0],\mathbf H)\bigr),$
where the path space carries the topology of uniform convergence on compact intervals and its Borel $\sigma$ field. The weighted moment estimate in 
Lemma~\ref{lem:backward-regularity-gaussian} gives $P_\xi(C_\eta^-(\mathbf H))=1$.
Let $\mathbb P_{\mathbb G}=\operatorname{Law}(\mathbb G,\mathbb Z)$, on the space of consistent local continuous restrictions given in 
Lemma~\ref{lem:gaussian-augmented-field-convolution}. By the independence in that lemma, the joint law with a limiting particle path is $\mathbb P_{\mathbb G}\otimes P_\xi$. 
Set $\mathbb X_{\eta,\xi}=L_\eta^\infty\bigl((-\infty,0],L^2(\mathbb P_{\mathbb G}\otimes P_\xi,\mathbf H)\bigr),$
with norm
\[
\|\mathbf Y\|_{\mathbb X_{\eta,\xi}}=\operatorname*{ess\,sup}_{r\le0}e^{-\eta r}\left(\int\!\int |\mathbf Y(r,g,z,\mathbf u)|_{\mathbf H}^2
\,P_\xi(d\mathbf u)\,\mathbb P_{\mathbb G}(dg,dz)\right)^{1/2}.
\]
The formula \eqref{eq:LP-convolution-operator} also defines $\mathscr K$ on $\mathbb X_{\eta,\xi}$, with the integrals
understood in $L^2(\mathbb P_{\mathbb G}\otimes P_\xi,\mathbf H)$.

For $\mathbf u=(\bar u,\widehat u)\in C_\eta^-(\mathbf H)$, write $u=\bar u+\widehat u$. Define
\begin{equation}\label{eq:driving-particle-forcing-limit}
\mathcal B(r,\xi,\mathbf u,\mathbb G) := \begin{pmatrix} \mathbb G_U(r,\xi) \\[1mm] -\mathbb G_U(r,\xi) + \gamma\mathbb G_K(r,\xi,u(r)) \end{pmatrix}.
\end{equation}
The value $\mathbb G_K(r,\xi,u(r))$ is defined by the consistent local versions fixed after Proposition~\ref{prop:driving-field-CLT}. 
By Lemma~\ref{lem:backward-regularity-gaussian}, this expression is well defined for $P_\xi$-almost every particle path.

Let $\mathbf U$ denote the particle coordinate with law $P_\xi$. For $\mathbf Y\in\mathbb X_{\eta,\xi}$, let $\widetilde{\mathbf U}$ be an independent 
copy of $\mathbf U$, with $(\mathbb G,\mathbb Z)$ fixed. Write
$\mathbf U=(\bar u,\widehat u), ~
\widetilde{\mathbf U}=(\widetilde{\bar u},\widetilde{\widehat u}),
\mathbf Y(\mathbb G,\mathbb Z,\mathbf U)=(\bar y,\widehat y), ~
\mathbf Y(\mathbb G,\mathbb Z,\widetilde{\mathbf U})=(\widetilde{\bar y},\widetilde{\widehat y}),$
and set $u=\bar u+\widehat u$, $\widetilde u=\widetilde{\bar u}+\widetilde{\widehat u}$, $y=\bar y+\widehat y$, and $\widetilde y
=\widetilde{\bar y}+\widetilde{\widehat y}$.
 Let $\widetilde{\mathbb E}$ denote expectation with respect to $\widetilde{\mathbf U}$ only, with $(\mathbb G,\mathbb Z)$ fixed. Define
\begin{equation}\label{eq:limit-linearized-LP-operator}
\mathscr L_\xi\mathbf Y :=\mathscr K\!\left[r\mapsto \begin{pmatrix} \displaystyle \widetilde{\mathbb E}\left[D\nabla U(\widetilde u(r))\widetilde y(r)\right]\\[3mm]
 \begin{aligned} 
&D\nabla U(u(r))y(r)-\widetilde{\mathbb E}\left[D\nabla U(\widetilde u(r))\widetilde y(r)\right]\\ 
&\qquad+\gamma\widetilde{\mathbb E}\left[D\nabla K(\widetilde u(r)-u(r))(\widetilde y(r)-y(r))\right] \end{aligned}
 \end{pmatrix}\right].
\end{equation}
The limiting fluctuation equation is
\begin{equation}\label{eq:gaussian-linear-fluctuation-equation}
\mathbf Y=\mathscr L_\xi\mathbf Y+\mathscr K\mathcal B(\cdot,\xi,\mathbf U,\mathbb G)+\mathbb Z.
\end{equation}

\begin{proposition}\label{prop:joint-backward-fluctuation}
Assume $(\mathbf{H}1)$--$(\mathbf{H}3)$, and let $\mathbb G^N$ and $\mathbb G$ be defined by
\eqref{eq:driving-field-triples}.
For each $\xi\in H_{cu}$, equation \eqref{eq:gaussian-linear-fluctuation-equation}
has a unique solution in $\mathbb X_{\eta,\xi}$.
Let $\mathbf Y_i$ denote this solution with $\mathbf U=\mathbf U_i^{*}$.
For every fixed $m\ge1$ and $T,R,M>0$, we have the joint convergence
\[
\begin{aligned}
\left(\mathbb G_U^N,\mathbb G_K^N,\bar W^N,\mathbb Z^N,(\mathbf U_i^{*})_{1\le i\le m},(\mathbf Y_i^N)_{1\le i\le m}\right)\Longrightarrow
\left(\mathbb G_U,\mathbb G_K,\bar W,\mathbb Z,(\mathbf U_i^{*})_{1\le i\le m},(\mathbf Y_i)_{1\le i\le m}\right)
\end{aligned}
\]
in
\[
\begin{aligned}
&C([-T,0]\times B_R^{cu},H)\times C([-T,0]\times B_R^{cu}\times K_M,H)\\
&\quad\times C([-T,0],H)\times C([-T,0],\mathbf H)\times C\bigl(B_R^{cu},C([-T,0],\mathbf H)\bigr)^{2m}.
\end{aligned}
\]
The pair $(\mathbb G,\mathbb Z)$ is independent of $(\mathbf U_i^{*})_{1\le i\le m}$.
Conditionally on $(\mathbb G,\mathbb Z)$, $(\mathbf U_i^{*},\mathbf Y_i)_{1\le i\le m}$ are independent and identically distributed.
\end{proposition}

\begin{proof}
\smallskip\noindent\emph{Step 1. The limiting linear equation.} By \eqref{eq:limit-linearized-LP-operator}, Jensen inequality, $(\mathbf{H}2)$--$(\mathbf{H}3)$, and \eqref{eq:standing-LP-gap},
\begin{equation}\label{eq:gaussian-linearized-LP-operator-bound}
\|\mathscr L_\xi\|_{\mathcal L(\mathbb X_{\eta,\xi})}
\le \ell_{\mathrm{LP}}\left(\frac1{\alpha-\eta}+\frac1{\eta-\beta}\right)=q<1.
\end{equation}
Proposition~\ref{prop:driving-field-CLT}, the independence in Lemma~\ref{lem:gaussian-augmented-field-convolution}, and \eqref{eq:backward-global-weighted-moment} 
give the weighted $L^2$ bound for $\mathcal B$. The weighted $L^2$ bound for $\mathbb Z$ follows from \eqref{eq:averaged-noise-weighted-theta-bound}, 
Lemma~\ref{lem:gaussian-augmented-field-convolution}, and Fatou's lemma. Therefore, 
$\mathscr K\mathcal B(\cdot,\xi,\mathbf U,\mathbb G)+\mathbb Z\in\mathbb X_{\eta,\xi}$.
Set $F_\xi:=\mathscr K\mathcal B(\cdot,\xi,\mathbf U,\mathbb G)+\mathbb Z$.
Since $\|\mathscr L_\xi\|\le q<1$, \eqref{eq:gaussian-linear-fluctuation-equation} has a unique solution in $\mathbb X_{\eta,\xi}$, given by
\begin{equation}\label{eq:gaussian-linear-fluctuation-resolvent}
\begin{aligned}
\mathbf Y=(I-\mathscr L_\xi)^{-1}F_\xi=\sum_{k=0}^{\infty}\mathscr L_\xi^k F_\xi,\qquad
\|(I-\mathscr L_\xi)^{-1}\|&\le(1-q)^{-1}.
\end{aligned}
\end{equation}
The series converges in $\mathbb X_{\eta,\xi}$ because
$\|\mathscr L_\xi^kF_\xi\|_{\mathbb X_{\eta,\xi}}\le q^k\|F_\xi\|_{\mathbb X_{\eta,\xi}}$.
For every $T>0$, the semigroup bounds also give
\[
\|\mathscr L_\xi W\|_{L^2(\mathbb P_{\mathbb G}\otimes P_\xi;C([-T,0],\mathbf H))}
\le C_T\|W\|_{\mathbb X_{\eta,\xi}}.
\]
The same bounds for $\mathscr K\mathcal B$, together with the local path moment bound for $\mathbb Z$, show that
$F_\xi\in L^2(\mathbb P_{\mathbb G}\otimes P_\xi;C([-T,0],\mathbf H))$.
For $k\ge1$, it follows that
\[
\|\mathscr L_\xi^kF_\xi\|_{L^2(\mathbb P_{\mathbb G}\otimes P_\xi;C([-T,0],\mathbf H))}
\le C_Tq^{k-1}\|F_\xi\|_{\mathbb X_{\eta,\xi}}.
\]
Thus the series also converges in this local path space and almost surely uniformly on each compact time interval.
The series defines a measurable solution map $\Phi_\xi(\mathbb G,\mathbb Z,\mathbf u)$ with continuous paths.

\smallskip\noindent\emph{Step 2. Tightness for a fixed center parameter.} Fix $\xi\in H_{cu}$ and define the empirical measure on the local path space by
\[
\Lambda_N=\frac1N\sum_{j=1}^N\delta_{(\mathbf U_j^{*},\mathbf Y_j^N)}\in\mathcal P\bigl(C_{\mathrm{loc}}(( -\infty,0],\mathbf H)^2\bigr).
\]
For $T>0$, define $\Lambda_N^T:=\frac1N\sum_{j=1}^N\delta_{\left(\mathbf U_j^{*}|_{[-T,0]},\mathbf Y_j^N|_{[-T,0]}\right)}\in\mathcal P\bigl(C([-T,0],\mathbf H)^2\bigr).$
Lemma~\ref{lem:gaussian-augmented-field-convolution} gives tightness of $(\mathbb G^N,\mathbb Z^N)$. Choose $S=T+1$, $0<\vartheta_1<\vartheta$, 
and $0<\rho_1<\min\{\rho,\vartheta_1\}$. Lemma~\ref{lem:backward-regularity-gaussian} gives the estimates for $\mathbf U_i^{*}$ below. For $\mathbf Y_i^N$, 
use \eqref{eq:exact-mean-value-LP-fluctuation}, the local moment bounds \eqref{eq:GU-local-pmoment}, \eqref{eq:particle-GK-random-state-pmoment}, 
and \eqref{eq:particle-fluctuation-local-higher-moment}, together with the semigroup estimates
 \eqref{eq:positive-fractional-semigroup-gaussian}--\eqref{eq:semigroup-increment-gaussian} and \eqref{eq:mean-noise-fluctuation-regularity-moments}. These give
\begin{equation}\label{eq:fluctuation-compactness-functional}
\begin{aligned}
 \sup_N\frac1N\sum_{i=1}^N\mathbb E\Bigl[ &\sup_{-T\le r\le0}\|\mathbf U_i^{*}(r)\|_{V_{\vartheta_1}\times V_{\vartheta_1}}^2 +[\mathbf U_i^{*}]_{C^{\rho_1}([-T,0],\mathbf H)}^2\\ 
 &+\sup_{-T\le r\le0}\|\mathbf Y_i^N(r)\|_{V_{\vartheta_1}\times V_{\vartheta_1}}^2 +[\mathbf Y_i^N]_{C^{\rho_1}([-T,0],\mathbf H)}^2 \Bigr]<\infty. 
 \end{aligned}
\end{equation}
Since $V_{\vartheta_1}\hookrightarrow H$ compactly, \eqref{eq:fluctuation-compactness-functional},
 the Arzel\`a--Ascoli theorem, Markov inequality, and Prokhorov theorem imply tightness of $(\Lambda_N^T)_N$ in $\mathcal P(C([-T,0],\mathbf H)^2)$, 
 see also \cite[Chapter~16]{Kallenberg2002}.  By exchangeability and
\eqref{eq:fluctuation-compactness-functional},
for every fixed $m\ge1$, the sequence
$\bigl((\mathbf U_i^{*},\mathbf Y_i^N)_{1\le i\le m}\bigr)_{N\ge\max\{2,m\}}$
is tight in $C([-T,0],\mathbf H)^{2m}$.

\smallskip\noindent\emph{Step 3. Identification of the subsequential limits.} By Step~2 and a diagonal argument over $T$, 
choose a subsequence such that, for every fixed $m$,
\begin{equation}\label{eq:gaussian-subseq-fluctuation-law}
\begin{aligned} &(\mathbb G^N,\mathbb Z^N,\Lambda_N, (\mathbf U_i^{*},\mathbf Y_i^N)_{1\le i\le m}) 
\Longrightarrow (\mathbb G,\mathbb Z,\Lambda, (\mathbf U_i^{*},\mathbf Y_i)_{1\le i\le m}). \end{aligned}
\end{equation}
The convergence is understood in the product of the local path topologies of Lemma~\ref{lem:gaussian-augmented-field-convolution},
 the weak topology on $\mathcal P\bigl(C_{\mathrm{loc}}(( -\infty,0],\mathbf H)^2\bigr)$, and the local path topology for the particle pairs.
The first marginal of $\Lambda$ is $P_\xi$, by the law of large numbers for the empirical measure of the i.i.d. paths $(\mathbf U_j^{*})_{j\ge1}$. 
Since the joint laws are invariant under permutations of the particle labels, the exchangeability
 result for empirical measures identifies the limiting particle pairs as conditionally i.i.d. with directing measure $\Lambda$, given 
 $(\mathbb G,\mathbb Z,\Lambda)$, see Sznitman~\cite[Proposition~2.2 and Remark~2.3]{Sznitman1991} and 
 Diaconis--Freedman~\cite{DiaconisFreedman1980}. Lemma~\ref{lem:gaussian-augmented-field-convolution} 
 implies that $(\mathbf U_i^{*})_{1\le i\le m}$ is independent of $(\mathbb G,\mathbb Z)$.
For a probability measure $\mu$ on $C_{\mathrm{loc}}((-\infty,0],\mathbf H)^2$ with finite local second moments in the second coordinate, 
write $\mathbf u=(\bar u,\hat u)$, $\mathbf y=(\bar y,\hat y)$,
$u=\bar u+\hat u$, and $y=\bar y+\hat y$. For the integration variables $\mathbf v=(\bar v,\hat v)$ and $\mathbf z=(\bar z,\hat z)$,
set $v=\bar v+\hat v$ and $z=\bar z+\hat z$. Define
\begin{equation}\label{eq:limit-linearized-particle-drift}
\mathcal D^\mu(r,\mathbf u,\mathbf y)=
\begin{pmatrix}
\displaystyle\int D\nabla U(v(r))z(r)\,\mu(d\mathbf v,d\mathbf z)\\[2mm]
\begin{aligned}
&D\nabla U(u(r))y(r)-\int D\nabla U(v(r))z(r)\,\mu(d\mathbf v,d\mathbf z)\\
&\qquad+\gamma\int D\nabla K(v(r)-u(r))(z(r)-y(r))\,\mu(d\mathbf v,d\mathbf z)
\end{aligned}
\end{pmatrix}.
\end{equation}
Then $\mathcal D^{\Lambda_N}(r,\mathbf U_i^{*},\mathbf Y_i^N)=\mathcal D_{N,i}(r,\mathbf Y^N)$.
Fix $S>T$. On $[-S,0]$, truncate the particle fluctuations $\mathbf Y_j^N$ in $\mathcal D^{\Lambda_N}$ by a continuous cutoff and localize the state argument of $\mathbb G_K^N$ to $K_M$.
The truncated drift integrands are bounded and continuous, and the localized field evaluation is continuous, hence \eqref{eq:gaussian-subseq-fluctuation-law} applies to the integrals over $[-S,0]$.
By \eqref{eq:particle-fluctuation-local-higher-moment} and Lemma~\ref{lem:backward-regularity-gaussian}, we remove the truncation of $\mathbf Y_j^N$ and let $M\to\infty$.
Moreover, \eqref{eq:LP-Taylor-remainder-vanish} shows that the Taylor remainder vanishes.
By \eqref{eq:fixed-N-uniform-dichotomy}, \eqref{eq:particle-fluctuation-weighted-higher-moment}, and
\eqref{eq:weighted-driving-field-parameter-bound}, the integral over $(-\infty,-S]$ satisfies
\[
\begin{aligned}
&\sup_{N\ge\max\{2,i\}}
\left\|t\mapsto\int_{-\infty}^{-S}e^{\mathbf L_s(t-r)}\mathbf P_s
\bigl[\mathcal D_{N,i}(r,\mathbf Y^N)+\mathcal B_i^N(r,\xi)\bigr]\,dr\right\|_{L^2(\Omega;C([-T,0],\mathbf H))}\\
&\qquad\le C_{T,\xi}e^{-(\eta-\beta)(S-T)}.
\end{aligned}
\]
The same estimate holds for the limiting integral by the moment bounds inherited through lower semicontinuity. Letting $S\to\infty$ gives
\begin{equation}\label{eq:gaussian-limit-fluctuation-with-Lambda}
\mathbf Y =\mathscr K\!\left[ \mathcal D^\Lambda(\cdot,\mathbf U,\mathbf Y) +\mathcal B(\cdot,\xi,\mathbf U,\mathbb G) \right]+\mathbb Z.
\end{equation}
This identity holds almost surely for $\Lambda$-almost every $(\mathbf U,\mathbf Y)$.
Let $\Phi_\xi$ be the solution map from Step~1. Since the first marginal of $\Lambda$ is $P_\xi$, the weighted moment bounds give
\[
D:=\sup_{r\le0}e^{-\eta r}\left(\mathbb E\int\left|\mathbf y(r)-\Phi_\xi(\mathbb G,\mathbb Z,\mathbf u)(r)\right|_{\mathbf H}^2
\,\Lambda(d\mathbf u,d\mathbf y)\right)^{1/2}<\infty.
\]
Using the first marginal identity, write the particle averages in the equation for $\Phi_\xi$ as integrals with respect to $\Lambda$.
Comparing with \eqref{eq:gaussian-limit-fluctuation-with-Lambda} and applying the same estimates as in Step~1 gives
\[
D\le\ell_{\mathrm{LP}}\left(\frac1{\alpha-\eta}+\frac1{\eta-\beta}\right)D=qD.
\]
Hence $D=0$. Equality at every rational time and continuity of both paths then give
\[
\Lambda(d\mathbf u,d\mathbf y)
=P_\xi(d\mathbf u)\,\delta_{\Phi_\xi(\mathbb G,\mathbb Z,\mathbf u)}(d\mathbf y)
\quad\text{a.s.}
\]
This identifies $\Lambda$ as a measurable function of $(\mathbb G,\mathbb Z)$. Together with the conditional i.i.d. property established
above, it shows that all subsequential limits have the same joint law. Hence the whole sequence converges. Therefore, for every fixed $m$ and $\xi$,
\begin{equation}\label{eq:fixed-center-joint-particle-limit}
\left( \mathbb G^N,\mathbb Z^N, (\mathbf U_i^{*},\mathbf Y_i^N)_{1\le i\le m} \right) \Longrightarrow \left( \mathbb G,\mathbb Z, (\mathbf U_i^{*},\mathbf Y_i)_{1\le i\le m} \right),
\qquad N \rightarrow \infty .
\end{equation}

\smallskip\noindent\emph{Step 4. Uniformity with respect to the center parameter.} Fix $T,R>0$. For every $S>0$ and $|\xi|,|\zeta|\le R$, 
\eqref{eq:fixed-N-particle-center-Lipschitz} and \eqref{eq:backward-xi-weighted-Lipschitz} give
\begin{equation}\label{eq:finite-limit-backward-xi-uniform-Lipschitz}
\sup_{N,i}\|\mathbf U_i^{N,*}(\xi)-\mathbf U_i^{N,*}(\zeta)\|_{\eta,S} +\sup_i\|\mathbf U_i^{*}(\xi)-\mathbf U_i^{*}(\zeta)\|_{\eta,S} \le C_{S,R}|\xi-\zeta|.
\end{equation}
Subtracting \eqref{eq:exact-mean-value-LP-fluctuation} at $\xi$ and $\zeta$ eliminates $\mathbb Z^N$. Split the drift difference as
\[
\mathfrak D_N^\xi\mathbf Y^N(\cdot,\xi)-\mathfrak D_N^\zeta\mathbf Y^N(\cdot,\zeta)
=\mathfrak D_N^\xi(\mathbf Y^N(\cdot,\xi)-\mathbf Y^N(\cdot,\zeta))+(\mathfrak D_N^\xi-\mathfrak D_N^\zeta)\mathbf Y^N(\cdot,\zeta),
\]
where $\mathfrak D_N^\xi=(\mathfrak D_{N,i}^\xi)_{i\le N}$.
Choose $\eta_-\in(\beta,\eta)$ close enough to $\eta$ that
$
q_-:=\ell_{\mathrm{LP}}\left(\frac1{\alpha-\eta_-}+\frac1{\eta_--\beta}\right)<1.
$
Set $\Pi_\delta=\{(\xi,\zeta)\in(B_R^{cu})^2:|\xi-\zeta|\le\delta\}$ and
\[
D_N(\delta):=\sup_{r\le0}e^{-\eta_-r}\left[\frac1N\sum_i\mathbb E\sup_{(\xi,\zeta)\in\Pi_\delta}
|\mathbf Y_i^N(r,\xi)-\mathbf Y_i^N(r,\zeta)|^2\right]^{1/2}.
\]
For fixed $S>T$, $(\mathbf{H}2)$, \eqref{eq:finite-limit-backward-xi-uniform-Lipschitz}, and \eqref{eq:particle-fluctuation-local-higher-moment} control the coefficient error 
on $[-S,0]$. Tightness in Proposition~\ref{prop:driving-field-CLT} and the $p$th moment bound in Lemma~\ref{lem:driving-field-local-pmoment} give local 
$L^2$ equicontinuity of $(\mathbb G_U^N)_N$. Together with \eqref{eq:particle-GK-Lp-equicontinuity}, these bounds control the coefficient error and forcing 
difference by $\varepsilon_S(\delta)\to0$, uniformly in $N$, with the supremum over $\Pi_\delta$ inside the expectation.
On $r<-S$, boundedness of the coefficients, \eqref{eq:weighted-driving-field-parameter-bound}, 
and \eqref{eq:particle-fluctuation-weighted-higher-moment} give a bound $Ce^{-(\eta-\eta_-)S}$ in the weaker weighted norm. Hence
\[
D_N(\delta)\le q_-D_N(\delta)+C_S\varepsilon_S(\delta)+Ce^{-(\eta-\eta_-)S}.
\]
The norms $D_N(\delta)$ are uniformly finite by \eqref{eq:particle-fluctuation-weighted-higher-moment}. Absorb the first term, let $\delta\downarrow0$ at fixed $S$, and then let $S\to\infty$. 
Returning to the difference equation on $[-T,0]$, with the same tail decomposition, gives
\begin{equation}\label{eq:fluctuation-xi-stochastic-equicontinuity}
\lim_{\delta\downarrow0}\sup_N\frac1N\sum_{i=1}^N\mathbb E \sup_{\substack{|\xi|,|\zeta|\le R\\|\xi-\zeta|\le\delta}} \|\mathbf Y_i^N(\cdot,\xi)-\mathbf Y_i^N(\cdot,\zeta)\|_{\eta,T}^2=0.
\end{equation}
By exchangeability, \eqref{eq:fluctuation-xi-stochastic-equicontinuity} also holds with the empirical average replaced by any fixed particle index. 
Step~3 applies simultaneously to every finite collection of center parameters. Therefore, \eqref{eq:fixed-center-joint-particle-limit},
 \eqref{eq:fluctuation-xi-stochastic-equicontinuity}, and the tightness criterion for continuous processes on compact
  parameter sets \cite[Theorem~16.5]{Kallenberg2002} give the joint convergence
\begin{align}
& \left( \mathbb G^N, \mathbb Z^N, (\mathbf U_i^{*}(\cdot,\cdot))_{1\le i\le m}, (\mathbf Y_i^N(\cdot,\cdot))_{1\le i\le m} \right) 
\Longrightarrow 
\left( \mathbb G, \mathbb Z, (\mathbf U_i^{*}(\cdot,\cdot))_{1\le i\le m}, (\mathbf Y_i(\cdot,\cdot))_{1\le i\le m} \right) \label{eq:joint-particle-LP-fluctuation-uniform}
\end{align}
in $C([-T,0]\times B_R^{cu},H)\times C([-T,0]\times B_R^{cu}\times K_M,H)\times C([-T,0],H)\times C([-T,0],\mathbf H)\times C\bigl(B_R^{cu},C([-T,0],\mathbf H)\bigr)^{2m}.$ 

\end{proof}

\begin{proof}[Proof of Theorem~\ref{thm:intro-LP-fluctuation}]
Proposition~\ref{prop:joint-backward-fluctuation} gives
\eqref{eq:joint-particle-LP-fluctuation-uniform}.
Taking the particle fluctuation marginal proves \eqref{eq:intro-joint-particle-LP-fluctuation}.
For fixed $\xi$, the operator $\mathscr L_\xi$ acts only on the particle variable and time. Thus, conditionally on $\mathbf U_i^{*}(\cdot,\xi)$,
 the partial sums of the series in \eqref{eq:gaussian-linear-fluctuation-resolvent} are linear functions of the independent Gaussian field 
 $(\mathbb G,\mathbb Z)$. Passing to the $L^2$ limit at each finite set of times proves the conditional Gaussianity.
The limits on $B_R^{cu}$, $R>0$, are consistent under restriction by the joint construction of the driving fields and uniqueness of
\eqref{eq:gaussian-linear-fluctuation-equation}. Since $H_{cu}$ is finite dimensional, \cite[Proposition~16.6]{Kallenberg2002} gives convergence in
$C_{\mathrm{loc}}(H_{cu},C([-T,0],\mathbf H))^m$. 
 
 \end{proof}

For $t\in\mathbb R$ and $r\le0$, set $\mathbb Z_\theta^N(t,r)=\mathbb Z^N(r,\theta_t\omega)$.
The limiting field $\mathbb Z_\theta(t,\cdot)$ is obtained from the formula for $\mathbb Z$ by replacing $\bar W(s)$ with
$\bar W(s+t)-\bar W(t)$.

\begin{corollary}\label{cor:theta-LP-Gaussian-fluctuation}
Under the assumptions of Theorem~\ref{thm:intro-LP-fluctuation}, let $T_0,T,R,M>0$ and fix $m\ge1$. There exist versions of the limiting fields satisfying
$\bigl((t,\xi)\mapsto\mathbf Y_{i,\theta}(t,\cdot,\xi)\bigr)\in C\bigl([0,T_0]\times B_R^{cu}, C([-T,0],\mathbf H)\bigr),~ 1\le i\le m,$
such that
\begin{align}\label{eq:theta-LP-fluctuation-joint}
&\Bigl(\mathbb G_{U,\theta}^N,\mathbb G_{K,\theta}^N,\bar W^N,\mathbb Z_\theta^N,\bigl((t,\xi)\mapsto
\mathbf U_i^{*}(\cdot,\xi,\theta_t\omega)\bigr)_{1\le i\le m},\bigl((t,\xi)\mapsto
\mathbf Y_i^N(\cdot,\xi,\theta_t\omega)\bigr)_{1\le i\le m}\Bigr)\notag\\
&\qquad\Longrightarrow
\Bigl(\mathbb G_{U,\theta},\mathbb G_{K,\theta},\bar W,\mathbb Z_\theta,\bigl((t,\xi)\mapsto\mathbf U_i^{*}(\cdot,\xi,\theta_t\omega)\bigr)_{1\le i\le m},
\bigl((t,\xi)\mapsto\mathbf Y_{i,\theta}(t,\cdot,\xi)\bigr)_{1\le i\le m}\Bigr).
\end{align}
in $C([0,T_0]\times[-T,0]\times B_R^{cu},H)\times C([0,T_0]\times[-T,0]\times B_R^{cu}\times K_M,H)\times C([0,T_0],H)\\
\times C([0,T_0],C([-T,0],\mathbf H)) \times C\bigl([0,T_0]\times B_R^{cu},C([-T,0],\mathbf H)\bigr)^{2m}$.

$(\mathbb G_{U,\theta},\allowbreak \mathbb G_{K,\theta},\allowbreak \bar W,\allowbreak \mathbb Z_\theta)$ is independent of 
$(\mathbf U_i^{*}(\cdot,\cdot,\theta_\cdot\omega))_{1\le i\le m}$.
Conditionally on $(\mathbb G_{U,\theta},\allowbreak \mathbb G_{K,\theta},\allowbreak \bar W,\allowbreak \mathbb Z_\theta)$, the pairs
 $(\mathbf U_i^{*}(\cdot,\cdot,\theta_\cdot\omega),\mathbf Y_{i,\theta})$, $1\le i\le m$, are i.i.d. For each $(t,\xi)$ and $1\le i\le m$,
$\mathbf Y_{i,\theta}(t,\cdot,\xi)$ is the unique solution of \eqref{eq:gaussian-linear-fluctuation-equation} in the weighted $L^2$ space obtained from $\mathbb X_{\eta,\xi}$ by
 replacing $\mathbb P_{\mathbb G}$ with the joint law of $(\mathbb G_{U,\theta},\allowbreak \mathbb G_{K,\theta},\allowbreak \bar W,\allowbreak \mathbb Z_\theta)$.
  It is obtained by taking
$\mathbf U=\mathbf U_i^{*}(\cdot,\xi,\theta_t\omega)$ and replacing $\mathbb G_U(r,\xi)$, $\mathbb G_K(r,\xi,z)$, and $\mathbb Z(r)$ by $\mathbb G_{U,\theta}(t,r,\xi)$, 
$\mathbb G_{K,\theta}(t,r,\xi,z)$, and $\mathbb Z_\theta(t,r)$, respectively.  Moreover,
\[
\bigl((t,\xi)\mapsto\mathbf Y_i^N(0,\xi,\theta_t\omega)\bigr)_{1\le i\le m}\Longrightarrow\bigl((t,\xi)\mapsto\mathbf Y_{i,\theta}(t,0,\xi)\bigr)_{1\le i\le m}
\]
in $C([0,T_0]\times B_R^{cu},\mathbf H)^m$.
\end{corollary}

\begin{proof}

With $Z_j$ as in the proof of Lemma~\ref{lem:gaussian-augmented-field-convolution},
$\mathbb Z_\theta^N(t,r)=N^{-1/2}\sum_{j=1}^N Z_j(r,\theta_t\omega)$.
The Hilbert space central limit theorem applied to the vectors formed by the summands defining 
$(\mathbb G_{U,\theta}^N,\allowbreak \mathbb G_{K,\theta}^N,\allowbreak \bar W^N,\allowbreak \mathbb Z_\theta^N)$ at finitely many parameter points gives joint finite dimensional convergence. 
The law of $\mathbb Z_\theta^N$ does not depend on $N$, and its paths are continuous in $(t,r)$; hence these fields are tight. Together
 with \eqref{eq:theta-driving-field-CLT}, this gives joint convergence of the first four fields in \eqref{eq:theta-LP-fluctuation-joint}. 
 Removing the first $m$ summands from every component, as in Lemma~\ref{lem:gaussian-augmented-field-convolution}, gives
\[
\left( \mathbb G_{U,\theta}^N,\mathbb G_{K,\theta}^N,\bar W^N,\mathbb Z_\theta^N, (\mathbf U_i^{*}(\cdot,\cdot,\theta_\cdot\omega))_{1\le i\le m} \right)
\Longrightarrow
\left( \mathbb G_{U,\theta},\mathbb G_{K,\theta},\bar W,\mathbb Z_\theta, (\mathbf U_i^{*}(\cdot,\cdot,\theta_\cdot\omega))_{1\le i\le m} \right)
\]
in the first five factors of the product space stated in the corollary, and proves independence of 
$(\mathbb G_{U,\theta},\allowbreak \mathbb G_{K,\theta},\allowbreak \bar W,\allowbreak \mathbb Z_\theta)$ from the fixed particle trajectories.
For every finite collection $(t_k,\xi_k)_{1\le k\le\ell}$, the invariance of $\mathbb P$ under $\theta_{t_k}$ and \eqref{eq:LP-Taylor-remainder-vanish} show that all the 
corresponding Taylor remainders converge to zero in probability. Apply Step~3 in the proof of Proposition~\ref{prop:joint-backward-fluctuation} simultaneously to 
\eqref{eq:exact-linearized-LP-N} at these parameter values. This identifies the finite dimensional distributions in \eqref{eq:theta-LP-fluctuation-joint}.

It remains to prove tightness with respect to $(t,\xi)$. Choose $1/p<\rho_0<\rho<\vartheta-1/p$ and set
\[d\bigl((t,\xi),(v,\zeta)\bigr) :=|t-v|^{\rho_0}+|\xi-\zeta|.\]
Compare \eqref{eq:LP-orbit-fixed-N} under $\theta_t$ and $\theta_v$. First use \eqref{eq:FN-uniform-Lipschitz} in the empirical norm, with contraction constant $q$, 
and then \eqref{eq:FNi-particle-Lipschitz} for each component, with contraction constant $q_s$. The increment bounds for the Ornstein--Uhlenbeck processes, 
as in Lemma~\ref{lem:backward-theta-regularity-gaussian}, give
\[
\sup_{N\ge2}\max_{i\le N}\mathbb E\sup_{|\xi|\le R}\bigl[t\mapsto\mathbf U_i^{N,*}(\cdot,\xi,\theta_t\omega)\bigr]_{C^{\rho_0}([0,T_0],C_\eta^-)}^p<\infty.
\]
The constants are uniform because the empirical norms of the Ornstein--Uhlenbeck increments have uniformly bounded moments. 
Together with $p>4$ and \eqref{eq:particle-fluctuation-local-higher-moment}, this controls the coefficient error in $L^2$ by H\"older's inequality. 
The additional noise difference $\mathbb Z_\theta^N(t)-\mathbb Z_\theta^N(v)$ tends to zero locally in $L^2$ as $t-v\to0$, uniformly in $N$, 
since these fields have an $N$ independent law, continuous paths, and the moment bounds in \eqref{eq:averaged-noise-local-theta-pmoment}.
 Their tails in the weaker weight are controlled by \eqref{eq:averaged-noise-weighted-theta-bound}.
Tightness in \eqref{eq:theta-driving-field-CLT} and the $p$th moment bound in Lemma~\ref{lem:driving-field-local-pmoment} give local
 $L^2$ equicontinuity of $(\mathbb G_{U,\theta}^N)_N$. Repeating Step~4 in the proof of Proposition~\ref{prop:joint-backward-fluctuation} with the metric $d$, and using 
Lemma~\ref{lem:backward-theta-regularity-gaussian}, \eqref{eq:theta-particle-GK-Lp-equicontinuity}, \eqref{eq:particle-fluctuation-local-higher-moment}, 
and the preceding noise estimates, gives, for every fixed $i$,
\begin{equation}\label{eq:theta-fluctuation-equicontinuity}
\lim_{\delta\downarrow0}\sup_N\mathbb E
\sup_{\substack{0\le t,v\le T_0,\ |\xi|,|\zeta|\le R\\
 d((t,\xi),(v,\zeta))\le\delta}}
\|\mathbf Y_i^N(\cdot,\xi,\theta_t\omega)
-\mathbf Y_i^N(\cdot,\zeta,\theta_v\omega)\|_{\eta,T}^2
=0.
\end{equation}
The convergence of the first five components, the finite dimensional convergence of the particle fields, \eqref{eq:theta-fluctuation-equicontinuity}, 
and the tightness criterion for continuous processes on compact parameter sets \cite[Theorem~16.5]{Kallenberg2002} give
 \eqref{eq:theta-LP-fluctuation-joint}. Independence of $(\mathbb G_{U,\theta},\allowbreak \mathbb G_{K,\theta},\allowbreak \bar W,\allowbreak \mathbb Z_\theta)$
  from the limiting particle fields and uniqueness of the weighted $L^2$ solution of 
 \eqref{eq:gaussian-linear-fluctuation-equation} give the conditional i.i.d. statement. The final convergence follows from the continuity of the map $f\mapsto f(0)$.
\end{proof}

\section{Gaussian approximation of the reduced dynamics}
\label{sec:gaussian-reduced-dynamics}

\subsection{Dynamics in the original variables}

Let $Y^N(t)$ be a trajectory of \eqref{eq:vector-particle-system} on $\mathcal M^N(\theta_t\omega)+z^N(\theta_t\omega)$, 
and write $V^N(t)=Y^N(t)-z^N(\theta_t\omega)$. Set $X_t^N=P_{cu}\bar u^N(t)=V_{cu}^N(t)+z_{cu}^N(\theta_t\omega)$, 
where $V_{cu}^N$ is the center component of $V^N$. Then
\begin{equation}\label{eq:phy-par-red}
u_i^N(t)=u_i^{N,*}\bigl(0,X_t^N-z_{cu}^N(\theta_t\omega),\theta_t\omega\bigr),
\end{equation}
and projection of \eqref{eq:bar-u-mean-fluctuation-system} gives
\begin{equation}\label{eq:finite-N-physical-reduced-equation}
\begin{aligned}
dX_t^N={}&\left[L_{cu}X_t^N+\frac1N\sum_{i=1}^N
P_{cu}\nabla U\bigl(u_i^{N,*}\bigl(0,X_t^N-z_{cu}^N(\theta_t\omega),\theta_t\omega\bigr)\bigr)\right]dt+\frac{\sigma}{\sqrt N}P_{cu}\,d\bar W^N(t,\omega).
\end{aligned}
\end{equation}
Let $X$ solve
\begin{equation}\label{eq:physical-limit-reduced-equation}
 \frac{d}{dt}X_t = L_{cu}X_t +\mathbb E\left[ P_{cu}\nabla U\bigl(u_1^{*}(0,X_t,\omega)\bigr) \right], \qquad X(0)=X_0,
\end{equation}
with  initial value $X_0\in H_{cu}$. By \eqref{eq:OU-mean-fluctuation-components}, for every $T>0$,
\begin{equation}\label{eq:center-OU-N-rate}
\mathbb E\sup_{0\le t\le T}|z_{cu}^N(\theta_t\omega)|^2\le \frac{C_T}{N}.
\end{equation}
In this subsection, for random variables $X_N$ with values in a normed space $E$, we write $X_N=O_{\mathbb P}(1)$ if
\[
\lim_{M\to\infty}\sup_{N\ge1}
\mathbb P\bigl(\|X_N\|_E>M\bigr)=0,
\]
and $X_N=o_{\mathbb P}(1)$ if, for every $\varepsilon>0$,
\[
\lim_{N\to\infty}\mathbb P\bigl(\|X_N\|_E>\varepsilon\bigr)=0.
\]

The next proposition gives a mean-square estimate for the difference between the finite particle and limiting reduced trajectories.

\begin{proposition}
\label{prop:reduced-mean-square-rate}
Assume $(\mathbf{H}1)$--$(\mathbf{H}3)$. Let $X^N$ solve \eqref{eq:finite-N-physical-reduced-equation} with initial value 
$X_0^N\in L^2(\Omega;H_{cu})$. Let $X$ solve \eqref{eq:physical-limit-reduced-equation}, where $X_0\in H_{cu}$ is deterministic. Then, for every $T>0$, there exists $C_T<\infty$, 
independent of $N$, such that
\begin{equation}\label{eq:reduced-mean-square-rate}
\mathbb E\sup_{0\le t\le T}|X_t^N-X_t|^2 \le C_T\left( \mathbb E|X_0^N-X_0|^2+\frac1N \right).
\end{equation}
In particular, if $\mathbb E|X_0^N-X_0|^2=O(N^{-1})$, then
\[
\mathbb E\sup_{0\le t\le T}|X_t^N-X_t|^2=O(N^{-1}).
\]
\end{proposition}

\begin{proof}
Fix $T>0$ and choose $R>1+\sup_{0\le t\le T}|X_t|$.
Adding and subtracting $u_i^{N,*}(0,x,\theta_t\omega)$ and using \eqref{eq:fixed-N-particle-center-Lipschitz}, Proposition~\ref{prop:backward-graph-N-rate},
 \eqref{eq:center-OU-N-rate}, and the invariance of $\mathbb P$ under $\theta_t$ gives
\begin{equation}\label{eq:physical-particle-state-fixed-center-N-rate}
\sup_{\substack{0\le t\le T\\ |x|\le R}} \mathbb E\left| u_i^{N,*}\bigl(0,x-z_{cu}^N(\theta_t\omega),\theta_t\omega\bigr) -u_i^{*}(0,x,\theta_t\omega) \right|^2 \le \frac{C_{T,R}}{N}.
\end{equation}
For all $x,y\in H_{cu}$, \eqref{eq:fixed-N-particle-center-Lipschitz} gives,  almost surely,
\begin{equation}\label{eq:finite-physical-particle-state-empirical-Lipschitz}
\frac1N\sum_{i=1}^N \left| u_i^{N,*}\bigl(0,x-z_{cu}^N(\theta_t\omega),\theta_t\omega\bigr) -u_i^{N,*}\bigl(0,y-z_{cu}^N(\theta_t\omega),\theta_t\omega\bigr) \right|^2
 \le C|x-y|^2.
\end{equation}
By $(\mathbf{H}2)$, Jensen inequality, and \eqref{eq:physical-particle-state-fixed-center-N-rate},
\[
\mathbb E\left| \frac1N\sum_{i=1}^N P_{cu}\left[ \nabla U\bigl(u_i^{N,*}\bigl(0,X_t-z_{cu}^N(\theta_t\omega),\theta_t\omega\bigr)\bigr) -
\nabla U\bigl(u_i^{*}(0,X_t,\theta_t\omega)\bigr) \right] \right|^2 \le \frac{C_T}{N}.
\]
Since $X$ is deterministic, the random variables $u_i^*(0,X_t,\theta_t\omega)$, $i\ge1$, are i.i.d. for each fixed $t\in[0,T]$.
Together with the invariance of $\mathbb P$ under $\theta_t$ and \eqref{eq:backward-global-weighted-moment},  this gives
\[
\mathbb E\left| \frac1N\sum_{i=1}^N P_{cu}\nabla U\bigl(u_i^{*}(0,X_t,\theta_t\omega)\bigr) - 
\mathbb E\left[ P_{cu}\nabla U\bigl(u_1^{*}(0,X_t,\omega)\bigr) \right] \right|^2 \le \frac{C_T}{N}.
\]
\eqref{eq:finite-physical-particle-state-empirical-Lipschitz} and $(\mathbf{H}2)$ imply that the finite particle drift in \eqref{eq:finite-N-physical-reduced-equation} 
is Lipschitz in $x$ with a constant independent of $N$. Subtracting \eqref{eq:physical-limit-reduced-equation} from \eqref{eq:finite-N-physical-reduced-equation} 
and applying the Cauchy--Schwarz and Burkholder--Davis--Gundy inequalities gives
\[
\begin{aligned}
\mathbb E\sup_{0\le r\le t}|X_r^N-X_r|^2\le C\mathbb E|X_0^N-X_0|^2+\frac{C_T}{N}+C_T\int_0^t\mathbb E\sup_{0\le r\le s}|X_r^N-X_r|^2\,ds.
\end{aligned}
\]
Gronwall inequality proves \eqref{eq:reduced-mean-square-rate}.
\end{proof}

Define the empirical fluctuation by
\begin{equation}\label{eq:Gamma-field-definition}
\Gamma_t^N(x) := \frac1{\sqrt N}\sum_{i=1}^N \left[ P_{cu}\nabla U\bigl(u_i^{*}(0,x,\theta_t\omega)\bigr) -
 \mathbb E\left[ P_{cu}\nabla U\bigl(u_1^{*}(0,x,\omega)\bigr) \right] \right].
\end{equation}
Since $\theta_t$ preserves $\mathbb P$, $\Gamma_t^N(x)=P_{cu}\mathbb G_{U,\theta}^N(t,0,x)$. Set
\begin{equation}\label{eq:Gamma-limit-definition}
\Gamma_t(x):=P_{cu}\mathbb G_{U,\theta}(t,0,x),
\end{equation}
where $\mathbb G_{U,\theta}$ is the limiting field in \eqref{eq:theta-driving-field-CLT}. Taking $I=[0,T_0]$ in \eqref{eq:theta-driving-field-CLT}, for every $T_0,T,R,M>0$,
\begin{equation}\label{eq:Gamma-field-limit}
\left( \mathbb G_{U,\theta}^N, \mathbb G_{K,\theta}^N, \bar W^N, \Gamma^N \right) \Longrightarrow \left( \mathbb G_{U,\theta}, \mathbb G_{K,\theta}, \bar W, \Gamma \right)
\end{equation}
in $C([0,T_0]\times[-T,0]\times B_R^{cu},H)\times
C([0,T_0]\times[-T,0]\times B_R^{cu}\times K_M,H)\times C([0,T_0],H)
\times C([0,T_0]\times B_R^{cu},H_{cu})$. The limiting fields are jointly Gaussian.
Define
\begin{equation}\label{eq:center-OU-averaged-Wiener}
\begin{aligned}
Z_{cu}^N(t)&:=\sqrt N\,z_{cu}^N(\theta_t\omega)=\sigma P_{cu}\int_{-\infty}^t e^{\Delta(t-r)}\,d\bar W^N(r),\\
Z_{cu}(t)&:=\sigma P_{cu}\int_{-\infty}^t e^{\Delta(t-r)}\,d\bar W(r).
\end{aligned}
\end{equation}
By \eqref{eq:wiener-center-at-zero} and Corollary~\ref{cor:theta-LP-Gaussian-fluctuation},
\begin{equation}\label{eq:Zcu-from-theta-Z}
\mathbf P_{cu}\mathbb Z_\theta^N(t,0) = \begin{pmatrix} Z_{cu}^N(t)\\ 0 \end{pmatrix},
\qquad
\mathbf P_{cu}\mathbb Z_\theta(t,0) = \begin{pmatrix} Z_{cu}(t)\\ 0 \end{pmatrix}.
\end{equation}
For every $T,R>0$, Lemmas~\ref{lem:backward-regularity-gaussian} and \ref{lem:backward-theta-regularity-gaussian} imply that $u_i^{*}(0,x,\theta_t\omega)$ is continuously 
Fr\'echet differentiable in $x$, and that its derivative is jointly continuous in $(t,x)\in[0,T]\times B_R^{cu}$. Estimate \eqref{eq:regularity-derivative-weighted-bound},
 together with \eqref{eq:derivative-theta-reduced-holder}, gives
\begin{equation}\label{eq:physical-particle-state-derivative-uniform-bound}
\sup_{0\le t\le T}\sup_{|x|\le R} \|D_xu_i^{*}(0,x,\theta_t\omega)\|_{\mathcal L(H_{cu},H)} \le C_{R,T} \qquad\text{a.s.},
\end{equation}
where $C_{R,T}$ is deterministic. Define the particle state fluctuation by
\begin{equation}\label{eq:physical-particle-state-fluctuation-definition}
\mathscr R_i^N(t,x)=\sqrt N\bigl[u_i^{N,*}\bigl(0,x-z_{cu}^N(\theta_t\omega),\theta_t\omega\bigr)-u_i^{*}(0,x,\theta_t\omega)\bigr].
\end{equation}
For the fields $\mathbf Y_{i,\theta}=(\bar y_{i,\theta},\hat y_{i,\theta})$ in Corollary~\ref{cor:theta-LP-Gaussian-fluctuation}, write 
$y_{i,\theta}=\bar y_{i,\theta}+\hat y_{i,\theta}$ and set
\begin{equation}\label{eq:physical-particle-state-fluctuation-limit}
\mathscr R_i(t,x) := y_{i,\theta}(t,0,x) - D_xu_i^{*}(0,x,\theta_t\omega)Z_{cu}(t).
\end{equation}
Fix $m\ge1$ and $T_0,T,R,M>0$, and set $J=[0,T_0]\times B_R^{cu}$. Corollary~\ref{cor:theta-LP-Gaussian-fluctuation} and 
Lemmas~\ref{lem:backward-regularity-gaussian} and \ref{lem:backward-theta-regularity-gaussian} give
\begin{align}
&\Bigl(\mathbb G_{U,\theta}^N,\mathbb G_{K,\theta}^N,\bar W^N,\mathbb Z_\theta^N,\Gamma^N,
\bigl(u_i^{*}(0,\cdot,\theta_\cdot\omega),\mathbf Y_i^N(0,\cdot,\theta_\cdot\omega),D_xu_i^{*}(0,\cdot,\theta_\cdot\omega)\bigr)_{1\le i\le m}\Bigr)\notag\\
&\qquad\Longrightarrow
\Bigl(\mathbb G_{U,\theta},\mathbb G_{K,\theta},\bar W,\mathbb Z_\theta,\Gamma,
\bigl(u_i^{*}(0,\cdot,\theta_\cdot\omega),\mathbf Y_{i,\theta}(\cdot,0,\cdot),D_xu_i^{*}(0,\cdot,\theta_\cdot\omega)\bigr)_{1\le i\le m}\Bigr)
\label{eq:physical-joint-convergence}
\end{align}
in
\[
\begin{aligned}
&C([0,T_0]\times[-T,0]\times B_R^{cu},H)\times C([0,T_0]\times[-T,0]\times B_R^{cu}\times K_M,H)\times C([0,T_0],H)\\
&\quad \times C([0,T_0],C([-T,0],\mathbf H))
\times C(J,H_{cu}) \times
\bigl[
C(J,H)\times C(J,\mathbf H)
\times C(J,\mathcal L(H_{cu},H))
\bigr]^m.
\end{aligned}
\]
The derivative fields are measurable with respect to the corresponding particle noises and are locally continuous by the two lemmas. The proof of 
Corollary~\ref{cor:theta-LP-Gaussian-fluctuation} applies with these derivative fields included among the fixed particle variables.  Finite summand
 removal gives independence of the limiting particle fields and their derivatives from 
 $(\mathbb G_{U,\theta},\allowbreak \mathbb G_{K,\theta},\allowbreak \bar W,\allowbreak \mathbb Z_\theta)$, and the same argument identifies the fluctuation limit.

The next proposition establishes the convergence of the particle fluctuations defined in \eqref{eq:physical-particle-state-fluctuation-definition}.

\begin{proposition}
\label{prop:physical-particle-state-fluctuation-uniform}
Assume $(\mathbf{H}1)$--$(\mathbf{H}3)$. Let $\mathscr R_i^N$ and $\mathscr R_i$ be defined by \eqref{eq:physical-particle-state-fluctuation-definition} and 
\eqref{eq:physical-particle-state-fluctuation-limit}. For every fixed particle index $i$ and every $T,R>0$,
\begin{equation}\label{eq:physical-particle-state-fluctuation-uniform-convergence}
\mathscr R_i^N \Longrightarrow \mathscr R_i \qquad \text{in } C([0,T]\times B_R^{cu},H).
\end{equation}
More precisely, for every fixed $m\ge1$,
\[
\bigl(u_i^{*}(0,\cdot,\theta_\cdot\omega),\mathscr R_i^N\bigr)_{1\le i\le m}
\Longrightarrow
\bigl(u_i^{*}(0,\cdot,\theta_\cdot\omega),\mathscr R_i\bigr)_{1\le i\le m}
\]
in $C([0,T]\times B_R^{cu},H)^{2m}$, jointly with the first five fields in \eqref{eq:physical-joint-convergence} in the spaces stated there. Conditionally on
 $(\mathbb G_{U,\theta},\allowbreak \mathbb G_{K,\theta},\allowbreak \bar W,\allowbreak \mathbb Z_\theta)$, the pairs $(u_i^{*}(0,\cdot,\theta_\cdot\omega),\mathscr R_i)$ are i.i.d.
\par
\end{proposition}

\begin{proof}
 \eqref{eq:center-OU-N-rate} gives
\[
\sup_{0\le t\le T}|z_{cu}^N(\theta_t\omega)|\longrightarrow0 \quad\text{in probability},\qquad \|Z_{cu}^N\|_{C([0,T],H_{cu})}=O_{\mathbb P}(1).
\]
 On the event $\{\sup_{0\le t\le T}|z_{cu}^N(\theta_t\omega)|\le1\}$, both $x$ and $x-z_{cu}^N(\theta_t\omega)$ belong to $B_{R+1}^{cu}$ 
whenever $|x|\le R$. Therefore, \eqref{eq:theta-fluctuation-equicontinuity}, applied on $B_{R+1}^{cu}$, gives
\begin{equation}\label{eq:random-center-fluctuation-replacement}
\sup_{\substack{0\le t\le T\\ |x|\le R}} \left\| \mathbf Y_i^N \bigl(\cdot,x-z_{cu}^N(\theta_t\omega),\theta_t\omega\bigr)
 - \mathbf Y_i^N(\cdot,x,\theta_t\omega) \right\|_{\eta,T} \longrightarrow0
\end{equation}
in probability. By \eqref{eq:physical-particle-state-fluctuation-definition} and \eqref{eq:particle-state-fluctuation},
\begin{equation}\label{eq:physical-fluctuation-two-terms}
\begin{aligned} \mathscr R_i^N(t,x) =  y_i^N \bigl(0,x-z_{cu}^N(\theta_t\omega),\theta_t\omega\bigr) 
+ \sqrt N\left[ u_i^{*}\bigl(0,x-z_{cu}^N(\theta_t\omega),\theta_t\omega\bigr) - u_i^{*}(0,x,\theta_t\omega) \right]. \end{aligned}
\end{equation}
By \eqref{eq:random-center-fluctuation-replacement}, the first term in \eqref{eq:physical-fluctuation-two-terms} differs from
 $y_i^N(0,x,\theta_t\omega)$ by $o_{\mathbb P}(1)$ in $C([0,T]\times B_R^{cu},H)$. For the second term, the fundamental theorem of calculus gives
\begin{equation}\label{eq:physical-fluctuation-center-correction-expansion}
\begin{aligned}
 &\sqrt N\left[ u_i^{*}\bigl(0,x-z_{cu}^N(\theta_t\omega),\theta_t\omega\bigr) - u_i^{*}(0,x,\theta_t\omega) \right] \\
 &\qquad= -\int_0^1 D_xu_i^{*}\bigl(0,x-a z_{cu}^N(\theta_t\omega),\theta_t\omega\bigr) Z_{cu}^N(t)\,da.
 \end{aligned}
\end{equation}
By Lemmas~\ref{lem:backward-regularity-gaussian} and \ref{lem:backward-theta-regularity-gaussian}, $D_xu_i^{*}(0,x,\theta_t\omega)$ is
 jointly continuous. Since $z_{cu}^N(\theta_\cdot\omega)\to0$ in probability in $C([0,T],H_{cu})$ and
  $\|Z_{cu}^N\|_{C([0,T],H_{cu})}=O_{\mathbb P}(1)$, \eqref{eq:physical-fluctuation-center-correction-expansion} gives
\[
\mathscr R_i^N(t,x)=y_i^N(0,x,\theta_t\omega)-D_xu_i^{*}(0,x,\theta_t\omega)Z_{cu}^N(t)+o_{\mathbb P}(1)
\]
in $C([0,T]\times B_R^{cu},H)$. Thus \eqref{eq:physical-joint-convergence} and \eqref{eq:Zcu-from-theta-Z} give
 \eqref{eq:physical-particle-state-fluctuation-uniform-convergence}, jointly with the stated fields. The argument 
 holds simultaneously for every fixed $i=1,\ldots,m$ and retains the limiting particle states, which proves the joint convergence of the pairs.
By Corollary~\ref{cor:theta-LP-Gaussian-fluctuation}, the pairs $(\mathbf U_i^{*}(\cdot,\cdot,\theta_\cdot\omega),\mathbf Y_{i,\theta})$ are conditionally i.i.d. given 
$(\mathbb G_{U,\theta},\allowbreak \mathbb G_{K,\theta},\allowbreak \bar W,\allowbreak \mathbb Z_\theta)$. 
The derivative $D_xu_i^{*}(0,\cdot,\theta_\cdot\omega)$ is determined by the same particle parameter field, while $Z_{cu}$ is 
determined by $\mathbb Z_\theta$ through \eqref{eq:Zcu-from-theta-Z}. Formula \eqref{eq:physical-particle-state-fluctuation-limit} therefore gives the stated
 conditional i.i.d. property.
 
\end{proof}

Define the corresponding corrections to the reduced drift by
\begin{equation}\label{eq:averaged-particle-state-correction-N}
\mathcal R_t^N(x) := \frac1{\sqrt N}\sum_{i=1}^N P_{cu}\left[ \nabla U\bigl(u_i^{N,*}\bigl(0,x-z_{cu}^N(\theta_t\omega),\theta_t\omega\bigr)\bigr)
- \nabla U\bigl(u_i^{*}(0,x,\theta_t\omega)\bigr) \right],
\end{equation}
\begin{equation}\label{eq:graph-average-fluctuation-formula}
\mathcal R_t(x) := \widetilde{\mathbb E}\left[ P_{cu}D\nabla U\bigl(u_1^{*}(0,x,\theta_t\omega)\bigr) \mathscr R_1(t,x) \right].
\end{equation}
Here $\widetilde{\mathbb E}$ denotes expectation over the limiting particle variable, equivalently conditional expectation given 
$\sigma(\mathbb G_{U,\theta},\allowbreak \mathbb G_{K,\theta},\allowbreak \bar W,\allowbreak \mathbb Z_\theta)$. The term $\mathcal R$ is the averaged 
contribution of the invariant graph fluctuation to the reduced drift.

The next proposition gives the first order expansion of $\mathcal R^N$ and the joint Gaussian limit of $(\Gamma^N,\mathcal R^N,\bar W^N)$.

\begin{proposition}
\label{prop:graph-average-fluctuation}
Assume $(\mathbf{H}1)$--$(\mathbf{H}3)$, and let $\mathcal R^N$ and $\mathcal R$ be defined by \eqref{eq:averaged-particle-state-correction-N} 
and \eqref{eq:graph-average-fluctuation-formula}. For every $T,R>0$,
\begin{equation}\label{eq:averaged-particle-state-linearization}
\mathcal R_t^N(x) = \frac1N\sum_{i=1}^N P_{cu}D\nabla U\bigl(u_i^{*}(0,x,\theta_t\omega)\bigr) \mathscr R_i^N(t,x) +o_{\mathbb P}(1)
\end{equation}
in $C([0,T]\times B_R^{cu},H_{cu})$. Moreover,
\begin{equation}\label{eq:graph-average-fluctuation-limit}
\mathcal R^N \Longrightarrow \mathcal R \qquad \text{in } C([0,T]\times B_R^{cu},H_{cu}).
\end{equation}
The convergence also holds jointly with $\Gamma^N$ and $\bar W^N$. Thus
\begin{equation}\label{eq:joint-direct-particle-state-fluctuation-limit}
(\Gamma^N,\mathcal R^N,\bar W^N) \Longrightarrow (\Gamma,\mathcal R,\bar W)
\end{equation}
in $C([0,T]\times B_R^{cu},H_{cu})^2 \times C([0,T],H)$, and the limiting fields are jointly Gaussian.
\end{proposition}

\begin{proof}
For $T,R>0$, set $\|f\|_{T,R}:=\sup_{0\le t\le T,\,|x|\le R}|f(t,x)|$. By \eqref{eq:physical-particle-state-fluctuation-definition} and \eqref{eq:particle-state-fluctuation},
\begin{equation}\label{eq:physical-correction-fixed-center-decomposition}
\begin{aligned} \mathscr R_i^N(t,x) =  y_i^N(0,x,\theta_t\omega)  + \sqrt N\left[ u_i^{N,*} \bigl(0,x-z_{cu}^N(\theta_t\omega),\theta_t\omega\bigr)
 - u_i^{N,*}(0,x,\theta_t\omega) \right]. \end{aligned}
\end{equation}
Hence \eqref{eq:fixed-N-particle-center-Lipschitz}, \eqref{eq:particle-fluctuation-local-higher-moment}, \eqref{eq:Zcu-from-theta-Z}, 
and \eqref{eq:averaged-noise-local-theta-pmoment} imply, for the fixed $p>4$,
\begin{align}
\sup_{N\ge2}\sup_{1\le i\le N} \mathbb E\|\mathscr R_i^N\|_{T,R}^p &\le C_{T,R,p}, \label{eq:physical-correction-uniform-pmoment}\\
 \sup_{N\ge2} \mathbb E\left[ \frac1N\sum_{i=1}^N\|\mathscr R_i^N\|_{T,R}^2 \right] &<\infty. \label{eq:physical-empirical-L2-field-tightness-result}
\end{align}
Taylor's formula and $(\mathbf{H}2)$ give
\[\left\| \mathcal R^N - \frac1N\sum_{i=1}^N P_{cu}D\nabla U(u_i^{*}(0,\cdot,\theta_\cdot\omega))\mathscr R_i^N \right\|_{T,R}
 \le \frac{C}{\sqrt N}\, \frac1N\sum_{i=1}^N\|\mathscr R_i^N\|_{T,R}^2.\]
By \eqref{eq:physical-empirical-L2-field-tightness-result}, the right hand side converges to zero in $L^1$, which proves \eqref{eq:averaged-particle-state-linearization}.
Proposition~\ref{prop:physical-particle-state-fluctuation-uniform} gives tightness of the laws of 
$(u_1^{*}(0,\cdot,\theta_\cdot\omega),\mathscr R_1^N)$ in $C([0,T]\times B_R^{cu},H)^2$. For every $\varepsilon,\delta>0$, choose a compact set $K$ such that
\[
\sup_{N\ge2}\mathbb P\bigl((u_1^{*}(0,\cdot,\theta_\cdot\omega),\mathscr R_1^N)\notin K\bigr)\le \varepsilon\delta.
\]
By exchangeability and Markov inequality,
\[
\mathbb P\left(\frac1N\sum_{i=1}^N
\mathbf 1_{\{(u_i^{*}(0,\cdot,\theta_\cdot\omega),\mathscr R_i^N)\notin K\}}>\delta\right)
\le \varepsilon.
\]
Hence the empirical measures $\frac1N\sum_{i=1}^N \delta_{(u_i^{*}(0,\cdot,\theta_\cdot\omega),\mathscr R_i^N)}$
are tight in $\mathcal P\bigl(C([0,T]\times B_R^{cu},H)^2\bigr)$ endowed with the weak topology, by \cite[Proposition~2.2(ii)]{Sznitman1991}. 
By the joint convergence and conditional i.i.d. property in Proposition~\ref{prop:physical-particle-state-fluctuation-uniform}, the exchangeability argument from
 Step~3 in the proof of Proposition~\ref{prop:joint-backward-fluctuation} applies.
Thus, in every joint subsequential limit, the empirical measure equals, almost surely, the conditional law of $(u_1^{*}(0,\cdot,\theta_\cdot\omega),\mathscr R_1)$ 
given $(\mathbb G_{U,\theta},\allowbreak \mathbb G_{K,\theta},\allowbreak \bar W,\allowbreak \mathbb Z_\theta)$.
  See \cite[Remark~2.3(1)]{Sznitman1991} and \cite{DiaconisFreedman1980}.

Let $\chi_L:[0,\infty)\to[0,1]$ be a continuous function equal to $1$ on $[0,L]$ and to $0$ on $[2L,\infty)$. The map
$(\widehat u,r) \longmapsto \chi_L(\|r\|_{T,R}) P_{cu}D\nabla U(\widehat u)r$
is bounded and continuous. The image of every compact set under this map is compact in $C([0,T]\times B_R^{cu},H_{cu})$. Approximate this compact image uniformly by a finite dimensional subspace and apply scalar weak convergence to its finite dimensional coordinates. Tightness of the measures and boundedness of the map control the contribution outside the 
compact set. Hence
\begin{align*}
 &\frac1N\sum_{i=1}^N \chi_L(\|\mathscr R_i^N\|_{T,R}) P_{cu}D\nabla U(u_i^{*}(0,\cdot,\theta_\cdot\omega))\mathscr R_i^N\\
 &\qquad\Longrightarrow \widetilde{\mathbb E}\left[ \chi_L(\|\mathscr R_1\|_{T,R}) P_{cu}D\nabla U(u_1^{*}(0,\cdot,\theta_\cdot\omega))\mathscr R_1 \right].
\end{align*}
By \eqref{eq:physical-empirical-L2-field-tightness-result}, the expected contribution from $\{\|\mathscr R_i^N\|_{T,R}>L\}$ is bounded 
by $C_{T,R}/L$, uniformly in $N$. 
Proposition~\ref{prop:physical-particle-state-fluctuation-uniform}, the Portmanteau theorem, and \eqref{eq:physical-correction-uniform-pmoment} give the same bound for the limit. 
Letting $L\to\infty$ proves \eqref{eq:graph-average-fluctuation-limit}. Since $\Gamma_t^N(x)=P_{cu}\mathbb G_{U,\theta}^N(t,0,x)$, the same joint convergence
 argument also gives \eqref{eq:joint-direct-particle-state-fluctuation-limit}.

To prove joint Gaussianity, let $\mathcal H_G$ be the closed linear span of all scalar coordinates of
 $(\mathbb G_{U,\theta},\allowbreak \mathbb G_{K,\theta},\allowbreak \bar W,\allowbreak \mathbb Z_\theta)$ in the $L^2$ space of its joint law. Let $P$ be the 
 law of the particle field $(t,r,x)\mapsto\mathbf U_1^{*}(r,x,\theta_t\omega)$, which is independent of these Gaussian fields by 
 Corollary~\ref{cor:theta-LP-Gaussian-fluctuation}. Approximation of the particle variable by simple random variables and $L^2$ continuity of the
Gaussian fields show that the forcing in \eqref{eq:gaussian-linear-fluctuation-equation}, with the replacements specified in that corollary, belongs at 
each time to $\mathcal H_G\widehat\otimes L^2(P,\mathbf H)$. The operator $\mathscr L_\xi$ acts only on the particle variable and time, so each
term of its Neumann series preserves this subspace. The contraction estimate gives convergence in weighted $L^2$. Closedness, 
\eqref{eq:physical-particle-state-derivative-uniform-bound}, \eqref{eq:Zcu-from-theta-Z}, and \eqref{eq:physical-particle-state-fluctuation-limit} therefore yield
$\mathscr R_1(t,x)\in\mathcal H_G\widehat\otimes L^2(P,H).$ For each fixed $(t,x)$, particle averaging
\[
F\longmapsto\widetilde{\mathbb E}\bigl[P_{cu}D\nabla U(u_1^{*}(0,x,\theta_t\omega))F\bigr]
\]
is a deterministic bounded linear map into $\mathcal H_G\widehat\otimes H_{cu}$, with norm at most $L_U$ by $(\mathbf{H}2)$. Thus every scalar coordinate of $\mathcal R_t(x)$ 
belongs to $\mathcal H_G$. Since this is a closed Gaussian linear space, every finite collection of coordinates of $(\Gamma,\mathcal R,\bar W)$ is jointly Gaussian.

\end{proof}

\subsection{Gaussian fluctuation of the center unstable dynamics}

We now use Proposition~\ref{prop:graph-average-fluctuation} to prove the Gaussian limit for the reduced trajectories.
 Define the deterministic drift in \eqref{eq:physical-limit-reduced-equation} by
\[
b(x):=L_{cu}x+\mathbb E\left[P_{cu}\nabla U\bigl(u_1^{*}(0,x,\omega)\bigr)\right],\qquad x\in H_{cu}.
\]
By $(\mathbf{H}2)$, Lemma~\ref{lem:backward-regularity-gaussian}, \eqref{eq:physical-particle-state-derivative-uniform-bound},
 and dominated convergence, $b$ is continuously Fr\'echet differentiable, with
\begin{equation}\label{eq:limiting-reduced-drift-derivative-formula}
\begin{aligned}
Db(x)=L_{cu}+\mathbb E\left[P_{cu}D\nabla U\bigl(u_1^{*}(0,x,\omega)\bigr)D_xu_1^{*}(0,x,\omega)\right].
\end{aligned}
\end{equation}
Since $H_{cu}$ is finite dimensional, the derivative in \eqref{eq:limiting-reduced-drift-derivative-formula} is uniformly continuous on every closed ball in $H_{cu}$.

\begin{proof}[Proof of Theorem~\ref{thm:physical-cu-gaussian}]

Proposition~\ref{prop:graph-average-fluctuation} gives \eqref{eq:intro-joint-direct-particle-state-fluctuation-limit} and the joint Gaussianity of its limit. 
Proposition~\ref{prop:reduced-mean-square-rate}, with $X_0^N=X_0$, gives \eqref{eq:intro-reduced-mean-square-rate}, and hence $X^N\to X$ uniformly in probability on 
$[0,T]$. It remains to prove \eqref{eq:physical-cu-Gaussian-convergence}. Set $Z^N=\sqrt N(X^N-X)$, so that $Z_0^N=0$. For $0\le s\le T$, set
$A_s^N:=\int_0^1 Db\bigl(X_s+\theta(X_s^N-X_s)\bigr)\,d\theta.$
By the mean value theorem, subtracting \eqref{eq:physical-limit-reduced-equation} from \eqref{eq:finite-N-physical-reduced-equation} and using \eqref{eq:Gamma-field-definition} 
and \eqref{eq:averaged-particle-state-correction-N} gives
\begin{equation}\label{eq:ZN-exact-decomposition}
\begin{aligned}
Z_t^N=\int_0^t A_s^N Z_s^N\,ds+\int_0^t\bigl[\Gamma_s^N(X_s^N)+\mathcal R_s^N(X_s^N)\bigr]ds+\sigma P_{cu}\bar W_t^N.
\end{aligned}
\end{equation}
Choose $R>1+\sup_{0\le t\le T}|X_t|$ and define $E_N:=\left\{\sup_{0\le t\le T}|X_t^N-X_t|\le1\right\}.$ By \eqref{eq:reduced-mean-square-rate}, $\mathbb P(E_N)\to1$. On $E_N$,
$$
X_s^N\in B_R^{cu},\qquad X_s+\theta(X_s^N-X_s)\in B_R^{cu}, \qquad 0\le s\le T,\quad 0\le\theta\le1.
$$
The uniform continuity of the derivative in \eqref{eq:limiting-reduced-drift-derivative-formula} then gives
\[
\sup_{0\le s\le T}\|A_s^N-Db(X_s)\|_{\mathcal L(H_{cu})}\longrightarrow0
\quad\text{in probability}.
\]
Since $X$ is deterministic, \eqref{eq:joint-direct-particle-state-fluctuation-limit}, \eqref{eq:reduced-mean-square-rate}, and Slutsky's theorem give
$$
\bigl(\Gamma^N,\mathcal R^N,\bar W^N,X^N\bigr)
\Longrightarrow
\bigl(\Gamma,\mathcal R,\bar W,X\bigr)
$$
in $C([0,T]\times B_R^{cu},H_{cu})^2\times C([0,T],H)\times C([0,T],H_{cu}).$
On $E_N$, $X^N([0,T])\subset B_R^{cu}$. The composition map $(f,x)\longmapsto \bigl[t\mapsto f(t,x(t))\bigr]$
is continuous on $C([0,T]\times B_R^{cu},H_{cu})\times C([0,T],B_R^{cu})$. Since $\mathbb P(E_N^c)\to0$, the continuous mapping theorem gives
$$
\bigl(\Gamma^N(\cdot,X_\cdot^N),\mathcal R^N(\cdot,X_\cdot^N),\bar W^N\bigr)
\Longrightarrow
\bigl(\Gamma(\cdot,X_\cdot),\mathcal R(\cdot,X_\cdot),\bar W\bigr).
$$
For given coefficients and driving terms, Gronwall inequality shows that the solution of \eqref{eq:ZN-exact-decomposition} 
depends continuously on them under uniform convergence. Therefore, Slutsky's theorem and the continuous mapping 
theorem \cite[Theorem~4.27]{Kallenberg2002} give \eqref{eq:physical-cu-Gaussian-convergence}, with limit equal to the unique solution
 of \eqref{eq:physical-cu-Gaussian-limit-equation}. Since $X$ and $Db(X)$ are deterministic, this solution is a linear function of the jointly
Gaussian fields $(\Gamma,\mathcal R,\bar W)$. Hence $Z$ is Gaussian.

\end{proof}

\begin{example}
\label{ex:linear-cancellation-test}

Let $\gamma=1$, $A=\Delta+I$, $H_{cu}=\operatorname{span}\{e_1\}$, $Q=e_1\otimes e_1$, and $\nabla K=0$. Let $\nabla U(u)=a\langle u,e_1\rangle e_1$ 
with $|a|<1/32$. Then $(\mathbf{H}1)$--$(\mathbf{H}3)$ hold with $\alpha=0$, $\beta=-1$, and $\eta=-1/2$.  Identifying the first mode with $\mathbb R$, the particle system is
$$
du_i^N=[a u_i^N+\gamma(\bar u^N-u_i^N)]dt+\sigma\,d\beta_i.
$$
All higher modes of the backward graph vanish. Define
$$
\zeta_i(t)=\sigma\int_{-\infty}^t e^{-(\gamma-a)(t-s)}d\beta_i(s), \qquad\bar\zeta_N(t)=\frac1N\sum_{i=1}^N\zeta_i(t).
$$
For $x\in H_{cu}$,
$$
u_i^{N,*}\bigl(0,x-z_{cu}^N(\theta_t\omega),\theta_t\omega\bigr) =x+\zeta_i(t)-\bar\zeta_N(t), \qquad u_i^{*}(0,x,\theta_t\omega)=x+\zeta_i(t).
$$
Hence
$$
\Gamma_t^N(x)=a\sqrt N\,\bar\zeta_N(t), \qquad\mathcal R_t^N(x)=-a\sqrt N\,\bar\zeta_N(t).
$$
Thus the two drift fluctuations cancel. If $\bar\beta^N=N^{-1/2}\sum_{i=1}^N\beta_i,$
then
$$
dX_t^N=aX_t^Ndt+\sigma N^{-1/2}d\bar\beta_t^N.
$$
If $X_0^N=X_0$ is deterministic, the limiting fluctuation satisfies
$$
dZ_t=aZ_tdt+\sigma d\bar\beta_t, \qquad Z_0=0.
$$
Hence the particle state correction cancels the direct empirical drift fluctuation in this linear case.
\end{example}

\begin{example}
\label{ex:nonlinear-covariance-example}
Let $e_1$ be the normalized first Dirichlet eigenfunction and consider

\begin{equation}\label{eq:sine-example-model}
\begin{aligned}
&A=\Delta+I,\qquad \gamma=1,\qquad Q=e_1\otimes e_1,\qquad \nabla K=0,\\
&\nabla U_\varepsilon(u)=\varepsilon\sin\langle u,e_1\rangle e_1,
\qquad \sigma>0.
\end{aligned}
\end{equation}
For $|\varepsilon|<1/32$, $(\mathbf{H}1)$--$(\mathbf{H}3)$ hold with
\[
\alpha=0,\qquad \beta=-1,\qquad \eta=-\frac12,
\qquad q=32|\varepsilon|.
\]
Write $u_i^{\varepsilon,*}$ for the limiting backward solution corresponding to $\nabla U_\varepsilon$. For each fixed $\varepsilon$, first let $N\to\infty$ and
 then consider the expansion as $\varepsilon\to0$. 
In this example, a remainder $R_\varepsilon(t)=O_{L^2}(\varepsilon^k)$ uniformly on bounded time intervals satisfies 
$\sup_{t\in I}\|R_\varepsilon(t)\|_{L^2}=O(|\varepsilon|^k)$ for every compact interval $I$.

All higher modes of the backward graph vanish. Identifying the first mode with $\mathbb R$ and taking $X_0=0$, symmetry and uniqueness give $X_t=0$ and
\begin{equation}\label{eq:sine-example-stationary}
u_1^{\varepsilon,*}(0,0,\theta_t\omega)=V_\varepsilon(t),\qquad dV_\varepsilon=(-\gamma V_\varepsilon+\varepsilon\sin V_\varepsilon)dt+\sigma dB_t,
\end{equation}
where $V_\varepsilon$ is stationary and symmetric. Let
$$
\zeta(t)=\sigma\int_{-\infty}^t e^{-\gamma(t-r)}dB_r, \qquad v=\frac{\sigma^2}{2\gamma}, \qquad m=e^{-v/2},
\qquad c(h)=v e^{-\gamma|h|}.
$$
Then
\begin{equation}\label{eq:sine-V-zeta-bound}
\sup_{t\in\mathbb R}|V_\varepsilon(t)-\zeta(t)| \le \frac{|\varepsilon|}{\gamma}.
\end{equation}
Let $P_0$ be the law of $B$ and let $\mathbb I$ be an isonormal Gaussian process over $L_0^2(P_0)$. Define
$$
G_{\sin}(t)=\mathbb I(\sin\zeta(t)), \qquad \mathcal Z_\gamma(t)=\mathbb I(\zeta(t)), \qquad B_t=\mathbb I(B_t).
$$
Then $\bar B$ is a Brownian motion and
$$
\mathcal Z_\gamma(t)=\sigma\int_{-\infty}^t e^{-\gamma(t-r)}d\bar B_r.
$$
For $F_\varepsilon(t)=\Gamma_t^\varepsilon(0)+\mathcal R_t^\varepsilon(0),$
we have, uniformly for $t$ in bounded intervals,
\begin{equation}\label{eq:sine-fluctuation-expansion}
\Gamma_t^\varepsilon(0) =\varepsilon G_{\sin}(t)+O_{L^2}(\varepsilon^2),\qquad
\mathcal R_t^\varepsilon(0) =-\varepsilon m\mathcal Z_\gamma(t)+O_{L^2}(\varepsilon^2).
\end{equation}
The first estimate follows from $\Gamma_t^\varepsilon(0)=\varepsilon\mathbb I(\sin V_\varepsilon(t))$
and \eqref{eq:sine-V-zeta-bound}. At $\varepsilon=0$, the operator $\mathscr L_\xi$ vanishes, 
and \eqref{eq:wiener-backward-contribution} and \eqref{eq:Zcu-from-theta-Z} give $\mathscr R_1^0(t,0)=-\mathcal Z_\gamma(t).$
Moreover, \eqref{eq:gaussian-linearized-LP-operator-bound},
\eqref{eq:gaussian-linear-fluctuation-resolvent}, and the variational equation \eqref{eq:regularity-variational-resolvent} give, uniformly for $t$ in bounded intervals,
$$
\mathbf Y_{1,\theta}^\varepsilon(t,\cdot,0)-\mathbb Z_\theta(t,\cdot) =O_{L^2}(|\varepsilon|), \qquad D_xu_1^{\varepsilon,*}(0,0,\theta_t\omega)=1+O_{L^2}(|\varepsilon|).
$$
Hence
$$
\mathscr R_1^\varepsilon(t,0)=-\mathcal Z_\gamma(t)+O_{L^2}(|\varepsilon|),
$$
and the second estimate in \eqref{eq:sine-fluctuation-expansion} follows from
\eqref{eq:graph-average-fluctuation-formula} and $|\mathbb E\cos V_\varepsilon-m|\le \frac{|\varepsilon|}{\gamma}.$
Since $(\zeta(t),\zeta(s))$ is Gaussian with variance $v$ and covariance $c(t-s)$,
$$
\mathbb E[\sin\zeta(t)\sin\zeta(s)]=e^{-v}\sinh c(t-s), \qquad \mathbb E[\sin\zeta(t)\zeta(s)] =m c(t-s).
$$
Therefore, uniformly for $t,s$ in bounded intervals,
\begin{align}
\operatorname{Cov}(F_\varepsilon(t),F_\varepsilon(s))&=\varepsilon^2e^{-v}\{\sinh c(t-s)-c(t-s)\}+O(|\varepsilon|^3),\label{eq:sine-corrected-covariance}\\
\operatorname{Var}(\mathcal R_t^\varepsilon(0))&=\varepsilon^2ve^{-v}+O(|\varepsilon|^3),\label{eq:sine-original-variable-variance}\\
\operatorname{Var}(F_\varepsilon(t))&=\varepsilon^2e^{-v}(\sinh v-v)+O(|\varepsilon|^3).
\label{eq:sine-positive-forcing-variance}
\end{align}
In particular, $\mathcal R^\varepsilon$ is nonzero for all sufficiently small nonzero $\varepsilon$.
Let $Z_\varepsilon^\Gamma$ solve \eqref{eq:physical-cu-Gaussian-limit-equation} at $X=0$ with the forcing term $\mathcal R^\varepsilon$ removed. 
Since the deterministic coefficient in \eqref{eq:physical-cu-Gaussian-limit-equation} at $X=0$ is $O(|\varepsilon|)$, variation of constants gives
$$
Z_\varepsilon^\Gamma(T)=\sigma\bar B_T+O_{L^2}(|\varepsilon|), \qquad
Z_\varepsilon(T)-Z_\varepsilon^\Gamma(T) =-\varepsilon m\int_0^T\mathcal Z_\gamma(s)ds+O_{L^2}(\varepsilon^2).
$$
Moreover,
$$
\mathbb E[\bar B_T\mathcal Z_\gamma(s)] =\frac{\sigma}{\gamma}(1-e^{-\gamma s}), \qquad 0\le s\le T.
$$
Thus, for every $T>0$,
\begin{equation}\label{eq:sine-trajectory-variance-error}
\operatorname{Var}Z_\varepsilon(T)-\operatorname{Var}Z_\varepsilon^\Gamma(T)
=-\frac{2\varepsilon m\sigma^2}{\gamma} \left[T-\frac{1-e^{-\gamma T}}{\gamma}\right] +O_T(\varepsilon^2).
\end{equation}
Hence the particle state correction changes the limiting trajectory variance.
\end{example}

\section*{Declaration of generative AI and AI-assisted technologies\\ in the manuscript preparation process}

The human authors initiated and led the research program, formulated the mathematical
problem, developed the main results, and independently verified all arguments. Through discussions with ChatGPT, 
the human authors obtained numerous references and explanations of their content, some of these references were relevant to the present work and proved useful
in developing the proofs. ChatGPT also assisted in improving notation, grammar, and presentation. The final manuscript was carefully revised and approved by the
human authors, who take full responsibility for its content.

\appendix

\section{Backward regularity estimates}
\label{app:backward-regularity}
\par
\subsection*{Proof of Lemma~\ref{lem:backward-regularity-gaussian}}
\begin{proof}
Define
\begin{equation}\label{eq:regularity-OU-path}
Z_{\mathrm{OU}}(t):=\sigma\int_{-\infty}^t e^{(A-\gamma I)(t-r)}\,dW_1(r).
\end{equation}
This is the stationary stochastic convolution in the centered equation of \eqref{eq:limiting-backward-LP-system}. By stationarity and the maximal inequality
 for stochastic convolutions \cite[Chapter~5]{DaPratoZabczyk2014},
\begin{align}\label{eq:regularity-OU-weighted-moment}
\left(\mathbb E\|Z_{\mathrm{OU}}\|_{C_\eta^-}^p\right)^{1/p} &\le \sum_{k=0}^\infty e^{\eta k} \left(\mathbb E\sup_{-k-1\le r\le-k}|Z_{\mathrm{OU}}(r)|^p\right)^{1/p} \\ 
&\le C_p\sum_{k=0}^\infty e^{\eta k}<\infty.
\end{align}
Since $(1+a)\nu>\eta$, the semigroup bounds
\eqref{eq:fixed-N-uniform-dichotomy} and
\eqref{eq:standing-LP-gap} yield

The Lyapunov--Perron contraction estimate and \eqref{eq:regularity-OU-weighted-moment} give
\begin{equation}\label{eq:regularity-LP-weighted-mixed-bound}
\sup_{|\xi|\le R} \left(\mathbb E\|\mathbf U_1^{*}(\cdot,\xi)\|_{C_\eta^-}^p\right)^{1/p} \le C_R\left(1+ \left(\mathbb E\|Z_{\mathrm{OU}}\|_{C_\eta^-}^p\right)^{1/p}\right)<\infty.
\end{equation}
Applying the Lyapunov--Perron contraction estimate to $\mathbf U_1^{*}(\cdot,\xi)-\mathbf U_1^{*}(\cdot,\zeta)$ gives
\begin{equation}\label{eq:regularity-center-mixed-Lipschitz}
\left(\mathbb E\|\mathbf U_1^{*}(\cdot,\xi)-\mathbf U_1^{*}(\cdot,\zeta)\|_{C_\eta^-}^p\right)^{1/p} \le C_R|\xi-\zeta|.
\end{equation}
The  equation for $\bar u^{*}$ in \eqref{eq:limiting-backward-LP-system} gives 
$\|\bar u^{*}(\cdot,\xi)-\bar u^{*}(\cdot,\zeta)\|_{C_\eta^-}\le C_R|\xi-\zeta|$.
 Hence, with $q_s=\ell_{\mathrm{LP}}/(\eta-\beta)<1$,
\[
\|\hat u_1^{*}(\cdot,\xi)-\hat u_1^{*}(\cdot,\zeta)\|_{C_\eta^-}
\le q_s\|\hat u_1^{*}(\cdot,\xi)-\hat u_1^{*}(\cdot,\zeta)\|_{C_\eta^-}
+C_R|\xi-\zeta|.
\]
Absorbing the first term gives \eqref{eq:backward-xi-weighted-Lipschitz}. Since $\sup_{-T\le r\le0}|f(r)|\le e^{-\eta T}\|f\|_{C_\eta^-},$
\eqref{eq:backward-xi-Lipschitz} follows. Together with \eqref{eq:regularity-LP-weighted-mixed-bound} at $\xi=0$,
 \eqref{eq:backward-xi-weighted-Lipschitz} gives \eqref{eq:backward-global-weighted-moment}.
For differentiability with respect to the center parameter, choose $\eta<\nu<\min(\alpha,0)$ sufficiently close to $\eta$ such that
\begin{equation}\label{eq:neighboring-weight-choice}
\ell_{\mathrm{LP}}\left(\frac1{\alpha-\nu}+\frac1{\nu-\beta}\right)<1, \qquad 0<a<\min\{1,(\nu-\eta)/(-\nu)\}.
\end{equation}
Then $(1+a)\nu>\eta$. Repeating the estimates leading to \eqref{eq:regularity-LP-weighted-mixed-bound} and \eqref{eq:regularity-center-mixed-Lipschitz}
with $\nu$ in place of $\eta$ gives solutions in $L^p(\Omega;C_\nu^-(\mathbf H))$. Since $C_\nu^-\subset C_\eta^-$, uniqueness in $L^p(\Omega;C_\eta^-(\mathbf H))$ 
identifies these solutions with $\mathbf U_1^{*}(\cdot,\xi)$.
Using \eqref{eq:limiting-backward-LP-system} at $\xi+h$ and $\xi$ in the $C_\nu^-$ norm, and using
\eqref{eq:neighboring-weight-choice}, we obtain
\[
\|\mathbf U_1^{*}(\cdot,\xi+h)-\mathbf U_1^{*}(\cdot,\xi)\|_{C_\nu^-}\le C_R|h| \qquad\text{a.s.}
\]
Let $V^\xi h$ solve the linear variational equation
\begin{equation}\label{eq:regularity-variational-resolvent}
V^\xi h=(r\mapsto e^{L_{cu}r}h\oplus0)+\mathscr L_\xi V^\xi h,
\end{equation}
where $\mathscr L_\xi$ is the linear operator obtained by differentiating the drift at $\mathbf U_1^{*}(\cdot,\xi)$. 
By $(\mathbf{H}2)$ and \eqref{eq:neighboring-weight-choice}, 
$\mathscr L_\xi$ is a contraction on $L^\infty(\Omega;C_\nu^-(\mathbf H))$, which gives a unique solution and
\begin{equation}\label{eq:regularity-derivative-weighted-bound}
\sup_{|\xi|\le R} \|V^\xi\|_{\mathcal L(H_{cu},C_\nu^-)}\le C_R \quad\text{a.s.}, \qquad C_\nu^-\subset C_\eta^-.
\end{equation}
To identify the derivative, subtract \eqref{eq:regularity-variational-resolvent} 
from the exact difference equation. By $(\mathbf{H}2)$, the remainder $R_h(r)$ satisfies
\[|R_h(r)| \le C_R\min\{|h|^2e^{2\nu r},|h|e^{\nu r}\} \le C_R|h|^{1+a}e^{(1+a)\nu r}, \qquad r\le0.\]
Since $(1+a)\nu>\eta$, the semigroup bounds \eqref{eq:fixed-N-uniform-dichotomy} and \eqref{eq:standing-LP-gap} yield
\begin{equation}\label{eq:parameter-differentiability-remainder}
\sup_{|\xi|\le R} \|\mathbf U_1^{*}(\cdot,\xi+h)-\mathbf U_1^{*}(\cdot,\xi)-V^\xi h\|_{C_\eta^-} \le C_R|h|^{1+a}, \qquad |h|\le1,
\end{equation}
almost surely. Thus $V^\xi=D_\xi\mathbf U_1^{*}(\cdot,\xi)$. Subtracting the variational equations at $\xi$ and $\zeta$ and using $(\mathbf{H}2)$,
the preceding pathwise Lipschitz estimate in $C_\nu^-$, \eqref{eq:regularity-derivative-weighted-bound}, and the Lyapunov--Perron estimate in 
$L^\infty(\Omega;C_\eta^-(\mathbf H))$ gives
\[
\|D_\xi\mathbf U_1^{*}(\cdot,\xi)-D_\xi\mathbf U_1^{*}(\cdot,\zeta)\|_{\mathcal L(H_{cu},C_\eta^-)} \le C_R|\xi-\zeta|^a.
\]
Since $H_{cu}$ is finite dimensional, a countable dense set argument gives a common full probability set on which the asserted $C^1$ regularity holds.
For the spatial and temporal regularity, \eqref{eq:positive-fractional-semigroup-gaussian}, \eqref{eq:fixed-N-uniform-dichotomy}, and the It\^o isometry give
\begin{align}
\mathbb E\|Z_{\mathrm{OU}}(t)\|_{V_\vartheta}^2 &=\sigma^2\int_0^\infty \|(I-\Delta)^\vartheta S_s(u)Q^{1/2}\|_{\mathcal L_2(H)}^2\,du \notag\\ 
&\le C\|Q^{1/2}\|_{\mathcal L_2(H)}^2 \int_0^\infty(1+u^{-2\vartheta})e^{-2cu}\,du<\infty, \label{eq:regularity-stochastic-convolution-Vtheta}
\end{align}
where $\vartheta<1/2$.  Choose $p_1>p$ so that $\vartheta+1/p_1<1/2$.
The factorization method for stochastic convolutions \cite[Section~5.3.1]{DaPratoZabczyk2014}, together with \eqref{eq:positive-fractional-semigroup-gaussian} 
and \eqref{eq:fixed-N-uniform-dichotomy}, gives
\begin{equation}\label{eq:OU-Vtheta-maximal-moment}
\mathbb E\sup_{-T\le t\le0} \|Z_{\mathrm{OU}}(t)\|_{V_\vartheta}^{p}<\infty.
\end{equation}
 For the  convolution,
\begin{equation}\label{eq:regularity-deterministic-Vtheta}
\sup_{-T\le t\le0} \left\|\int_{-\infty}^tS_s(t-r)f(r)\,dr\right\|_{V_\vartheta} \le C_{T,\vartheta}\sup_{r\le0}e^{-\eta r}|f(r)|,
\end{equation}
by \eqref{eq:positive-fractional-semigroup-gaussian}, \eqref{eq:fixed-N-uniform-dichotomy}, and $\beta<\eta$. Equations \eqref{eq:backward-global-weighted-moment},
\eqref{eq:regularity-derivative-weighted-bound}, \eqref{eq:OU-Vtheta-maximal-moment}, and \eqref{eq:regularity-deterministic-Vtheta}
imply \eqref{eq:backward-Vtheta-moment} and \eqref{eq:backward-xi-derivative-Vtheta-moment}.

Let $h=t-s\in(0,1]$. For the stochastic convolution,
\[Z_{\mathrm{OU}}(t)-Z_{\mathrm{OU}}(s) =(S_s(h)-I)Z_{\mathrm{OU}}(s) +\sigma\int_s^tS_s(t-r)\,dW_1(r).\]
By \eqref{eq:semigroup-increment-gaussian}, \eqref{eq:regularity-stochastic-convolution-Vtheta}, and the BDG inequality,
\begin{equation}\label{eq:regularity-stochastic-increment-bound}
\mathbb E|Z_{\mathrm{OU}}(t)-Z_{\mathrm{OU}}(s)|^p \le C_{T,p}\left(h^{p\vartheta}+h^{p/2}\right) \le C_{T,p}h^{p\vartheta}.
\end{equation}
Taking time increments in \eqref{eq:limiting-backward-LP-system} and \eqref{eq:regularity-variational-resolvent}, and using
\eqref{eq:semigroup-increment-gaussian}, \eqref{eq:regularity-stochastic-increment-bound} and the preceding moment bounds, we obtain
\begin{align*}
\mathbb E\sup_{|\xi|\le R} |u_1^{*}(t,\xi)-u_1^{*}(s,\xi)|^p
 &\le C_{T,R,p}|t-s|^{p\vartheta},\\ \mathbb E\sup_{|\xi|\le R} \|D_\xi u_1^{*}(t,\xi)-D_\xi u_1^{*}(s,\xi)\|^p &\le C_{T,R,p}|t-s|^{p\vartheta}.
\end{align*}
The Garsia--Rodemich--Rumsey inequality \cite{GarsiaRodemichRumsey1970}, applied to these increment bounds with values in $C(B_R^{cu},H)$ and 
$C(B_R^{cu},\mathcal L(H_{cu},H))$, respectively, gives the required $L^p$ estimates for every $0<\rho<\vartheta-1/p$.
 Hence $r\mapsto u_1^{*}(r,\xi)$ and $r\mapsto D_\xi u_1^{*}(r,\xi)$ admit $C^\rho$ versions, consistently with the Kolmogorov
 continuity criterion \cite[Theorem~3.23]{Kallenberg2002}. This proves \eqref{eq:backward-time-holder} and \eqref{eq:backward-xi-derivative-time-holder}.
\end{proof}

\par
\subsection*{Proof of Lemma~\ref{lem:backward-theta-regularity-gaussian}}
\begin{proof}
Fix $1/p<\rho_0<\rho<\vartheta-1/p$. Let $Z_{\mathrm{OU},i}$ be defined by \eqref{eq:regularity-OU-path} with $W_i$ in place of $W_1$. By $(\mathbf{H}1)$,
$Z_{\mathrm{OU},i}(r,\theta_t\omega)=Z_{\mathrm{OU},i}(r+t,\omega).$
With $Z_{\mathrm{OU}}=Z_{\mathrm{OU},1}$, set
\[K_{\mathrm{OU}} :=\sup_{r\le0}e^{-\eta r} \bigl[t\mapsto Z_{\mathrm{OU}}(r,\theta_t\omega)\bigr]_ {C^{\rho_0}([0,T_0],H)}.\]
By \eqref{eq:regularity-stochastic-increment-bound} and the Garsia--Rodemich--Rumsey inequality, the $C^{\rho_0}$ seminorm on every
 bounded time interval has a finite $p$-moment. Stationarity and a decomposition of $(-\infty,0]$ into the intervals $[-k-1,-k]$, 
 together with $\eta<0$, then give $\mathbb E K_{\mathrm{OU}}^p<\infty$.

For fixed $(r,\xi)$, $\bar u^*(r,\xi)$ and $\mathbb E\nabla U(u_1^*(r,\xi))$ are deterministic, and the map $x\mapsto\widetilde{\mathbb E} \nabla K(\widetilde u_1^*(r,\xi)-x)$
 is deterministic and $L_K$-Lipschitz. Subtracting the third equation of \eqref{eq:limiting-backward-LP-system}
at $\theta_t\omega$ and $\theta_v\omega$,
and using $q_s<1$, gives
\begin{equation}\label{eq:theta-LP-weighted-holder-bound}
\sup_{|\xi|\le R} \|\mathbf U_1^{*}(\cdot,\xi,\theta_t\omega) -\mathbf U_1^{*}(\cdot,\xi,\theta_v\omega)\|_{C_\eta^-} \le \frac{C_R}{1-q_s}K_{\mathrm{OU}}|t-v|^{\rho_0}.
\end{equation}
Taking moments proves \eqref{eq:backward-theta-weighted-holder}. Since $\sup_{-T\le r\le0}|f(r)| \le e^{-\eta T}\|f\|_{C_\eta^-},$ \eqref{eq:backward-theta-holder} follows.
For \eqref{eq:backward-joint-theta-regularity}, \eqref{eq:backward-global-weighted-moment} and \eqref{eq:theta-LP-weighted-holder-bound} give
\[
M_i:=1+\sup_{\substack{0\le t\le T_0\\|\xi|\le R}} \|\mathbf U_i^{*}(\cdot,\xi,\theta_t\omega)\|_{C_\eta^-}
\]
satisfies $\mathbb EM_i^p<\infty$. By $(\mathbf{H}2)$ and the definition of $M_i$,
the nonlinear term in the equation for $\hat u_i^*(\cdot,\xi,\theta_t\omega)$ in
\eqref{eq:limiting-backward-LP-system} has norm at most $C_RM_i e^{\eta r}$,
uniformly for $r\le0$, $t\in[0,T_0]$, and $|\xi|\le R$.
By \eqref{eq:positive-fractional-semigroup-gaussian} and \eqref{eq:regularity-deterministic-Vtheta}, the  
 convolution terms are bounded in $V_\vartheta$. Equation~\eqref{eq:semigroup-increment-gaussian} gives their $C^\rho$ bound on $[-T,0]$.
For the stochastic convolution, $r+t\in[-T,T_0]$. Applying the factorization method used to obtain \eqref{eq:OU-Vtheta-maximal-moment} on $[-T,T_0]$ gives the $V_\vartheta$
 supremum in $L^p$. The increment estimate \eqref{eq:regularity-stochastic-increment-bound} and the Garsia--Rodemich--Rumsey inequality 
 give the $C^\rho$ seminorm of $r\mapsto Z_{\mathrm{OU},i}(r,\theta_t\omega)$. This proves \eqref{eq:backward-joint-theta-regularity}.
Combining the estimates in $r$, $\xi$, and $t$ gives \eqref{eq:backward-joint-increment-bound}. The random bounds depend on the 
particle index only through $W_i$, hence the constants in \eqref{eq:backward-joint-increment-bound} may be chosen i.i.d. in $i$. 
Since $u_i^{*}=\bar u^{*}+\hat u_i^{*}$ and $\bar u^{*}$ is deterministic, the bounds in \eqref{eq:backward-joint-theta-regularity}
 and \eqref{eq:backward-joint-increment-bound} hold for both $u_i^{*}$ and $\mathbf U_i^{*}$.
For the derivative estimate, use the stronger weight $\nu$ from \eqref{eq:neighboring-weight-choice}. Repeating the preceding argument with $\nu$ in place of $\eta$ 
gives i.i.d. random variables $K_i$ with $\mathbb EK_i^p<\infty$ and
\[
\sup_{|\xi|\le R}
\|\mathbf U_i^{*}(\cdot,\xi,\theta_t\omega)-\mathbf U_i^{*}(\cdot,\xi,\theta_v\omega)\|_{C_\nu^-}
\le C_RK_i|t-v|^{\rho_0}.
\]
Since $\theta_t$ preserves $\mathbb P$, the joint law of $(\widetilde u_1^{*},D_\xi\widetilde u_1^{*})$ is unchanged. 
The derivative $D_\xi\bar u^*$ is deterministic, and the joint law
of $(\widetilde u_1^*,D_\xi\widetilde u_1^*)$ is invariant under $\theta_t$.
Subtracting the equations for $D_\xi\hat u_i^*$ in \eqref{eq:regularity-variational-resolvent} at $\theta_t\omega$ and $\theta_v\omega$,
we collect the terms linear in $D_\xi\hat u_i^*(\cdot,\xi,\theta_t\omega)-D_\xi\hat u_i^*(\cdot,\xi,\theta_v\omega)$.
By $(\mathbf{H}2)$, \eqref{eq:regularity-derivative-weighted-bound} and the preceding $C_\nu^-$ estimate, the remaining integrand is bounded by
\[
C_Re^{\nu r}\min\{1,K_i|t-v|^{\rho_0}e^{\nu r}\}
\le C_RK_i^a|t-v|^{a\rho_0}e^{(1+a)\nu r}.
\]
Hence the nonhomogeneous term in the difference equation is bounded by
$C_R(K_i^a+\mathbb EK_1^a)|t-v|^{a\rho_0}e^{(1+a)\nu r}$. Since $(1+a)\nu>\eta$, its $C_\eta^-$ norm is bounded by $C_R(K_i^a+\mathbb EK_1^a)|t-v|^{a\rho_0}$.
Estimating the difference equation in $C_\eta^-$ and absorbing the term with coefficient $q_s<1$, we obtain
\[
\sup_{|\xi|\le R}
\|D_\xi\mathbf U_i^{*}(\cdot,\xi,\theta_t\omega)-D_\xi\mathbf U_i^{*}(\cdot,\xi,\theta_v\omega)\|_{\mathcal L(H_{cu},C_\eta^-)}
\le C_R(K_i^a+\mathbb EK_1^a)|t-v|^{a\rho_0}.
\]
Taking moments proves \eqref{eq:derivative-theta-reduced-holder}. A common version on a countable dense subset of $[0,T_0]$, together 
with this H\"older estimate and the continuity in the center parameter from Lemma~\ref{lem:backward-regularity-gaussian}, yields joint continuity in $(t,\xi)$.
\end{proof}

\section{Auxiliary empirical field moment estimates}
\label{app:empirical-moments}
\par

\subsection*{Proof of Lemma~\ref{lem:driving-field-local-pmoment}}
\begin{proof}
Let $\rho$ satisfy \eqref{eq:rho-gaussian}, and set
\[\Theta_{T,R,M}:=[-T,0]\times B_R^{cu}\times K_M, \qquad d((r,\xi,z),(s,\zeta,y)):=|r-s|^\rho+|\xi-\zeta|+|z-y|,\]
and, for $\theta=(r,\xi,z)\in\Theta_{T,R,M}$, let
\[F_j^K(\theta):=\nabla K(u_j^{*}(r,\xi)-z), \qquad d_{N,K}(\theta,\theta') :=\left[\frac1N\sum_{j=1}^N |F_j^K(\theta)-F_j^K(\theta')|^2\right]^{1/2}.\]
By \eqref{eq:gaussian-LN-bound} and \eqref{eq:gaussian-Hj-moment}, $d_{N,K}\le L_N^K d$ and $\sup_N\mathbb E(L_N^K)^p<\infty$. Hence, on $\{L_N^K>0\}$,
\[
N(\Theta_{T,R,M},d_{N,K},u)
\le
N\!\left(\Theta_{T,R,M},d,\frac{u}{L_N^K}\right).
\]
Fix $\theta_0\in\Theta_{T,R,M}$. By $(\mathbf{H}2)$ and \eqref{eq:backward-Vtheta-moment}, $F_1^K(\theta_0)$ has finite $p$th moment. 
The Rosenthal inequality for Hilbert space valued sums \cite[Theorem~6.20]{LedouxTalagrand1991} gives
\[
\sup_{N\ge1}\mathbb E|\mathbb G_K^N(\theta_0)|^p<\infty.
\]
Let $(\varepsilon_j)_{j=1}^N$ be independent Rademacher variables, independent of $(F_j^K)_{j=1}^N$, and let $\mathbb E_\varepsilon$ denote expectation with respect to them.
 Conditionally on $(F_j^K)_{j=1}^N$, the entropy estimate and the moment comparison for Rademacher series \cite[Theorems~4.7 and~11.2, and Remark~11.5]{LedouxTalagrand1991} give
\[
\begin{aligned}
 &\left(\mathbb E_\varepsilon\sup_{\theta\in\Theta_{T,R,M}}\left| N^{-1/2}\sum_{j=1}^N\varepsilon_j\bigl[F_j^K(\theta)-F_j^K(\theta_0)\bigr]\right|^p\right)^{1/p}\\
 &\qquad\le C_pL_N^K\int_0^D\bigl[\log N(\Theta_{T,R,M},d,u)\bigr]^{1/2}\,du,
\end{aligned}
\]
where $D=\operatorname{diam}(\Theta_{T,R,M},d)$.
The entropy integral is finite by \eqref{eq:KM-entropy-integral} and the logarithmic entropy bound for $[-T,0]\times B_R^{cu}$ with respect to the metric $|r-s|^\rho+|\xi-\zeta|$.
Taking $p$th powers and expectations, using
$\sup_N\mathbb E(L_N^K)^p<\infty$, and applying symmetrization
\cite[Lemma~6.3]{LedouxTalagrand1991} together with the bound at
$\theta_0$ yield \eqref{eq:GK-local-pmoment}.
For $\mathbb G_U^N$, apply the same estimates to $\nabla U(u_j^{*}(r,\xi))$, with the parameter $z$ removed. The bounds \eqref{eq:backward-time-holder}
and \eqref{eq:backward-xi-Lipschitz}, together with the logarithmic entropy bound for $[-T,0]\times B_R^{cu}$ with respect to the metric
$|r-s|^\rho+|\xi-\zeta|$, give \eqref{eq:GU-local-pmoment}.
Choose $1/p<\rho_0<\rho$.
For the fields taken at $\theta_t\omega$, use the parameter sets $[0,T_0]\times[-T,0]\times B_R^{cu}$ and $[0,T_0]\times[-T,0]\times B_R^{cu}\times K_M$. 
By \eqref{eq:backward-joint-increment-bound}, the Lipschitz constants have uniformly bounded $p$ moments for the metric $|t-v|^{\rho_0}+|r-s|^\rho+|\xi-\zeta|+|z-y|.$
The covering numbers of $[0,T_0]$ for the metric $|t-v|^{\rho_0}$ have polynomial growth. Together with \eqref{eq:KM-entropy-integral}, this gives a
 finite entropy integral on each enlarged parameter set and proves the two estimates for $\mathbb G_{U,\theta}^N$ and $\mathbb G_{K,\theta}^N$ 
 in Lemma~\ref{lem:driving-field-local-pmoment}.
\end{proof}

\subsection*{Proof of Lemma~\ref{lem:particle-driving-forcing-pmoment}}
\begin{proof}
We use $1/p<\rho<\vartheta-1/p$ and, for the time parameter, $1/p<\rho_0<\rho$.
Fix $i$ and set $\mathcal G_i:=\sigma(W_i(s):s\in\mathbb R)$. We condition on $\mathcal G_i$. For $\theta=(r,\xi)\in[-T,0]\times B_R^{cu}$, write
\begin{equation}\label{eq:particle-GK-conditional-decomposition}
\mathbb G_K^N\bigl(r,\xi,u_i^{*}(r,\xi)\bigr) = N^{-1/2}\sum_{j\ne i}F_{j|i}(\theta) + N^{-1/2} \left[ \nabla K(0)
-\widetilde{\mathbb E}\nabla K\bigl(\widetilde u_1^{*}(r,\xi)-u_i^{*}(r,\xi)\bigr) \right],
\end{equation}
where
$F_{j|i}(\theta) := \nabla K\bigl(u_j^{*}(r,\xi)-u_i^{*}(r,\xi)\bigr) - 
\widetilde{\mathbb E}\nabla K\bigl(\widetilde u_1^{*}(r,\xi)-u_i^{*}(r,\xi)\bigr).$
Conditionally on $\mathcal G_i$, the variables $F_{j|i}(\theta)$, $j\ne i$, are independent and satisfy $\mathbb E[F_{j|i}(\theta)\mid\mathcal G_i]=0$. If $u_j^{*\prime}$ 
is an independent copy of $u_j^{*}$, conditional Jensen inequality, $(\mathbf{H}2)$, and \eqref{eq:backward-Vtheta-moment} give
\[
\begin{aligned}
\mathbb E\bigl[|F_{j|i}(\theta)|^p\mid\mathcal G_i\bigr]
&\le\mathbb E\mathbb E'\left[\left|\nabla K(u_j^{*}-u_i^{*})-\nabla K(u_j^{*\prime}-u_i^{*})\right|^p\,\middle|\,\mathcal G_i\right]\\
&\le C\,\mathbb E\mathbb E'|u_j^{*}-u_j^{*\prime}|^p \le C_{T,R,p}.
\end{aligned}
\]
 Moreover, \eqref{eq:backward-time-holder} and \eqref{eq:backward-xi-Lipschitz} imply
\begin{equation}\label{eq:conditional-particle-K-Lipschitz}
|F_{j|i}(r,\xi)-F_{j|i}(s,\zeta)| \le C(1+H_j+H_i) \bigl(|r-s|^\rho+|\xi-\zeta|\bigr),
\end{equation}
where $H_j$ are i.i.d. nonnegative random variables with $\mathbb EH_j^p<\infty$.
Define  $d_{N,i}\bigl((r,\xi),(s,\zeta)\bigr):=\left[\frac1N\sum_{j\ne i}|F_{j|i}(r,\xi)-F_{j|i}(s,\zeta)|^2\right]^{1/2}.$
By \eqref{eq:conditional-particle-K-Lipschitz},
\[
d_{N,i}\bigl((r,\xi),(s,\zeta)\bigr)\le C\left[\frac1N\sum_{j\ne i}(1+H_j+H_i)^2\right]^{1/2}\bigl(|r-s|^\rho+|\xi-\zeta|\bigr).
\]
The random multiplier on the right-hand side has a uniformly bounded $p$th moment. The fixed parameter estimate above and the conditional Rosenthal inequality \cite[Theorem~6.20]
{LedouxTalagrand1991} control the normalized sum at one fixed parameter point. Conditional symmetrization, the entropy estimate, and the moment comparison for Rademacher series 
\cite[Lemma~6.3, Theorems~4.7 and~11.2, and Remark~11.5]{LedouxTalagrand1991} then give, after taking expectations,
\[
\sup_{N\ge2}\max_{1\le i\le N}\mathbb E\sup_{\substack{-T\le r\le0\\|\xi|\le R}}\left|N^{-1/2}\sum_{j\ne i}F_{j|i}(r,\xi)\right|^p<\infty.
\] The $j=i$ term in \eqref{eq:particle-GK-conditional-decomposition} is bounded uniformly
 in $N$ and $i$ in $L^p$ after multiplication by $N^{-1/2}$, since $\nabla K(0)=0$ and
\[\left|\widetilde{\mathbb E}\nabla K\bigl(\widetilde u_1^{*}(r,\xi)-u_i^{*}(r,\xi)\bigr)\right| \le C\left( \mathbb E|u_1^{*}(r,\xi)| 
+ |u_i^{*}(r,\xi)| \right).\]
Taking the supremum in $(r,\xi)$ and using \eqref{eq:backward-Vtheta-moment} proves \eqref{eq:particle-GK-random-state-pmoment}.

Set $d_0((r,\xi),(s,\zeta)):=|r-s|^\rho+|\xi-\zeta|$.
Apply the local conditional entropy bound
\cite[Theorem~11.6 and Remark~11.5]{LedouxTalagrand1991}
to the increments of $N^{-1/2}\sum_{j\ne i}F_{j|i}$ with
$d_0((r,\xi),(s,\zeta))\le\delta$.
The moment comparison and integration over the sample give
the following bound for the $L^p$ norm of the supremum of their $H$ norms
\[C_{T,R,p} \int_0^{\delta} \sqrt{\log N([-T,0]\times B_R^{cu},d_0,u)}\,du,\]
which tends to zero as $\delta\downarrow0$. The $j=i$ term is controlled by \eqref{eq:backward-time-holder} and \eqref{eq:backward-xi-Lipschitz}. 
Hence \eqref{eq:particle-GK-Lp-equicontinuity} follows.
For \eqref{eq:theta-particle-GK-random-state-pmoment} and \eqref{eq:theta-particle-GK-Lp-equicontinuity}, condition again on $\mathcal G_i$ and take each trajectory 
at $\theta_t\omega$. The terms with $j\ne i$ remain conditionally independent and centered. By \eqref{eq:backward-joint-increment-bound}, their increments are 
bounded by
\[
C(\mathcal L_j+\mathcal L_i+\mathbb E\mathcal L_1)\bigl(|t-v|^{\rho_0}+|r-s|^\rho+|\xi-\zeta|\bigr),
\]
where $\sup_j\mathbb E\mathcal L_j^p<\infty$. Since $\theta_t$ preserves $\mathbb P$, the fixed parameter moment estimate above also holds 
for the trajectories $u_j^{*}(r,\xi,\theta_t\omega)$. Conditional symmetrization and the entropy estimate for Rademacher processes on the enlarged 
parameter set then give \eqref{eq:theta-particle-GK-random-state-pmoment} and \eqref{eq:theta-particle-GK-Lp-equicontinuity}.
\end{proof}

\end{document}